\documentclass [11pt]{article}
\usepackage{amsmath,amssymb}
\usepackage{cite}
\begin{document}

\title{ Cramer's Rule for Generalized Inverse Solutions of
Some Matrices Equations}
\author {\textbf{IVAN I. KYRCHEI}\\Pidstrygach Institute for Applied Problems of Mechanics and Mathematics\\ NAS of Ukraine, Ukraine,  kyrchei@lms.lviv.ua}

\date{}
 \maketitle


\begin{abstract}
By a generalized inverse of a given matrix, we mean a matrix  that exists for a larger class of matrices than the nonsingular matrices, that has some of the properties of the usual inverse, and that agrees with inverse when given matrix happens to be nonsingular. In theory, there are many different generalized inverses that exist. We shall consider the Moore Penrose, weighted Moore-Penrose,  Drazin and  weighted Drazin inverses.

New determinantal representations of these generalized inverse based
on their limit representations are introduced in this paper. Application of
this new method allows us to obtain analogues classical adjoint matrix.   Determinantal representations for projection matrices ${\rm {\bf A}}^{
+} {\rm {\bf A}}$, ${\rm {\bf A}} {\rm {\bf A}}^{ +}$, and ${\rm
{\bf A}}^{ D} {\rm {\bf A}}$ are presented as well. Using the obtained analogues
of the adjoint matrix, we get Cramer's rules for the least squares solution with the minimum norm
and for the Drazin inverse solution of singular linear systems.
The least squares solutions with the minimum norm and the Drazin inverse solutions of the matrix
equations ${\rm {\bf A}}{\rm {\bf X}} =
{\rm {\bf D}}$, ${\rm {\bf X}}{\rm {\bf B}} = {\rm {\bf D}}$ and ${\rm {\bf A}}{\rm {\bf X}}{\rm {\bf B}} ={\rm {\bf D}} $ are considered. Thus  ${\rm {\bf A}}$
${\rm {\bf B}}$ can be  singular matrices of appropriate size.  We  get analogs of the Cramer rule for the generalized inverse solutions of these matrix equations and using their for determinantal representations  of solutions of some differential matrix equations, ${\bf X}'+ {\bf A}{\bf X}={\bf B}$ and ${\bf X}'+{\bf X}{\bf A}={\bf B}$, where the matrix ${\bf A}$ is singular.


\end{abstract}

\section{Preface}
\newtheorem{corollary}{Corollary}[section]
\newtheorem{theorem}{Theorem}[section]
\newtheorem{lemma}{Lemma}[section]
\newtheorem{definition}{Definition}[section]
\newtheorem{remark}{Remark}[section]
\newcommand{\rank}{\mathop{\rm rank}\nolimits}
\newtheorem{proposition}{Proposition}[section]

It's well-known in linear algebra, an $n$-by-$n$ square matrix ${\bf A}$ is called invertible (also nonsingular or nondegenerate) if there exists an $n$-by-$n$ square matrix ${\bf X}$ such that
\[{\rm {\bf A}}{\bf X}={\bf X}{\bf A}={\bf I}_{n}.\]
If this is the case, then the matrix ${\bf X}$ is uniquely determined by ${\bf A}$ and is called the inverse of ${\bf A}$, denoted by ${\bf A}^{-1}$.

By a generalized inverse of a given matrix, we mean a matrix  that exists for a larger class of matrices than the nonsingular matrices, that has some of the properties of the usual inverse, and that agrees with inverse when given matrix happens to be nonsingular.

For any matrix ${\rm {\bf A}}\in {\mathbb
C}^{m\times n} $ consider the following equations in ${\bf X}$:
 \begin{equation}\label{eq1:MP}  {\rm {\bf A}}{\bf X}
{\rm {\bf A}} = {\rm {\bf A}}; \end{equation}
                                \begin{equation}\label{eq2:MP}  {\bf X} {\rm {\bf
A}}{\bf X}  = {\bf X};
                               \end{equation}
\begin{equation}\label{eq3:MP}
                                  \left( {\rm {\bf A}}{\bf X} \right)^{ *}  = {\rm
{\bf A}}{\bf X}; \end{equation}
                               \begin{equation}\label{eq4:MP} \left( {{\bf X} {\rm {\bf A}}} \right)^{ *}  ={\bf X} {\rm {\bf A}}. \end{equation}
                                and if $m = n$, also
                                \begin{equation}\label{eq5:MP}
                                   {\rm {\bf A}}{\bf X}  = {\bf X}{\rm
{\bf A}}; \end{equation}
                                 \begin{equation}\label{eq6:MP}
                                   {\rm {\bf A}}^{k+1}{\rm
{\bf X}}={\rm {\bf
         A}}^{k}. \end{equation}
For a sequence $\mathcal{G}$ of $\{1,2,3,4,5\}$ the set of matrices obeying the equations
represented in $\mathcal{G}$ is denoted by ${\bf A}\{\mathcal{G}\}$.
 A matrix from ${\bf A}\{\mathcal{G}\}$ is called an $\mathcal{G}$-inverse of ${\bf A}$
 and denoted by ${\bf A}^{(\mathcal{G})}$.

Consider some principal cases.

 If ${\bf X}$ satisfies the all equations
(\ref{eq1:MP})-(\ref{eq4:MP}) is said to be \textbf{the Moore-Penrose
inverse} of ${\bf A}$ and denote ${\bf A}^{+}={\bf A}^{(1, 2, 3, 4)}$. The
Moore–Penrose inverse was independently described by E. H. Moore \cite{mo}
in 1920, Arne Bjerhammar \cite{bj} in 1951 and Roger Penrose \cite{pen} in
1955. R. Penrose introduced  the characteristic equations
(\ref{eq1:MP})-(\ref{eq4:MP}). If $\det {\bf A}\neq 0$, then ${\bf A}^{+}={\bf A}^{-1}$.

\textbf{The group inverse} ${\bf A}^{g}$ is the unique ${\bf A}^{(1, 2, 5)}$
inverse of ${\bf A}$, and exists if and only if $Ind\,{\bf A}= {\rm min}\{k: \rank{\bf A}^{k+1}=\rank{\bf A}^{k}\}=1$.

A matrix ${\bf X}={\bf A}^{D}$ is said to be \textbf{the Drazin inverse} of ${\bf A}$ if (\ref{eq6:MP})
(for some positive integer $k$), (\ref{eq2:MP}) and (\ref{eq5:MP}) are satisfied, where $k=Ind\,{\bf A}$. It is  named after Michael P. Drazin \cite{dr}. In particular, when $Ind{\kern 1pt} {\rm {\bf A}}=1$,
then the matrix ${\rm {\bf X}}$  is
the group inverse, ${\rm {\bf X}}={\rm {\bf
A}}^{g }$.
If $Ind{\kern 1pt} {\rm {\bf A}}=0$, then ${\rm {\bf A}}$ is
nonsingular, and ${\rm {\bf A}}^{D}\equiv {\rm {\bf A}}^{-1}$.

Let Hermitian positive definite matrices ${\bf M}$ and ${\bf N}$ of order $m$ and $n$, respectively, be given. For
any  matrix ${\bf A}\in {\mathbb C}^{m\times n}$, \textbf{the weighted Moore-Penrose inverse} of ${\bf A}$ is the unique
solution ${\bf X}={\bf A}^{+}_{M,N}$ of the matrix equations (\ref{eq1:MP}) and (\ref{eq2:MP}) and the following
equations in ${\bf X}$ \cite{pr}:
\[(3M)\,\,({\bf M}{\rm {\bf A}}{\bf X})^{*}  = {\bf M}{\rm
{\bf A}}{\bf X};\,\,(4N)\,\,({\bf N}{\bf X}{\rm {\bf A}})^{*}  = {\bf N}{\bf X}{\rm
{\bf A}}.\]
In particular, when ${\bf M}={\bf I}_m$ and ${\bf N}={\bf I}_n$, the matrix ${\bf X}$ satisfying the  equations (\ref{eq1:MP}), (\ref{eq2:MP}), (3M), (4N) is the  Moore-Penrose inverse ${\bf A}^{+}$.

The weighted Drazin inverse is being  considered as well.

To determine the inverse and to give its analytic solution, we calculate a matrix of cofactors,
 known as an adjugate matrix or a classical adjoint matrix. The classical adjoint of ${\bf A}$, denote $Adj[ {\bf A}]$, is the transpose of the cofactor matrix, then ${\bf A}^{-1}=\frac{Adj[ {\bf A}]}{\| {\bf A}  \|}$. Representation an inverse matrix by its classical adjoint matrix also plays a key role for Cramer's rule of  systems of linear equations or matrices equations.

 Obviously, the important question is the following:
  what are the analogues  for the adjoint matrix of generalized inverses and, consequently, for Cramer's rule of generalized inverse solutions of matrix equations?

This is the main goal of the paper.
In the paper we shall adopt the following notation. Let ${\mathbb
C}^{m\times n} $ be the set of $m$ by $n$ matrices with complex
entries, ${\mathbb C}^{m\times n}_{r} $ be a subset of ${\mathbb
C}^{m\times n} $  in which any matrix has rank $r$, ${\rm \bf
I}_{m}$ be the identity matrix of order $m$, and $\|.\|$ be the
Frobenius norm of a matrix.

Denote by ${\rm {\bf a}}_{.j} $ and ${\rm {\bf a}}_{i.} $ the
$j$th column  and the $i$th row of ${\rm {\bf A}}\in {\mathbb
C}^{m\times n} $, respectively. Then   ${\rm {\bf a}}^{\ast}_{.j}
$ and ${\rm {\bf a}}^{\ast}_{i.} $ denote the $j$th column  and
the $i$th row of a conjugate and transpose matrix ${\rm {\bf A}}^{\ast}$
as well. Let ${\rm {\bf A}}_{.j} \left( {{\rm {\bf b}}} \right)$
denote the matrix obtained from ${\rm {\bf A}}$ by replacing its
$j$th column with the  vector ${\rm {\bf b}}$, and by ${\rm {\bf
A}}_{i.} \left( {{\rm {\bf b}}} \right)$ denote the matrix
obtained from ${\rm {\bf A}}$ by replacing its $i$th row with
${\rm {\bf b}}$.

Let $\alpha : = \left\{ {\alpha _{1} ,\ldots ,\alpha _{k}} \right\}
\subseteq {\left\{ {1,\ldots ,m} \right\}}$ and $\beta : = \left\{ {\beta
_{1} ,\ldots ,\beta _{k}}  \right\} \subseteq {\left\{ {1,\ldots ,n}
\right\}}$ be subsets of the order $1 \le k \le \min {\left\{ {m,n}
\right\}}$. Then ${\left| {{\rm {\bf A}}_{\beta} ^{\alpha} } \right|}$
denotes the minor of ${\rm {\bf A}}$ determined by the rows indexed by
$\alpha$ and the columns indexed by $\beta$. Clearly, ${\left| {{\rm {\bf
A}}_{\alpha} ^{\alpha} } \right|}$ denotes a principal minor determined by
the rows and columns indexed by $\alpha$.  The cofactor of $a_{ij} $ in
${\rm {\bf A}}\in {\mathbb C}^{n\times n} $ is denoted by
${\frac{{\partial }}{{\partial a_{ij}} }}{\left| {{\rm {\bf A}}}
\right|}$.

For $1 \leq k\leq n$,  $\textsl{L}_{ k, n}: = {\left\{
{\,\alpha :\alpha = \left( {\alpha _{1} ,\ldots ,\alpha _{k}}
\right),\,{\kern 1pt} 1 \le \alpha _{1} \le \ldots \le \alpha _{k} \le n}
\right\}}$ denotes the collection of strictly increasing sequences of $k$ integers
chosen from $\left\{ {1,\ldots ,n} \right\}$. Let $N_{k} : =
\textsl{L}_{k, m} \times \textsl{L}_{ k, n} $. For fixed $\alpha \in
\textsl{L}_{ p, m}$, $\beta \in \textsl{L}_{ p, n}$, $1\leq p\leq k$, let
\[\begin{array}{c}
   I_{k,\,m} \left( {\alpha} \right): = {\left\{ {I:\,I \in
\textsl{L}_{ k,\, m} ,I \supseteq \alpha} \right\}},\\
  J_{k,\,n}
\left( {\beta} \right): = {\left\{ {J:\,J \in \textsl{L}_{ k,\, n}
,J \supseteq \beta} \right\}},\\
N_{k} \left( {\alpha ,\beta} \right): = I_{k,\,m} \left( {\alpha}
\right)\times J_{k,\,n} \left( {\beta}  \right)
\end{array}
\]
For case $i \in \alpha $ and $j \in \beta$, we denote
\[\begin{array}{c}
I_{k,m} {\left\{ {i} \right\}}: = {\left\{ {\alpha:\,\alpha \in L_{k, m},
i \in \alpha} \right\}},
 J_{k,\,n} {\left\{ {j}
\right\}}: =  {\left\{ {\beta:\,\beta\in  L_{k, n}, j \in  \beta}
\right\}},\\
N_{k} {\left\{ {i,j} \right\}}: = I_{k,\, m} {\left\{ {i} \right\}}\times
J_{k,\, n} {\left\{ {j} \right\}}.
\end{array}
\]
 The paper is
organized as follows.
In Section 2    determinantal representations by analogues of the classical adjoint matrix for the Moore Penrose, weighted Moore-Penrose,  Drazin and  weighted Drazin inverses are obtained.

In Section 3 we show that the obtained analogues of the adjoint matrix for the generalized inverse matrices enable us to obtain
 natural analogues of Cramer's rule for generalized inverse solutions of systems of  linear equations and demonstrate it in two examples.

In Section 4, we  obtain analogs of the
Cramer rule for generalized inverse solutions of the matrix equations, ${\bf A}{\bf X}={\bf B}$, ${\bf X}{\bf A}={\bf B}$ and ${\bf A}{\bf X}{\bf B}={\bf D}$,
 namely for the minimum norm least squares solutions  and Drazin inverse solutions.
We show  numerical examples to illustrate the main results as well.

In Section 5, we  use the  determinantal representations of Drazin inverse solution to solutions of the following differential matrix equations, ${\bf X}'+ {\bf A}{\bf X}={\bf B}$ and ${\bf X}'+{\bf X}{\bf A}={\bf B}$, where  ${\bf A}$ is singular. It is demonstrated in the example.

Facts set forth in Sections 2 and 3 were partly published in \cite{ky1}, in
Section 4  were published in \cite{ky2,ky3} and in Sections 5
were published in \cite{ky3}.

Note that  we obtained some of the submitted results   for matrices over the quaternion skew field within the framework of the theory of the column and row determinants (\cite{ky4,ky41,ky5,ky6,ky7,ky8,ky9,ky10}).

\section{ Analogues of the classical adjoint matrix for the generalized inverse matrices}
For determinantal representations of the generalized inverse matrices as analogues
of the classical adjoint matrix, we apply the method,
 which consists  on the limit representation of the generalized inverse matrices,
  lemmas on rank of some matrices and on characteristic polynomial.
  We   used this method at first in \cite{ky1} and then in \cite{ky3}.
  Liu et al. in  \cite{liu1} deduce the new determinantal
representations of the outer inverse ${\bf A}_{T,S}^{(2)}$  based on these principles as well.
In  this chapter we obtain  detailed determinantal representations by analogues of the classical adjoint matrix for the Moore Penrose, weighted Moore-Penrose,  Drazin and  weighted Drazin inverses.

\subsection{ Analogues of the classical adjoint matrix for the Moore
- Penrose inverse}

 Determinantal representation of the
Moore - Penrose inverse was studied in \cite{mo},\cite{ba,ben1,ga,st1}.
The main result  consists in the following theorem.
\begin{theorem}\label{theor:st} The Moore - Penrose inverse
 ${\rm {\bf A}}^{+}=(a_{ij}^{+})\in {\mathbb C}^{n\times m}$
of ${\rm {\bf A}}\in {\mathbb C}_{r}^{m\times n} $ has the
following determinantal representation
\[
a_{ij}^{+} = {\frac{{{\sum\limits_{\left( {\alpha,\,\beta} \right)
\in N_{r}  {\left\{ {j,\,i} \right\}}} {{\left| {\left( {{\rm {\bf
A}}^{\ast}} \right)_{\alpha} ^{\beta} } \right|}{\frac{{\partial}
}{{\partial a_{j\,i} }}}{\left| {{\rm {\bf A}}_{\beta} ^{\alpha} }
\right|}} }}}{{{\sum\limits_{\left( {\gamma,\,\delta}  \right) \in
N_{r} } {{\left| {\left( {{\rm {\bf A}}^{\ast}} \right)_{\gamma}
^{\delta} } \right|}\,{\left| {{\rm {\bf A}}_{\delta} ^{\gamma} }
\right|}}} }}}, \,\,1 \leq i, j \leq n.
\]
 \end{theorem}
This determinantal representation  of the Moore - Penrose inverse is
based on corresponding full-rank representation \cite{ba}: if ${\bf A}={\bf P}{\bf Q}$, where ${\rm {\bf P}}\in {\mathbb C}_{r}^{m\times r} $ and ${\rm {\bf Q}}\in {\mathbb C}_{r}^{r\times n} $, then
\[
{\bf A}^{+}={\bf Q}^{*}({\bf P}^{*}{\bf A}{\bf Q}^{*})^{-1}{\bf P}^{*}.
\]

For a better understanding of the structure of the Moore - Penrose inverse we consider it  by singular value decomposition of ${\bf A}$. Let
\[\begin{array}{cc}
    {\bf A}{\bf A}^{*}{\bf u}_{i}=\sigma^{2}_{i}{\bf u}_{i}, & i=\overline{1, m} \\
    {\bf A}^{*}{\bf A}{\bf v}_{i}=\sigma^{2}_{i}{\bf v}_{i}, & i=\overline{1, n},
  \end{array}
\]
\[\sigma_{1}\leq\sigma_{2}\leq...\sigma_{r}>0=\sigma_{r+1}=\sigma_{r+2}=...\]
and the singular value decomposition (SVD) of ${\bf A}$ is ${\bf A}= {\bf U}{\bf \Sigma}{\bf V}^{*}$, where
\[\begin{array}{cc}
    {\bf U}=[ {\bf u}_{1}\,{\bf u}_{2}...{\bf u}_{m}]\in {\mathbb C}^{m\times m}, & {\bf U}^{*}{\bf U}={\bf I}_m, \\
    {\bf V}=[ {\bf v}_{1}\,{\bf v}_{2}...{\bf v}_{n}]\in {\mathbb C}^{n\times n}, & {\bf V}^{*}{\bf V}={\bf I}_n,
  \end{array}
\]
\[{\bf \Sigma}={\rm diag}(\sigma_{1},\sigma_{2},...,\sigma_{r})\in {\mathbb C}^{m\times n}.\] Then \cite{pen},
${\bf A}^{+}= {\bf V}{\bf \Sigma}^{+}{\bf U}^{*},
$
where ${\bf \Sigma}^{+}={\rm diag}(\sigma^{-1}_{1},\sigma^{-1}_{2},...,\sigma^{-1}_{r})$.

 We need the following
 limit representation of the Moore-Penrose inverse.
 \begin{lemma}\cite{ba1}\label{kyrc2}
If ${\rm {\bf A}}\in {\mathbb C}^{m\times n} $, then
\[{\rm {\bf
A}}^{ +} = {\mathop {\lim }\limits_{\lambda \to 0}} {\rm {\bf
A}}^{ *} \left( {{\rm {\bf A}}{\rm {\bf A}}^{ *}  + \lambda {\rm
{\bf I}}} \right)^{ - 1} = {\mathop {\lim }\limits_{\lambda\to 0}}
\left( {{\rm {\bf A}}^{ *} {\rm {\bf A}} + \lambda {\rm {\bf I}}}
\right)^{ - 1}{\rm {\bf A}}^{ *}, \]
 where $\lambda \in {\mathbb
R} _{ +}  $, and ${\mathbb R} _{ +} $ is the set of  positive real
numbers.
\end{lemma}

 \begin{corollary}\cite{ho}\label{cor:A+A-1} If
${\rm {\bf A}}\in {\mathbb C}^{m\times n} $, then the following
statements are true.
 \begin{itemize}
\item [ i)] If $\rm{rank}\,{\rm {\bf A}} = n$, then ${\rm {\bf A}}^{ +}
= \left( {{\rm {\bf A}}^{ *} {\rm {\bf A}}} \right)^{ - 1}{\rm
{\bf A}}^{ * }$ .
\item [ ii)] If $\rm{rank}\,{\rm {\bf A}} =
m$, then ${\rm {\bf A}}^{ +}  = {\rm {\bf A}}^{ * }\left( {{\rm
{\bf A}}{\rm {\bf A}}^{ *} } \right)^{ - 1}.$
\item [ iii)] If $\rm{rank}\,{\rm {\bf A}} = n = m$, then ${\rm {\bf
A}}^{ +}  = {\rm {\bf A}}^{ - 1}$ .
\end{itemize}
\end{corollary}

We need the following well-known theorem about the characteristic polynomial and lemmas on rank of some matrices.
\begin{theorem} \cite{la}\label{th:char_pol} Let $d_{r}$ be the sum of  principal minors
of order $r$ of  ${\rm {\bf A}}\in {\mathbb C}^{n\times n}$. Then
its characteristic polynomial $ p_{{\rm {\bf A}}}\left( {t}
\right)$ can be expressed as  $ p_{{\rm {\bf A}}}\left( {t}
\right) = \det \left( { t{\rm {\bf I}} - {\rm {\bf A}}} \right) =
t^{n} - d_{1} t^{n - 1} + d_{2} t^{n - 2} - \ldots + \left( { - 1}
\right)^{n}d_{n}. $
\end{theorem}

\begin{lemma} \label{kyrc5} If ${\rm {\bf A}}\in {\mathbb C}^{m\times n}_{r}
$, then $
 \rank\,\left( {{\rm {\bf A}}^{ *} {\rm {\bf A}}}
\right)_{.\,i} \left( {{\rm {\bf a}}_{.j}^{ *} }  \right) \le r.
$
\end{lemma}
 {\textit{Proof}}.
Let ${\rm {\bf P}}_{i\,k} \left( {-a_{j\,k}}  \right)\in {\mathbb
C}^{n\times n} $, $(k \ne i )$, be the matrix with $-a_{j\,k} $ in
the $(i, k)$ entry, 1 in all diagonal entries, and 0 in others. It
is the  matrix of an elementary transformation. It follows that
\[
\left( {{\rm {\bf A}}^{ *} {\rm {\bf A}}} \right)_{.\,i} \left(
{{\rm {\bf a}}_{.\,j}^{ *} }  \right) \cdot {\prod\limits_{k \ne
i} {{\rm {\bf P}}_{i\,k} \left( {-a_{j\,k}}  \right) = {\mathop
{\left( {{\begin{array}{*{20}c}
 {{\sum\limits_{k \ne j} {a_{1k}^{ *}  a_{k1}} } } \hfill & {\ldots}  \hfill
& {a_{1j}^{ *} }  \hfill & {\ldots}  \hfill & {{\sum\limits_{k \ne
j } {a_{1k}^{ *}  a_{kn}}}} \hfill \\
 {\ldots}  \hfill & {\ldots}  \hfill & {\ldots}  \hfill & {\ldots}  \hfill &
{\ldots}  \hfill \\
 {{\sum\limits_{k \ne j} {a_{nk}^{ *}  a_{k1}} } } \hfill & {\ldots}  \hfill
& {a_{nj}^{ *} }  \hfill & {\ldots}  \hfill & {{\sum\limits_{k \ne
j } {a_{nk}^{ *}  a_{kn}}}} \hfill \\
\end{array}} }
\right)}\limits_{i-th}}}}.
\]
 The obtained above matrix  has the following factorization.
\[
{\mathop {\left( {{\begin{array}{*{20}c}
 {{\sum\limits_{k \ne j} {a_{1k}^{ *}  a_{k1}} } } \hfill & {\ldots}  \hfill
& {a_{1j}^{ *} }  \hfill & {\ldots}  \hfill & {{\sum\limits_{k \ne
j } {a_{1k}^{ *}  a_{kn}} } } \hfill \\
 {\ldots}  \hfill & {\ldots}  \hfill & {\ldots}  \hfill & {\ldots}  \hfill &
{\ldots}  \hfill \\
 {{\sum\limits_{k \ne j} {a_{nk}^{ *}  a_{k1}} } } \hfill & {\ldots}  \hfill
& {a_{nj}^{ *} }  \hfill & {\ldots}  \hfill & {{\sum\limits_{k \ne
j } {a_{nk}^{ *}  a_{kn}} } } \hfill \\
\end{array}} } \right)}\limits_{i-th}}  =
\]
\[
 = \left( {{\begin{array}{*{20}c}
 {a_{11}^{ *} }  \hfill & {a_{12}^{ *} }  \hfill & {\ldots}  \hfill &
{a_{1m}^{ *} }  \hfill \\
 {a_{21}^{ *} }  \hfill & {a_{22}^{ *} }  \hfill & {\ldots}  \hfill &
{a_{2m}^{ *} }  \hfill \\
 {\ldots}  \hfill & {\ldots}  \hfill & {\ldots}  \hfill & {\ldots}  \hfill
\\
 {a_{n1}^{ *} }  \hfill & {a_{n2}^{ *} }  \hfill & {\ldots}  \hfill &
{a_{nm}^{ *} }  \hfill \\
\end{array}} } \right){\mathop {\left( {{\begin{array}{*{20}c}
 {a_{11}}  \hfill & {\ldots}  \hfill & {0} \hfill & {\ldots}  \hfill &
{a_{n1}}  \hfill \\
 {\ldots}  \hfill & {\ldots}  \hfill & {\ldots}  \hfill & {\ldots}  \hfill &
{\ldots}  \hfill \\
 {0} \hfill & {\ldots}  \hfill & {1} \hfill & {\ldots}  \hfill & {0} \hfill
\\
 {\ldots}  \hfill & {\ldots}  \hfill & {\ldots}  \hfill & {\ldots}  \hfill &
{\ldots}  \hfill \\
 {a_{m1}}  \hfill & {\ldots}  \hfill & {0} \hfill & {\ldots}  \hfill &
{a_{mn}}  \hfill \\
\end{array}} } \right)}\limits_{i-th}} j-th.
\]
 Denote by ${\rm {\bf \tilde {A}}}: = {\mathop
{\left( {{\begin{array}{*{20}c}
 {a_{11}}  \hfill & {\ldots}  \hfill & {0} \hfill & {\ldots}  \hfill &
{a_{1n}}  \hfill \\
 {\ldots}  \hfill & {\ldots}  \hfill & {\ldots}  \hfill & {\ldots}  \hfill &
{\ldots}  \hfill \\
 {0} \hfill & {\ldots}  \hfill & {1} \hfill & {\ldots}  \hfill & {0} \hfill
\\
 {\ldots}  \hfill & {\ldots}  \hfill & {\ldots}  \hfill & {\ldots}  \hfill &
{\ldots}  \hfill \\
 {a_{m1}}  \hfill & {\ldots}  \hfill & {0} \hfill & {\ldots}  \hfill &
{a_{mn}}  \hfill \\
\end{array}} } \right)}\limits_{i-th}} j-th$. The matrix
${\rm {\bf \tilde {A}}}$ is obtained from ${\rm {\bf A}}$ by
replacing all entries of the $j$th row  and of the $i$th column
with zeroes except that the $(j, i)$ entry equals 1. Elementary
transformations of a matrix do not change its rank. It follows
that $\rank\left( {{\rm {\bf A}}^{ *} {\rm {\bf A}}}
\right)_{.\,i} \left( {{\rm {\bf a}}_{.j}^{ *} } \right) \le \min
{\left\{ {\rank{\rm {\bf A}}^{ * },\rank{\rm {\bf \tilde {A}}}}
\right\}}$. Since $\rank{\rm {\bf \tilde {A}}} \ge \rank\,{\rm
{\bf A}} = \rank{\rm {\bf A}}^{ *} $ and $\rank{\rm {\bf A}}^{
* }{\rm {\bf A}} = \rank{\rm {\bf A}}$ the proof is completed.$\blacksquare$

The following lemma can be proved in the same way.
\begin{lemma} If ${\rm {\bf A}}\in {\mathbb C}^{m\times n}_{r} $,
then $ \rank\left( {{\rm {\bf A}}{\rm {\bf A}}^{ *} }
\right)_{i\,.} \left( {{\rm {\bf a}}_{j\,.}^{ *} }  \right) \le
r.$
\end{lemma}
Analogues of the characteristic polynomial are considered in the following
two lemmas.
\begin{lemma}\label{lem:char_mpdet}
If ${\bf A}\in
{\mathbb{C}}^{m\times n}$ and $\lambda \in \mathbb{R}$, then

\begin{equation}
\label{eq:char_mpdet}\det \left(  \left(\lambda {\rm {\bf I}}_n + {\bf A}^{
*}{\bf A}  \right)_{.\,i} \left( { {{\bf a}}_{.j}^{ *} } \right)\right) =
c_{1}^{\left( {ij} \right)} \lambda^{n - 1} + c_{2}^{\left( {ij} \right)}
\lambda^{n - 2} + \ldots + c_{n}^{\left( {ij} \right)},
\end{equation}

\noindent where $c_{n}^{\left( {ij} \right)} = \left| \left( {{\bf A}^{
*}{\bf A} } \right)_{.\,i} \left( { {{\bf a}}_{.j}^{ *} } \right)\right|$
and $c_{s}^{\left( {ij} \right)} = {\sum\limits_{\beta \in J_{s,\,n}
{\left\{ {i} \right\}}} {  \left|\left( {\left( {{\bf A}^{ *}{\bf A} }
\right)_{.\,i} \left( { {{\bf a}}_{.j}^{ *} }  \right)}\right)_{\beta} ^{\beta} \right|
 } }$ for all $s = \overline {1, n- 1} $, $i = \overline
{1,n}$, and $j = \overline {1,m}$.
\end{lemma}
{\textit{Proof.}} Denote ${\bf A}^{ *}{\bf A}={\bf V}=\left(v_{ij} \right)\in
{\mathbb{C}}^{n\times n}$. Consider  $\left(\lambda {\rm {\bf I}}_n +
{\bf V}  \right)_{.\,i} \left( { {{\bf v}}_{.i} } \right)
\in {\rm {\mathbb{C}}}^{n\times n}$.  Taking into account Theorem
\ref{th:char_pol} we obtain
\begin{equation}
\label{eq:char_mpdet1}  \left|\left(\lambda {\rm {\bf I}}_n +{\bf V}  \right)_{.\,i} \left( { {{\bf v}}_{.i} } \right)\right| =
d_{1} \lambda^{n - 1} + d_{2} \lambda^{n - 2} + \ldots + d_{n},
\end{equation}
where $d_{s} = {\sum\limits_{\beta \in J_{s,\,n} {\left\{ {i} \right\}}}
| ( {{\rm {\bf V}} })_{\beta} ^{\beta}   | } $ is the
sum of all principal minors of order $s$ that contain the $i$-th column
for all $s = \overline {1, n - 1} $ and $d_{n} = \det  {{\rm {\bf
V}} }$. Since ${\rm {\bf v}}_{. {\kern 1pt} i} =
 {\sum\limits_{l} {{\rm {\bf a}}_{.\,l}^{*}  a_{li}}
}$, where ${\rm {\bf a}}_{.{\kern 1pt}  l}^{*}  $ is the $l$th
column-vector  of ${\rm {\bf {A}}}^{*}$ for all
$l=\overline{1,n}$, then we have on the one hand
\begin{equation}
\label{eq:char_mpdet2}
\begin{array}{c}
  \left| \left( {\lambda{\rm {\bf I}} + {\rm {\bf V}} }
\right)_{.{\kern 1pt} i} \left( {\rm {\bf v}}_{. {\kern 1pt} i}
\right)\right|
   = \sum\limits_{l} \left|  \left( {\lambda{\rm {\bf I}} + {\rm {\bf
   V}}} \right)_{.{\kern 1pt} l} \left( {{\rm {\bf {a}}}_{.\,l}^{*}
{a}_{li}} \right) \right|= \\ {\sum\limits_{l} \left| {\left(
{\lambda{\rm {\bf I}} + {\rm {\bf V}} } \right)_{.{\kern 1pt} i}
\left( {{\rm {\bf a}}_{.{\kern 1pt} {\kern 1pt} l}^{*} }  \right)
} \right|} \cdot {a}_{li}
\end{array}
\end{equation}
 Having changed the order of summation,  we obtain on the other hand
for all $s = \overline {1, n - 1} $
\begin{equation}
\label{eq:char_mpdet3}
\begin{array}{c}
  d_{s} = {\sum\limits_{\beta \in J_{s,\,n} {\left\{ {i} \right\}}}
{\left| \left( {{\rm {\bf V}} } \right)_{\beta}
^{\beta}  \right|} }  =
  \sum\limits_{\beta \in J_{s,\,n}
{\left\{ {i} \right\}}} \sum\limits_{l} \left|\left({ { {{\rm {\bf
V}} }_{.\, i} \left( {{\rm {\bf a}}_{.\,
l}^{*}a_{l\,i}} \right)} }\right)_{\beta} ^{\beta}
\right|=\\
  \sum\limits_{l}  \sum\limits_{\beta \in J_{s,\,n} {\left\{ {i}
\right\}}} \left|\left({ {{\rm {\bf V}} }
_{.{\kern 1pt} i} \left( {{\rm {\bf a}}_{.{\kern 1pt}
l}^{*} }  \right)} \right)_{\beta} ^{\beta}  \right|  \cdot
a_{l{\kern 1pt} i}.
\end{array}
\end{equation}
By substituting (\ref{eq:char_mpdet2}) and (\ref{eq:char_mpdet3}) in
(\ref{eq:char_mpdet1}), and equating factors at $a _ {l \, i} $ when
$l = j $, we obtain the equality (\ref{eq:char_mpdet}). $\blacksquare$

 By analogy can be proved
the following lemma.
\begin{lemma}
If ${\bf A}\in {\mathbb{C}}^{m\times
n}$  and $\lambda \in \mathbb{R}$,
then
\[ \det\left(  {( \lambda{\rm {\bf I}}_{m} + {\bf A} {\bf A}^{*} )_{j\,.\,} ({\rm {\bf a}}_{i.}^{*}
)}\right)
   = r_{1}^{\left( {ij} \right)} \lambda^{m
- 1} +r_{2}^{\left( {ij} \right)} \lambda^{m - 2} + \ldots + r_{m}^{\left(
{ij} \right)},
\]

\noindent where  $r_{m}^{\left( {ij} \right)} = \left| {({\bf A} {\bf A}^{*} )_{j\,.\,} ({\rm {\bf a}}_{i.\,}^{*} )}\right|$ and
$r_{s}^{\left( {ij} \right)} = \sum\limits_{\alpha \in I_{s,m} {\left\{
{j} \right\}}} \left| \left( {({\bf A} {\bf A}^{*} )_{j\,.\,} ({\rm {\bf
a}}_{i.\,}^{*} )}\right)_{\alpha} ^{\alpha} \right| $ for all
$s = \overline {1,n - 1} $,  $i = \overline {1,m}$, and $j = \overline
{1,n}$.
\end{lemma}
The following theorem and remarks introduce the determinantal representations of the Moore-Penrose by analogs of the classical adjoint matrix.

\begin{theorem}\label{theor:det_repr_MP}
If ${\rm {\bf A}} \in {\rm {\mathbb{C}}}_{r}^{m\times n} $ and
$r<min\{m,n\}$, then the Moore-Penrose inverse  ${\rm {\bf A}}^{
+} = \left( {a_{ij}^{ +} } \right) \in {\rm
{\mathbb{C}}}_{}^{n\times m} $ possess the following determinantal
representations:
\begin{equation}
\label{eq:det_repr_A*A}
 a_{ij}^{ +}  = {\frac{{{\sum\limits_{\beta
\in J_{r,\,n} {\left\{ {i} \right\}}} {{\left| \left( {\left(
{{\rm {\bf A}}^{ *} {\rm {\bf A}}} \right)_{\,. \,i} \left( {{\rm
{\bf a}}_{.j}^{ *} }  \right)} \right) {\kern 1pt}
_{\beta} ^{\beta} \right|}} }  }}{{{\sum\limits_{\beta \in
J_{r,\,\,n}} {{\left| {\left( {{\rm {\bf A}}^{ *} {\rm {\bf A}}}
\right){\kern 1pt} _{\beta} ^{\beta} }  \right|}}} }}},
\end{equation}
or
\begin{equation}
\label{eq:det_repr_AA*} a_{ij}^{ +}  =
{\frac{{{\sum\limits_{\alpha \in I_{r,m} {\left\{ {j} \right\}}}
{{\left| \left( {({\rm {\bf A}}{\rm {\bf A}}^{ *} )_{j\,.\,} ({\rm
{\bf a}}_{i.\,}^{ *} )} \right)\,_{\alpha} ^{\alpha}\right|} }
}}}{{{\sum\limits_{\alpha \in I_{r,\,m}}  {{\left| {\left( {{\rm
{\bf A}}{\rm {\bf A}}^{ *} } \right){\kern 1pt} _{\alpha}
^{\alpha} } \right|}}} }}}.
\end{equation}
for all $i = \overline {1,n} $, $j = \overline {1,m} $.
\end{theorem}

{\textit{Proof}}. At first we shall obtain the representation
(\ref{eq:det_repr_A*A}). If $\lambda \in {\mathbb R} _{ +}  $, then the matrix
$\left( {\lambda {\rm {\bf I}} + {\rm {\bf A}}^{ *} {\rm {\bf A}}}
\right)\in {\mathbb C}^{n\times n} $ is Hermitian and $\rank
\left( {\lambda {\rm {\bf I}} + {\rm {\bf A}}^{ *} {\rm {\bf A}}}
\right)=n$.
 Hence,   there exists its  inverse
\[
\left( {\lambda {\rm {\bf I}} + {\rm {\bf A}}^{ *} {\rm {\bf A}}}
\right)^{ - 1} = {\frac{{1}}{{\det \left( {\lambda {\rm {\bf I}} +
{\rm {\bf A}}^{ * }{\rm {\bf A}}} \right)}}}\left(
{{\begin{array}{*{20}c}
 {L_{11}}  \hfill & {L_{21}}  \hfill & {\ldots}  \hfill & {L_{n\,1}}  \hfill
\\
 {L_{12}}  \hfill & {L_{22}}  \hfill & {\ldots}  \hfill & {L_{n\,2}}  \hfill
\\
 {\ldots}  \hfill & {\ldots}  \hfill & {\ldots}  \hfill & {\ldots}  \hfill
\\
 {L_{1\,n}}  \hfill & {L_{2\,n}}  \hfill & {\ldots}  \hfill & {L_{n\,n}}  \hfill
\\
\end{array}} } \right),
\]
 where $L_{ij} $ $(\forall i,j=\overline{1,n})$ is a cofactor
  in $\lambda {\rm {\bf
I}} + {\rm {\bf A}}^{ *} {\rm {\bf A}}$.   By Lemma \ref{kyrc2},
${\rm {\bf A}}^{ +}  = {\mathop {\lim }\limits_{\lambda \to 0}}
\left( {\lambda {\rm {\bf I}} + {\rm {\bf A}}^{ * }{\rm {\bf A}}}
\right)^{ - 1}{\rm {\bf A}}^{ *} $, so that
\begin{equation}
\label{eq3} {\rm {\bf A}}^{ +}  ={\mathop {\lim} \limits_{\lambda
\to 0}}
\begin{pmatrix}
 \frac{\det \left( {\lambda {\rm {\bf I}} + {\rm {\bf A}}^{ *} {\rm {\bf A}}}
\right)_{.1} \left( {{\rm {\bf a}}_{.\,1}^{ *} }  \right)}{{\det
\left( {\lambda {\rm {\bf I}} + {\rm {\bf A}}^{ *} {\rm {\bf A}}}
\right)}}
 & \ldots & \frac{\det \left( {\lambda {\rm {\bf I}} + {\rm {\bf
A}}^{ *} {\rm {\bf A}}} \right)_{.\,1} \left( {{\rm {\bf
a}}_{.\,m}^{ *} }  \right)}{{\det \left( {\lambda {\rm {\bf I}} +
{\rm {\bf
A}}^{ *} {\rm {\bf A}}} \right)}}\\
  \ldots & \ldots & \ldots \\
\frac{\det \left( {\lambda {\rm {\bf I}} + {\rm {\bf A}}^{ *} {\rm
{\bf A}}} \right)_{.\,n} \left( {{\rm {\bf a}}_{.\,1}^{ *} }
\right)}{{\det \left( {\lambda {\rm {\bf I}} + {\rm {\bf A}}^{ *}
{\rm {\bf A}}} \right)}} &
 \ldots &
 \frac{\det \left( {\lambda {\rm {\bf I}} + {\rm {\bf A}}^{ *} {\rm {\bf
A}}} \right)_{.\,n} \left( {{\rm {\bf a}}_{.\,m}^{ *} }
\right)}{{\det \left( {\lambda {\rm {\bf I}} + {\rm {\bf A}}^{ *}
{\rm {\bf A}}} \right)}}
\end{pmatrix}.
\end{equation}
 From Theorem \ref{th:char_pol} we get
\[ \det \left( {\lambda{\rm {\bf I}} + {\rm {\bf A}}^{ *}
{\rm {\bf A}}} \right) = \lambda ^{n} + d_{1} \lambda ^{n - 1} +
d_{2} \lambda ^{n - 2} + \ldots + d_{n},
\]
 where $d_{r} $ $(\forall r=\overline{1,n-1})$ is a
sum of  principal minors  of ${\rm {\bf A}}^{ *} {\rm {\bf A}}$ of
order $r$  and $d_{n}=\det  {\rm {\bf A}}^{ *} {\rm {\bf A}}$.
Since $\rank{\rm {\bf A}}^{ *} {\rm {\bf A}} = \rank{\rm {\bf A}}
= r$, then $d_{n} = d_{n - 1} = \ldots = d_{r + 1} = 0$ and
\begin{equation}
\label{eq4}\det \left( {\lambda {\rm {\bf I}} + {\rm {\bf A}}^{ *}
{\rm {\bf A}}} \right) = \lambda ^{n} + d_{1} \lambda ^{n - 1} +
d_{2} \lambda ^{n - 2} + \ldots + d_{r} \lambda ^{n - r}.
 \end{equation}
In the same way,  we have for  arbitrary $ 1\leq i \leq n$ and
$1\leq j\leq m $ from Lemma \ref{lem:char_mpdet}
\[ \det \left( {\lambda {\rm {\bf I}} + {\rm {\bf
A}}^{ *} {\rm {\bf A}}} \right)_{.\,i} \left( {{\rm {\bf
a}}_{.j}^{ *} } \right) = l_{1}^{\left( {ij} \right)} \lambda ^{n
- 1} + l_{2}^{\left( {ij} \right)} \lambda ^{n - 2} + \ldots +
l_{n}^{\left( {ij} \right)}, \]
 where for an arbitrary $1\leq k \leq n - 1$,\,
 $l_{k}^{\left( {ij}
\right)} = {\sum\limits_{\beta \in J_{k,\,n} {\left\{ {i}
\right\}}} {{\left| {\left( {({\rm {\bf A}}^{ *} {\rm {\bf
A}})_{.\,i} ({\rm {\bf a}}_{.j}^{ *}  )} \right)
}_{\beta} ^{\beta}  \right|}}}$, and \,\,$l_{n}^{\left( {i\,j} \right)} = \det
\left( {{\rm {\bf A}}^{ *} {\rm {\bf A}}} \right)_{.\,i} \left(
{{\rm {\bf a}}_{.\,j}^{ *} } \right)$. By Lemma \ref{kyrc5},
$\rank\left( {{\rm {\bf A}}^{ *} {\rm {\bf A}}} \right)_{.\,i}
\left( {{\rm {\bf a}}_{.\,j}^{ *} } \right) \le r$ so that if $
k>r$, then ${{\left| {\left( {({\rm {\bf A}}^{ *} {\rm {\bf
A}})_{\,.\,i} ({\rm {\bf a}}_{.j}^{ *}  )} \right) }_{\beta}
^{\beta} \right|}}= 0$, $(\forall \beta \in J_{k,\,n} {\left\{
{i} \right\}}, \forall i = \overline {1,n}, \forall j = \overline
{1,m})$. Therefore if $r + 1 \le k < n$, then $l_{k}^{\left( {ij}
\right)} = {\sum\limits_{\beta \in J_{k,\,n} {\left\{ {i}
\right\}}} {{\left| {\left( {({\rm {\bf A}}^{ *} {\rm {\bf
A}})_{\,.\,i} ({\rm {\bf a}}_{.j}^{ *}  )} \right) }_{\beta}
^{\beta} \right|}}}= 0$ and $l_{n}^{\left( {i\,j} \right)} =
\det \left( {{\rm {\bf A}}^{ *} {\rm {\bf A}}} \right)_{.\,i}
\left( {{\rm {\bf a}}_{.\,j}^{ *} } \right) = 0$, $\left( {\forall
i = \overline {1,n},\, \forall j = \overline {1,m}} \right)$.
Finally we obtain
\begin{equation}\label{eq5}
  \det \left( {\lambda {\rm {\bf I}} + {\rm {\bf A}}^{ *} {\rm {\bf
A}}} \right)_{.\,i} \left( {{\rm {\bf a}}_{.\,j}^{ *} } \right) =
l_{1}^{\left( {i\,j} \right)} \lambda ^{n - 1} + l_{2}^{\left(
{i\,j} \right)} \lambda ^{n - 2} + \ldots + l_{r}^{\left( {ij}
\right)} \lambda ^{n - r}.
\end{equation}

 By replacing the denominators and the numerators of the fractions in entries
 of  matrix (\ref{eq3}) with  the
 expressions (\ref{eq4}) and (\ref{eq5}) respectively, we get
\[
   {\rm {\bf A}}^{ +}  = {\mathop {\lim} \limits_{\lambda \to 0}}
\begin{pmatrix}
 {{\frac{{l_{1}^{\left( {11} \right)} \lambda ^{n - 1} + \ldots +
l_{r}^{\left( {11} \right)} \lambda^{n - r}}}{{\lambda ^{n} +
d_{1} \lambda ^{n - 1} + \ldots + d_{r} \lambda ^{n - r}}}}}
 & \ldots & {{\frac{{l_{1}^{\left( {1m} \right)} \lambda ^{n
- 1} + \ldots + l_{r}^{\left( {1m} \right)} \lambda ^{n -
r}}}{{\lambda ^{n} + d_{1} \lambda
^{n - 1} + \ldots + d_{r} \lambda ^{n - r}}}}} \\
  \ldots & \ldots & \ldots \\
  {{\frac{{l_{1}^{\left( {n1} \right)} \lambda ^{n - 1} + \ldots +
l_{r}^{\left( {n1} \right)} \lambda ^{n - r}}}{{\lambda ^{n} +
d_{1} \lambda ^{n - 1} + \ldots + d_{r} \lambda ^{n - r}}}}} &
\ldots & {{\frac{{l_{1}^{\left( {nm} \right)} \lambda ^{n - 1} +
\ldots + l_{r}^{\left( {nm} \right)} \lambda ^{n - r}}}{{\lambda
^{n} + d_{1} \lambda ^{n - 1} + \ldots + d_{r} \lambda ^{n -
r}}}}}
\end{pmatrix}=\]
  \[ = \left( {{\begin{array}{*{20}c}
 {{\frac{{l_{r}^{\left( {11} \right)}} }{{d_{r}} }}} \hfill & {\ldots}
\hfill & {{\frac{{l_{r}^{\left( {1m} \right)}} }{{d_{r}} }}} \hfill \\
 {\ldots}  \hfill & {\ldots}  \hfill & {\ldots}  \hfill \\
 {{\frac{{l_{r}^{\left( {n1} \right)}} }{{d_{r}} }}} \hfill & {\ldots}
\hfill & {{\frac{{l_{r}^{\left( {nm} \right)}} }{{d_{r}} }}} \hfill \\
\end{array}} } \right).
\]
From here it follows (\ref{eq:det_repr_A*A}).

 We can prove   (\ref{eq:det_repr_AA*}) in the same way. $\blacksquare$
 \begin{corollary}\label{kyrc7}
If ${\rm {\bf A}}\in {\mathbb C}^{m\times n}_{r}$ and $r < \min
{\left\{ {m,n} \right\}}$  or $r = m < n $, then the projection
matrix ${\rm {\bf P}} = {\rm {\bf A}}^{ +} {\rm {\bf A}}$ can be
represented as
\[
{\rm {\bf P}} = \left({\frac{{p_{ij}}}{{d_{r} \left( {{\rm {\bf
A}}^{ *} {\rm {\bf A}}} \right)}}}\right)_{n\times n},
\]
 where ${\rm {\bf d}}_{.\,j} $ denotes the
$j$th column of $({\rm {\bf A}}^{ *} {\rm {\bf A}})$ and, for
arbitrary $1\leq i,j \leq n $,
   $p_{ij} ={\sum\limits_{\beta \in
J_{r,n} {\left\{ {i} \right\}}} {{\left| {\left( {({\rm {\bf A}}^{
*} {\rm {\bf A}})_{\,.\,i} ({\rm {\bf d}}_{.j}  )} \right)_{\beta}
^{\beta} }  \right|}}}$.
\end{corollary}

{\textit{Proof}}. Representing the Moore - Penrose inverse ${\rm
{\bf A}}^{+}$ by (\ref{eq:det_repr_A*A}), we obtain
\[
{\rm {\bf P}}={\frac{{1}}{{d_{r} \left( {{\rm {\bf A}}^{ *} {\rm
{\bf A}}} \right)}}}\begin{pmatrix}
  l_{11} & l_{12} & \ldots & l_{1m} \\
  l_{21} & l_{22} & \ldots & l_{2m} \\
  \ldots & \ldots & \ldots & \ldots \\
  l_{n1} & l_{n2} & \ldots & l_{nm}
\end{pmatrix}\begin{pmatrix}
  a_{11} & a_{12} & \ldots & a_{1n} \\
  a_{21} & a_{22} & \ldots & a_{2n} \\
  \ldots & \ldots & \ldots & \ldots \\
  a_{m\,1} & a_{m\,2} & \ldots & a_{m\,n}
\end{pmatrix}.
\]
Therefore, for arbitrary $1\leq i,j\leq n$ we get
 \[\begin{array}{c}
  p_{i\,j}=\sum\limits_{k}
 {{\sum\limits_{\beta \in J_{r,\,
n} {\left\{ {i} \right\}}} {{\left| {\left( {({\rm {\bf A}}^{ *}
{\rm {\bf A}})_{.\,i} ({\rm {\bf a}}_{.\,k}^{ *}  )}
\right)_{\beta} ^{\beta} }  \right|}}}} \cdot a_{k\,j}=\\
  =\sum\limits_{\beta \in J_{r,\,
n} {\left\{ {i} \right\}}} \sum\limits_{k}
 {{ {{\left|
 {\left( {({\rm {\bf A}}^{ *}
{\rm {\bf A}})_{.\,i} ({\rm {\bf a}}_{.\,k}^{ *}\cdot a_{k\,j}  )}
\right)_{\beta} ^{\beta} } \right|}}}}={\sum\limits_{\beta \in
J_{r,\,n} {\left\{ {i} \right\}}} {{\left| {\left( {({\rm {\bf
A}}^{ *} {\rm {\bf A}})_{.\,i} ({\rm {\bf d}}_{.j}^{ *}  )}
\right)_{\beta} ^{\beta} } \right|}}}.\blacksquare
\end{array}
\]

Using the representation (\ref{eq:det_repr_AA*}) of the Moore - Penrose
inverse the following corollary can be proved in the same way.
 \begin{corollary} \label{kyrc8} If ${\rm
{\bf A}}\in {\mathbb C}^{m\times n}_{r}$, where $r < \min {\left\{
{m,n} \right\}}$ or $r = n < m$,
 then a projection matrix ${\rm {\bf Q}} = {\rm
{\bf A}}{\rm {\bf A}}^{ +} $  can be represented as
\[
{\rm {\bf Q}} =  \left({\frac{{q_{ij}}}{{d_{r} \left( {{\rm {\bf
A}}{\rm {\bf A}}^{ *} } \right)}}}\right)_{m\times m},
\]
where ${\rm {\bf g}}_{i.} $ denotes the $i$th row of $({\rm {\bf
A}}{\rm {\bf A}}^{ *} )$ and,  for arbitrary $1\leq i,j \leq m $,
$q_{i\,j} ={\sum\limits_{\alpha \in I_{r,m} {\left\{ {j}
\right\}}} {{\left| {\left( {({\rm {\bf A}} {\rm {\bf A}}^{
*})_{j.}\, ({\rm {\bf g}}_{i.}\,  )} \right)_{\alpha} ^{\alpha} }
\right|}}} $.
\end{corollary}
\begin{remark}\label{rem:det_repr_MP_A*A} If $\rank{\rm {\bf A}} =
n$, then   from   Corollary \ref{cor:A+A-1} we get ${\rm {\bf A}}^{ +}  =
\left( {{\rm {\bf A}}^{ *} {\rm {\bf A}}} \right)^{ - 1}{\rm {\bf
A}}^{ *} $. Representing $({\rm {\bf A}}^{ *} {\rm {\bf A}})^{-1}$
by the classical adjoint matrix, we  have
 \begin{equation}
\label{eq:A*A_ful} {\rm {\bf A}}^{ +}  = {\frac{{1}}{{\det ({\rm {\bf
A}}^{ *} {\rm {\bf A}})}}}
\begin{pmatrix}
   {\det ({\rm {\bf A}}^{ *} {\rm {\bf A}})_{.1} \left( {{\rm {\bf a}}_{.1}^{ *} }  \right)}&
   \ldots & {\det ({\rm
{\bf A}}^{ *} {\rm {\bf A}})_{.1} \left( {{\rm {\bf
a}}_{.\;m}^{ *} }  \right)}\\
  \ldots & \ldots & \ldots \\
  {\det ({\rm {\bf A}}^{ *} {\rm {\bf A}})_{.\,n} \left( {{\rm {\bf a}}_{.\,1}^{ *} }  \right)}
   & \ldots
   & {\det ({\rm {\bf A}}^{ *} {\rm {\bf A}})_{.\,n} \left( {{\rm {\bf a}}_{.\,m}^{ *} }  \right)}
\end{pmatrix}.\end{equation}
 If $n < m$, then  (\ref{eq:det_repr_A*A}) is
 valid.
 \end{remark}
 \begin{remark}\label{rem2}
 As above, if $\rank{\rm {\bf A}} = m$, then
\begin{equation}
\label{eq:AA*_ful} {\rm {\bf A}}^{ +}  = {\frac{{1}}{{\det ({\rm {\bf
A}}{\rm {\bf A}}^{ *} )}}}\begin{pmatrix}
  {\det ({\rm {\bf A}}{\rm {\bf A}}^{ *} )_{1\,.} \left( {{\rm {\bf a}}_{1\,.}^{ *} }  \right)}
   & \ldots &
   {\det ({\rm {\bf A}}{\rm {\bf A}}^{ *} )_{m\,.} \left( {{\rm {\bf a}}_{1\,.}^{ *} }  \right)} \\
 \ldots & \ldots & \ldots \\
 {\det ({\rm {\bf A}}{\rm {\bf A}}^{ *} )_{1\,.} \left( {{\rm {\bf a}}_{n\,.}^{ *} }  \right)}
  & \ldots &
{\det ({\rm {\bf A}}{\rm {\bf A}}^{ *} )_{m\,.} \left( {{\rm {\bf
a}}_{n\,.}^{ *} }  \right)}
\end{pmatrix}.
\end{equation}
If  $n > m$, then  (\ref{eq:det_repr_AA*}) is valid as well.
 \end{remark}

 \begin{remark} By definition of the classical adjoint $Adj({\rm {\bf A}})$
  for an arbitrary invertible matrix ${\rm {\bf A}}\in {\mathbb C}^{n\times
  n}$ one may put, $Adj({\rm {\bf A}})\cdot{\rm {\bf A}}=\det{\rm {\bf A}}\cdot{\rm {\bf I}}_{n}$.

  If ${\rm {\bf A}}\in {\mathbb C}^{m\times n}$ and $\rank{\rm {\bf
A}} = n$, then by  Corollary \ref{cor:A+A-1}, ${\rm {\bf A}}^{ +} {\rm
{\bf A}} = {\rm {\bf I}}_{n}$.
  Representing the matrix ${\rm {\bf A}}^{ +} $
 by  (\ref{eq:A*A_ful}) as ${\rm
{\bf A}}^{ +}  = {\frac{{{\rm {\bf L}}}}{{\det \left( {{\rm {\bf
A}}^{ *} {\rm {\bf A}}} \right)}}}$, we obtain
 ${\rm {\bf L}}{\rm {\bf A}} = \det \left(
{{\rm {\bf A}}^{ *} {\rm {\bf A}}} \right) \cdot {\rm {\bf I}}_{n}
$. This means that the matrix ${\rm {\bf L}}=(l_{ij})\in {\mathbb C}^{n\times m}$
 is a
 left  analogue of $Adj({\rm {\bf A}})$ , where ${\rm {\bf A}}\in {\mathbb C}^{m\times n}_{n}$, and $l_{ij}={\det ({\rm {\bf A}}^{ *} {\rm {\bf A}})_{.i} \left( {{\rm {\bf a}}_{.j}^{ *} }  \right)}$ for all $i = \overline {1,n} $, $j = \overline {1,m} $.

  If $\rank{\rm {\bf A}} = m$, then by  Corollary \ref{cor:A+A-1}, ${\rm
{\bf A}}{\rm {\bf A}}^{ +} = {\rm {\bf I}}_{m}$. Representing the
matrix ${\rm {\bf A}}^{ +} $
 by  (\ref{eq:AA*_ful}) as ${\rm {\bf A}}^{ +}  = {\frac{{{\rm
{\bf R}}}}{{\det \left( {{\rm {\bf A}}{\rm {\bf A}}^{ *} }
\right)}}}$, we obtain ${\rm {\bf A}}{\rm {\bf R}} = {\rm {\bf
I}}_{m} \cdot \det \left( {{\rm {\bf A}}{\rm {\bf A}}^{ *} }
\right)$. This means that the matrix ${\rm {\bf R}}=(r_{ij})\in {\mathbb C}^{m\times n}$ is a
 right analogue of $Adj({\rm {\bf A}})$, where ${\rm {\bf A}}\in {\mathbb C}^{m\times n}_{m}$, and $r_{ij}={\det ({\rm {\bf A}}{\rm {\bf A}}^{ *} )_{j\,.} \left( {{\rm {\bf a}}_{i\,.}^{ *} }  \right)}$ for all $i = \overline {1,n} $, $j = \overline {1,m} $.

  If
${\rm} {\rm {\bf A}}\in {\mathbb C}^{m\times n}_{r} $ and
$r<\min\{m, n\}$, then by (\ref{eq:det_repr_A*A}) we have ${\rm {\bf A}}^{ +}
= {\frac{{{\rm {\bf L}}}}{{d_{r} \left( {{\rm {\bf A}}^{ *} {\rm
{\bf A}}} \right)}}}$,
 where ${\rm {\bf L}} = \left( {l_{ij}}  \right)\in {\mathbb C}^{n\times m}$ and $l_{ij}={{{\sum\limits_{\beta
\in J_{r,\,n} {\left\{ {i} \right\}}} {{\left| \left( {\left(
{{\rm {\bf A}}^{ *} {\rm {\bf A}}} \right)_{\,. \,i} \left( {{\rm
{\bf a}}_{.j}^{ *} }  \right)} \right) {\kern 1pt}
_{\beta} ^{\beta} \right|}} }  }}$ for all $i = \overline {1,n} $, $j = \overline {1,m} $.
  From Corollary \ref{kyrc7}  we get ${\rm
{\bf L}}{\rm {\bf A}} = d_{r} \left( {\rm {\bf A}}^{ *} {\rm {\bf
A}} \right) \cdot {\rm {\bf P}}$.  The matrix ${\rm {\bf P}}$ is
idempotent. All eigenvalues of an idempotent matrix chose from 1
or 0 only. Thus, there exists an unitary matrix ${\rm {\bf U}}$
such that
\[{\rm {\bf L}}{\rm {\bf A}} = d_{r} \left( {\rm {\bf
A}}^{ *} {\rm {\bf A}} \right) {\rm {\bf U}} {{\rm {\bf
diag}}\left( {{\rm 1},\ldots,{\rm 1},{\rm 0},\ldots ,{\rm 0}}
\right)}  {\rm {\bf U}}^{ *} ,\]
where $
 {{\rm {\bf diag}}\left( {{\rm 1},\ldots,{\rm 1},{\rm
0},\ldots ,{\rm 0}} \right)} \in {\mathbb C}^{n\times n}$ is a
diagonal matrix.  Therefore,  the matrix ${\rm {\bf L}}$
 can be considered  as a left
analogue of $Adj({\rm {\bf A}})$, where ${\rm {\bf A}}\in {\mathbb
C}^{m\times n}_{r}$.

 In the same
way, if ${\rm} {\rm {\bf A}}\in {\mathbb C}^{m\times n}_{r} $ and
$r<\min\{m, n\}$, then by   (\ref{eq:det_repr_A*A}) we have ${\rm {\bf A}}^{
+}  = {\frac{{{\rm {\bf R}}}}{{d_{r} \left( {{\rm {\bf A}}{\rm
{\bf A}}^{ *} } \right)}}}$, where ${\rm {\bf R}} = \left(
{r_{ij}} \right)\in {\mathbb C}^{n\times m}$, $r_{ij}={{{\sum\limits_{\alpha \in I_{r,m} {\left\{ {j} \right\}}}
{{\left| \left( {({\rm {\bf A}}{\rm {\bf A}}^{ *} )_{j\,.\,} ({\rm
{\bf a}}_{i.\,}^{ *} )} \right)\,_{\alpha} ^{\alpha}\right|} }
}}}$for all $i = \overline {1,n} $, $j = \overline {1,m} $. From Corollary
\ref{kyrc8}  we get  ${\rm {\bf A}}{\rm {\bf R}} = d_{r} \left(
{\rm {\bf A}}{\rm {\bf A}}^{ *} \right) \cdot {\rm {\bf Q}}$. The
matrix ${\rm {\bf Q}}$ is idempotent. There exists an unitary
matrix ${\rm {\bf V}}$ such that
\[{\rm {\bf A}}{\rm {\bf R}} =
d_{r} \left( {\rm {\bf A}} {\rm {\bf A}}^{ *} \right) {\rm {\bf
V}} {{\rm {\bf diag}}\left( {{\rm 1},\ldots ,{\rm 1},{\rm
0},\ldots,{\rm 0}} \right)}{\rm {\bf V}}^{ *} ,\]
 where $ {{\rm
{\bf diag}}\left( {{\rm 1},\ldots ,{\rm 1},{\rm 0},\ldots,{\rm 0}}
\right)} \in {\mathbb C}^{m\times m}$. Therefore,  the matrix
${\rm {\bf R}}$ can be considered  as a right analogue of
$Adj({\rm {\bf A}})$ in this case.
\end{remark}

 \begin{remark} To obtain an entry of ${\rm {\bf A}}^{ +}$ by
 Theorem \ref{theor:st} one calculates
 $(C^{r}_{n}C^{r}_{m}+C^{r-1}_{n-1}C^{r-1}_{m-1})$ determinants of
 order $r$. Whereas by the equation (\ref{eq:det_repr_A*A}) we calculate as much
 as
 $(C^{r}_{n}+C^{r-1}_{n-1})$ determinants of
 order $r$ or  we calculate the total of
 $(C^{r}_{m}+C^{r-1}_{m-1})$ determinants by (\ref{eq:det_repr_AA*}). Therefore the
 calculation of entries of ${\rm {\bf A}}^{ +} $ by Theorem \ref{theor:det_repr_MP} is easier
  than by Theorem \ref{theor:st}.
  \end{remark}
\subsection{Analogues of the classical adjoint matrix for the weighted Moore-Penrose inverse}
Let Hermitian positive definite matrices ${\bf M}$ and ${\bf N}$ of order $m$ and $n$, respectively, be given. The weighted Moore-Penrose inverse ${\bf X}={\bf A}^{+}_{M,N}$ can be explicitly expressed from the weighted singular value decomposition due to Van Loan \cite{lo}
\begin{lemma}Let ${\bf A}\in {\mathbb C}^{m\times n}_{r}$. There exists ${\bf U}\in {\mathbb C}^{m\times m}$, ${\bf V}\in {\mathbb C}^{n\times n}$ satisfying ${\bf U}^{*}{\bf M}{\bf U}={\bf I}_{m}$ and ${\bf V}^{*}{\bf N}^{-1}{\bf V}={\bf I}_{n}$ such that
\[
{\bf A}={\bf U}\left(
                 \begin{array}{cc}
                   {\bf D} & {\bf 0} \\
                   {\bf 0} & {\bf 0} \\
                 \end{array}
               \right)
{\bf V}^{*}
\]
Then the weighted Moore-Penrose inverse ${\bf A}^{+}_{M,N}$ can be represented
\[
{\bf A}^{+}_{M,N}={\bf N}^{-1}{\bf V}\left(
                 \begin{array}{cc}
                   {\bf D}^{-1} & {\bf 0} \\
                   {\bf 0} & {\bf 0} \\
                 \end{array}
               \right)
{\bf U}^{*}{\bf M},
\]
where ${\bf D}=diag(\sigma_{1}, \sigma_{2},...,\sigma_{r})$, $\sigma_{1}\geq \sigma_{2}\geq...\geq\sigma_{r}>0$ and $\sigma^{2}$ is the nonzero eigenvalues of ${\bf N}^{-1}{\bf A}^{*}{\bf M}{\bf A}$.
\end{lemma}
For the weighted Moore-Penrose inverse ${\bf X}={\bf A}^{+}_{M,N}$  we have  the following limit representation
\begin{lemma}(\cite{wei}, Corollary 3.4.) Let ${\bf A}\in {\mathbb C}^{m\times n}$, ${\bf A}^{\sharp}={\bf N}^{-1}{\bf A}^{*}{\bf M}$. Then
\[{\bf A}^{+}_{M,N}=\lim_{\lambda \to 0}(\lambda {\bf I}+{\bf A}^{\sharp}{\bf A})^{-1}{\bf A}^{\sharp}
\]
\end{lemma}

Denote by ${\rm {\bf a}}_{.j}^{\sharp} $ and ${\rm {\bf
a}}_{i.}^{\sharp} $ the $j$th column  and the $i$th row of  ${\rm
{\bf A}}^{\sharp} $ respectively. By putting ${\bf A}^{\sharp}$ instead ${\bf A}^{*}$, we obtain the proofs of the following  two lemmas and theorem similar to the proofs of Lemmas \ref{kyrc5}, \ref{lem:char_mpdet} and Theorem \ref{theor:det_repr_MP}, respectively.
\begin{lemma} \label{lem:rank_wmpcol} If ${\bf A}^{\sharp}{\bf A}\in {\mathbb{C}}^{n\times n}$, then
\[
 \rank\,\left( {{\bf A}^{\sharp}{\bf A} }
\right)_{.\,i} \left( { {{\bf a}}_{.j}^{ \sharp} }  \right) \le
\rank\,\left( {{\bf A}^{\sharp}{\bf A} } \right).
\]
for all  $i = \overline{1,n}$ and $j = \overline {1,m}$
\end{lemma}

Analogues of the characteristic polynomial are considered in the following lemma.
\begin{lemma}\label{lem:char_wmpdet}
If ${\bf A}\in
{\mathbb{C}}^{m\times n}$ and $\lambda \in \mathbb{R}$, then
\[\det \left(  \left(\lambda {\rm {\bf I}}_n + {\bf A}^{
\sharp}{\bf A}  \right)_{.\,i} \left( { {{\bf a}}_{.j}^{ \sharp} } \right)\right) =
c_{1}^{\left( {ij} \right)} \lambda^{n - 1} + c_{2}^{\left( {ij} \right)}
\lambda^{n - 2} + \ldots + c_{n}^{\left( {ij} \right)},
\]
\noindent where $c_{n}^{\left( {ij} \right)} = \left| \left( {{\bf A}^{
\sharp}{\bf A} } \right)_{.\,i} \left( { {{\bf a}}_{.j}^{ \sharp} } \right)\right|$
and $c_{s}^{\left( {ij} \right)} = {\sum\limits_{\beta \in J_{s,\,n}
{\left\{ {i} \right\}}} {  \left|\left( {\left( {{\bf A}^{ \sharp}{\bf A} }
\right)_{.\,i} \left( { {{\bf a}}_{.j}^{ \sharp} }  \right)}\right)_{\beta} ^{\beta} \right|
 } }$ for all $s = \overline {1, n- 1} $, $i = \overline
{1,n}$, and $j = \overline {1,m}$.
\end{lemma}

The following theorem  introduce the determinantal representations of the weighted Moore-Penrose by analogs of the classical adjoint matrix.

\begin{theorem}\label{theor:det_repr_wMP}
If ${\rm {\bf A}} \in {\rm {\mathbb{C}}}_{r}^{m\times n} $ and
$r<min\{m,n\}$, then the weighted Moore-Penrose inverse  ${\rm {\bf A}}^{
+}_{M,N} = \left( {\tilde{a}_{ij}^{ +} } \right) \in {\rm
{\mathbb{C}}}_{}^{n\times m} $ possess the following determinantal
representation:
\begin{equation}
\label{eq:det_repr_wA*A}
 \tilde{a}_{ij}^{ +}  = {\frac{{{\sum\limits_{\beta
\in J_{r,\,n} {\left\{ {i} \right\}}} {{\left| \left( {\left(
{{\rm {\bf A}}^{ \sharp} {\rm {\bf A}}} \right)_{\,. \,i} \left( {{\rm
{\bf a}}_{.j}^{ \sharp} }  \right)} \right) {\kern 1pt}
_{\beta} ^{\beta} \right|}} }  }}{{{\sum\limits_{\beta \in
J_{r,\,\,n}} {{\left| {\left( {{\rm {\bf A}}^{ \sharp} {\rm {\bf A}}}
\right){\kern 1pt} _{\beta} ^{\beta} }  \right|}}} }}},
\end{equation}
for all $i = \overline {1,n} $, $j = \overline {1,m} $.
\end{theorem}

\subsection{Analogues of the classical adjoint matrix for the Drazin
inverse}
The Drazin inverse can be represented explicitly by the Jordan
canonical form as follows.
\begin{theorem}\cite{cam} If ${\rm {\bf A}} \in {\mathbb C}^{n\times n}$ with
$Ind{\kern 1pt} {\rm {\bf A}} = k $ and
\[
{\rm {\bf A}} = {\rm {\bf P}}\begin{pmatrix}
  {\rm {\bf C}} & {\rm {\bf 0}} \\
  {\rm {\bf 0}}& {\rm {\bf N}}
\end{pmatrix} {\rm {\bf P}}^{ - 1}
\]
where ${\rm {\bf C}}$ is nonsingular and $\rank{\rm {\bf C}} =
\rank{\rm {\bf A}}^{k}$, and ${\rm {\bf N}}$ is nilpotent of order
$k$, then
\begin{equation}
\label{eq:jord_dr}
{\rm {\bf A}}^{D} = {\rm {\bf P}}\begin{pmatrix}
  \rm {\bf C}^{ - 1} & {\rm {\bf 0}} \\
  {\rm {\bf 0}}& {\rm {\bf 0}}
\end{pmatrix} {\rm {\bf P}}^{ - 1}.
\end{equation}
\end{theorem}

 Stanimirovic' \cite{st2} introduced a determinantal representation of the Drazin inverse
by the following theorem.
\begin{theorem}
The Drazin inverse ${\rm {\bf A}}^{D}=\left(a_{ij}^{D}\right)$ of
an arbitrary matrix ${\rm {\bf A}}\in {\mathbb C}^{n\times n} $
with $Ind{\rm {\bf A}}=k $ possesses the following determinantal
representation
\begin{equation}\label{eq:det_dr_st}
a_{ij}^{D} = {\frac{{{\sum\limits_{\left( {\alpha,\,\beta} \right)
\in N_{r_{k}}  {\left\{ {j,\,i} \right\}}} {{\left| {\left( {{\rm
{\bf A}}^{s}} \right)_{\alpha} ^{\beta} }
\right|}{\frac{{\partial} }{{\partial a_{j\,i} }}}{\left| {{\rm
{\bf A}}_{\beta} ^{\alpha} }  \right|}} }}}{{{\sum\limits_{\left(
{\gamma,\,\delta}  \right) \in N_{r_{k}} } {{\left| {\left( {{\rm
{\bf A}}^{s}} \right)_{\gamma} ^{\delta} } \right|}\,{\left| {{\rm
{\bf A}}_{\delta} ^{\gamma} }  \right|}}} }}}, \,\,1 \leq i, j
\leq n;
\end{equation}
where $s \ge k$ and $r_{k} = rank{\rm {\bf A}}^{s}$. \end{theorem}
This determinantal representations  of the Drazin  inverse is
based on a full-rank representation.

We use the following  limit representation of the
Drazin inverse.
\begin{lemma} \cite{ca}\label{theor:lim_draz} If ${\rm {\bf A}}
 \in {\mathbb C}^{n\times n}$, then
\[
{\rm {\bf A}}^{D} = {\mathop {\lim} \limits_{\lambda \to 0}}
\left( {\lambda {\rm {\bf I}}_n + {\rm {\bf A}}^{k + 1}} \right)^{
- 1}{\rm {\bf A}}^{k},
\]
where $k = Ind{\kern 1pt} {\rm {\bf A}}$, $\lambda \in {\mathbb R}
_{ +}  $, and ${\mathbb R} _{ +} $ is a set of the real positive
numbers.
\end{lemma}
Since the equation  (\ref{eq6:MP})  can be replaced by follows
\[ (6a)\,\,  {\rm {\bf X}}{\rm {\bf A}}^{k+1}={\rm {\bf
         A}}^{k},\]
the following  lemma can be obtained by analogy to Lemma
\ref{theor:lim_draz}.
\begin{lemma} \label{theor:lim_A_k_A_k+1} If ${\rm {\bf A}}
 \in {\mathbb C}^{n\times n}$, then
\[
{\rm {\bf A}}^{D} = {\mathop {\lim} \limits_{\lambda \to 0}}{\rm
{\bf A}}^{k} \left( {\lambda {\rm {\bf I}}_n + {\rm {\bf A}}^{k +
1}} \right)^{ - 1},
\]
where $k = Ind{\kern 1pt} {\rm {\bf A}}$, $\lambda \in {\mathbb R}
_{ +}  $, and ${\mathbb R} _{ +} $ is a set of the real positive
numbers.
\end{lemma}
 Denote by ${\rm {\bf a}}_{.j}^{(k)} $ and ${\rm {\bf
a}}_{i.}^{(k)} $ the $j$th column  and the $i$th row of ${\rm {\bf
A}}^{k} $ respectively.

We
consider the following auxiliary lemma.
\begin{lemma} \label{lem:rank_A_row} If ${\rm {\bf A}} \in {\mathbb C}^{n\times n}$
with $Ind{\kern 1pt} {\rm {\bf A}} = k $, then for all $i,j =
\overline {1,n}$
\[ \rank{\rm {\bf A}}_{i{\kern 1pt} .}^{k + 1}
\left( {{\rm {\bf a}}_{j.}^{\left( {k} \right)}}  \right) \le
\rank{\rm {\bf A}}^{k + 1}.
\]
\end{lemma}

{\textit{Proof.}}
  The matrix  $ {\rm {\bf A}}^{ k+1} _{i\,.} \left({{\rm {\bf a}}_{j\,.}^{ (k)} }  \right)$ may by represent as follows
\[\left(
    \begin{array}{ccc}
      {{\sum\limits_{s=1}^{n}a_{1s} {a_{s1}^{ (k)}  } } } & \ldots & {{\sum\limits_{s=1}^{n}a_{1s} {a_{sn}^{ (k)}  } } } \\
      \ldots & \ldots & \ldots \\
      a_{j1}^{ (k)}  & \ldots & a_{j\,n}^{ (k)}  \\
      \ldots & \ldots & \ldots \\
      {{\sum\limits_{s=1}^{n}a_{ns} {a_{s1}^{ (k)}  } } } & \ldots & {{\sum\limits_{s=1}^{n}a_{ns} {a_{sn}^{ (k)}  } } } \\
    \end{array}
  \right)
\]
Let ${\rm {\bf P}}_{l\,i} \left( {-a_{l\,j}}  \right)\in {\mathbb
C}^{n\times n} $, $(l \ne i )$, be a matrix with $-a_{l\,j} $ in
the $(l, i)$ entry, 1 in all diagonal entries, and 0 in others.
It is a matrix of an
elementary transformation. It follows that
\[
{\rm {\bf A}}^{ k+1} _{i\,.} \left({{\rm {\bf a}}_{j\,.}^{ (k)} }  \right) \cdot \prod\limits_{l \ne
i}{\rm {\bf P}}_{l\,i} \left( {-a_{l\,j}}  \right) =\left({
    \begin{array}{ccc}
      {{\sum\limits_{s\neq j}^{n} a_{1s} {a_{s1}^{ (k)}  } } } & \ldots & {{\sum\limits_{s\neq j}^{n}a_{1s} {a_{sn}^{ (k)}  } } } \\
      \ldots & \ldots & \ldots \\
      a_{j1}^{ (k)}  & \ldots & a_{j\,n}^{ (k)}  \\
      \ldots & \ldots & \ldots \\
      {{\sum\limits_{s\neq j}^{n} a_{ns} {a_{s1}^{ (k)}  } } } & \ldots & {{\sum\limits_{s\neq j}^{n}a_{ns} {a_{sn}^{ (k)}  } } } \\
    \end{array}}
  \right)\, ith
\]
 The obtained above matrix  has the following factorization.
\[
 \begin{pmatrix}
      {{\sum\limits_{s\neq j}^{n} a_{1s} {a_{s1}^{ (k)}  } } } & \ldots & {{\sum\limits_{s\neq j}^{n}a_{1s} {a_{sn}^{ (k)}  } } } \\
      \ldots & \ldots & \ldots \\
      a_{j1}^{ (k)}  & \ldots & a_{j\,n}^{ (k)}  \\
      \ldots & \ldots & \ldots \\
      {{\sum\limits_{s\neq j}^{n} a_{ns} {a_{s1}^{ (k)}  } } } & \ldots & {{\sum\limits_{s\neq j}^{n}a_{ns} {a_{sn}^{ (k)}  } } } \\
 \end{pmatrix} =
\]
\[ \begin{pmatrix}
      a_{11} &  \ldots & 0 &  \ldots & a_{1n} \\
       \ldots &  \ldots &  \ldots &  \ldots &  \ldots \\
      0 &  \ldots & 1 &  \ldots & 0 \\
       \ldots &  \ldots &  \ldots &  \ldots &  \ldots \\
      a_{n1} &  \ldots & 0 &  \ldots & a_{nn} \\
  \end{pmatrix}
    \begin{pmatrix}
      a_{11}^{(k)} &  a_{12}^{(k)}  & \ldots &  a_{1n}^{(k)}  \\
      a_{21}^{(k)}  &  a_{22}^{(k)}  & \ldots &  a_{2n}^{(k)}  \\
      \ldots & \ldots & \ldots & \ldots \\
       a_{n1}^{(k)}  &  a_{n2}^{(k)}  & \ldots &  a_{nn}^{(k)}  \\
    \end{pmatrix}
\]
 Denote the first matrix by
  \[{\rm {\bf \tilde {A}}}: = {\mathop
{\left( {{\begin{array}{*{20}c}
 {a_{11}}  \hfill & {\ldots}  \hfill & {0} \hfill & {\ldots}  \hfill &
{a_{1n}}  \hfill \\
 {\ldots}  \hfill & {\ldots}  \hfill & {\ldots}  \hfill & {\ldots}  \hfill &
{\ldots}  \hfill \\
 {0} \hfill & {\ldots}  \hfill & {1} \hfill & {\ldots}  \hfill & {0} \hfill
\\
 {\ldots}  \hfill & {\ldots}  \hfill & {\ldots}  \hfill & {\ldots}  \hfill &
{\ldots}  \hfill \\
 {a_{n1}}  \hfill & {\ldots}  \hfill & {0} \hfill & {\ldots}  \hfill &
{a_{nn}}  \hfill \\
\end{array}} } \right)}\limits_{jth}} ith.\]
The matrix
${\rm {\bf \tilde {A}}}$ is obtained from ${\rm {\bf A}}$ by
replacing all entries of the $i$th row  and the $j$th column with
zeroes except for 1 in the $(i, j)$ entry. Elementary
transformations of a matrix do not change its rank. It follows that
$\rank{\rm {\bf A}}^{ k+1} _{i\,.} \left({{\rm {\bf a}}_{j\,.}^{ (k)} }  \right) \le \min {\left\{
{\rank{\rm {\bf A}}^{ k },\rank{\rm {\bf \tilde {A}}}} \right\}}$.
Since $\rank{\rm {\bf \tilde {A}}} \ge \rank{\rm {\bf A}}^{ k} $ the proof is completed. $\blacksquare$

The following lemma is proved similarly.
\begin{lemma} \label{lem:rank_A_col} If ${\rm {\bf A}} \in {\mathbb C}^{n\times n}$
with $Ind{\kern 1pt} {\rm {\bf A}} = k $, then for all $i,j =
\overline {1,n}$
\[
\label{eq:rank_column} \rank{\rm {\bf A}}_{.{\kern 1pt} i}^{k + 1}
\left( {{\rm {\bf a}}_{.j}^{\left( {k} \right)}}  \right) \le
\rank{\rm {\bf A}}^{k + 1}.
\]
\end{lemma}

\begin{lemma}\label{lem:char_pol_dr}
If ${\bf A}\in {\mathbb{C}}^{n\times
n}$  and $\lambda \in \mathbb{R}$,
then
\begin{equation}
\label{eq:char_drdet}\det\left(  {( \lambda{\rm {\bf I}}_{n} +  {\bf A}^{k+1} )_{j\,.\,} ({\rm {\bf a}}_{i.}^{(k)}
)}\right)
   = r_{1}^{\left( {ij} \right)} \lambda^{n
- 1} +r_{2}^{\left( {ij} \right)} \lambda^{n - 2} + \ldots + r_{n}^{\left(
{ij} \right)},
\end{equation}
\noindent where  $r_{n}^{\left( {ij} \right)} = \left| {{\bf A}^{k+1} _{j\,.\,} ({\rm {\bf a}}_{i.\,}^{(k)} )}\right|$ and
$r_{s}^{\left( {ij} \right)} = \sum\limits_{\alpha \in I_{s,n} {\left\{
{j} \right\}}} \left| \left( {{\bf A}^{k+1}_{j\,.\,} ({\rm {\bf
a}}_{i.\,}^{(k)} )}\right)_{\alpha} ^{\alpha} \right| $ for all
$s = \overline {1,n - 1} $ and $i, j = \overline {1,n}$.
\end{lemma}
{\textit{Proof.}} Consider the  matrix $\left(  {( \lambda{\rm {\bf I}}_{n} +  {\bf A}^{k+1} )_{j\,.\,} ({\rm {\bf a}}_{j.}^{(k)}
)}\right)
\in {\rm {\mathbb{C}}}^{n\times n}$.  Taking into account Theorem
\ref{th:char_pol} we obtain
\begin{equation}
\label{eq:char_drdet1}  \left|\left(  {( \lambda{\rm {\bf I}}_{n} +  {\bf A}^{k+1} )_{j\,.\,} ({\rm {\bf a}}_{j.}^{(k)}
)}\right)\right| =
d_{1} \lambda^{n - 1} + d_{2} \lambda^{n - 2} + \ldots + d_{n},
\end{equation}
where $d_{s} = {\sum\limits_{\alpha \in I_{s,\,n} {\left\{ {j} \right\}}}
| ( {{\bf A}^{k+1} })_{\alpha} ^{\alpha}   | } $ is the
sum of all principal minors of order $s$ that contain the $j$-th row
for all $s = \overline {1, n - 1} $ and $d_{n} = \det  {{\bf A}^{k+1} }$. Since ${\rm {\bf a}}_{j.}^{(k+1)} =
 {\sum\limits_{l} { a_{jl}{\rm {\bf a}}_{l.}^{(k)} }
}$, where ${\rm {\bf a}}_{l.}^{(k)}  $ is the $l$th
row-vector  of ${\rm {\bf {A}}}^{k}$ for all
$l=\overline{1,n}$, then we have on the one hand
\begin{equation}
\label{eq:char_drdet2}
\begin{array}{c}
   \left|\left(  {( \lambda{\rm {\bf I}}_{n} +  {\bf A}^{k+1} )_{j\,.} ({\rm {\bf a}}_{j.}^{(k)}
)}\right)\right|
   = \sum\limits_{l} \left|  \left( {\lambda{\rm {\bf I}} + {\bf A}^{k+1} } \right)_{l.} \left( { a_{jl}{\rm {\bf a}}_{l.}^{(k)} } \right) \right|= \\ \sum\limits_{l}{ a_{jl} \cdot \left| {\left(
{\lambda{\rm {\bf I}} +  {\bf A}^{k+1} } \right)_{l.}
\left( {{\rm {\bf a}}_{l.}^{(k)} }  \right)
} \right|}
\end{array}
\end{equation}
 Having changed the order of summation,  we obtain on the other hand
for all $s = \overline {1, n - 1} $
\begin{equation}
\label{eq:char_drdet3}
\begin{array}{c}
  d_{s} = {\sum\limits_{\alpha \in I_{s,\,n} {\left\{ {j} \right\}}}
{\left| \left( {{\bf A}^{k+1} } \right)_{\alpha}
^{\alpha}  \right|} }  =
  \sum\limits_{\alpha \in I_{s,\,n} {\left\{ {j} \right\}}} \sum\limits_{l} \left|\left({ { {\bf A}^{k+1} _{j.}  \left( { a_{jl}{\rm {\bf a}}_{l.}^{(k)} } \right)} }\right)_{\alpha}
^{\alpha}
\right|=\\
  \sum\limits_{l} a_{jl}  \cdot \sum\limits_{\alpha \in I_{s,\,n} {\left\{ {j} \right\}}} \left|\left({ { {\bf A}^{k+1} _{j.}  \left( { {\rm {\bf a}}_{l.}^{(k)} } \right)} }\right)_{\alpha}
^{\alpha}  \right|
\end{array}
\end{equation}
By substituting (\ref{eq:char_drdet2}) and (\ref{eq:char_drdet3}) in
(\ref{eq:char_drdet1}), and equating factors at $a _ {jl} $ when
$l = i $, we obtain the equality (\ref{eq:char_drdet}). $\blacksquare$

\begin{theorem}\label{theor:dr_repr_row}
 If $Ind{\kern 1pt} {\rm {\bf A}} = k $ and $\rank{\rm {\bf
A}}^{k + 1} = \rank{\rm {\bf A}}^{k}=r \le n$ for  ${\rm {\bf A}}\in {\mathbb C}^{n\times n} $, then the
Drazin inverse  ${\rm {\bf A}}^{D} = \left( {a_{ij}^{ D} } \right)
\in {\rm {\mathbb{C}}}_{}^{n\times n} $ possess the following
determinantal representations:
\begin{equation}
\label{eq:dr_repr_row}
  a_{ij}^{ D} = \frac{
  {\sum\limits_{\alpha \in I_{r,n} {\left\{ {j} \right\}}}
{{\left| {\left( {{\rm {\bf A}}_{j.}^{k + 1} \left( {{\rm {\bf
a}}_{i.}^{\left( {k} \right)}}  \right)} \right)_{\alpha}
^{\alpha} } \right|}}}  }
{  {\sum\limits_{\alpha \in I_{r,n}} {{\left| {\left( {{\rm {\bf
A}}^{k + 1}} \right)_{\alpha} ^{\alpha }} \right|}}}},
[ \end{equation}
and
\begin{equation}
\label{eq:dr_repr_col}
  a_{ij}^{ D} = \frac{
  {\sum\limits_{\beta \in J_{r,n} {\left\{ {i} \right\}}}
{{\left| {\left( {{\rm {\bf A}}_{.\,i}^{k + 1} \left( {{\rm {\bf
a}}_{.j}^{\left( {k} \right)}}  \right)} \right)_{\beta} ^{\beta}
} \right|}}}  }{
  {\sum\limits_{\beta \in J_{r,n}} {{\left| {\left( {{\rm {\bf
A}}^{k + 1}} \right)_{\beta} ^{\beta }} \right|}}}},
 \end{equation}
for all $i,j = \overline {1,n}$.
\end{theorem}
{\textit{Proof.}} At first we shall prove the equation (\ref{eq:dr_repr_row}).

If $\lambda\in {\mathbb{R}}_{+}$, then  $\rank \left( {\lambda {\rm {\bf I}} + {\rm {\bf A}}^{ k+1} } \right)=n$.
 Hence,   there exists the  inverse matrix
\[
\left( {\lambda {\rm {\bf I}} + {\rm {\bf A}}^{ k+1} }
\right)^{ - 1} = {\frac{{1}}{{\det \left( {\lambda {\rm {\bf I}} +
{\rm {\bf A}}^{ k+1 }} \right)}}}\left(
{{\begin{array}{*{20}c}
 {R_{11}}  \hfill & {R_{21}}  \hfill & {\ldots}  \hfill & {R_{n\,1}}  \hfill
\\
 {R_{12}}  \hfill & {R_{22}}  \hfill & {\ldots}  \hfill & {R_{n\,2}}  \hfill
\\
 {\ldots}  \hfill & {\ldots}  \hfill & {\ldots}  \hfill & {\ldots}  \hfill
\\
 {R_{1\,n}}  \hfill & {R_{2\,n}}  \hfill & {\ldots}  \hfill & {R_{n\,n}}  \hfill
\\
\end{array}} } \right),
\]
 where  $R_{ij} $  is a cofactor
  in $\lambda {\rm {\bf
I}} + {\rm {\bf A}}^{ k+1}$ for all $ i,j=\overline{1,n}$.   By Theorem \ref{theor:lim_A_k_A_k+1},
${\rm {\bf A}}^{D} = {\mathop {\lim} \limits_{\lambda \to 0}}{\rm
{\bf A}}^{k} \left( {\lambda {\rm {\bf I}}_n + {\rm {\bf A}}^{k +
1}} \right)^{ - 1}$, so that
\[{\rm {\bf A}}^{ D}  ={\mathop {\lim} \limits_{\lambda \to 0}}{\frac{{1}}{{\det \left( {\lambda {\rm {\bf I}} +
{\rm {\bf A}}^{ k+1 }} \right)}}}
                                       \begin{pmatrix}
                                         \sum_{s=1}^{n} a_{1s}^{(k)} R_{1s} & \ldots & \sum_{s=1}^{n} a_{1s}^{(k)} R_{ns}  \\
                                         \ldots &  \ldots &  \ldots \\
                                         \sum_{s=1}^{n} a_{ns}^{(k)} R_{1s} &  \ldots & \sum_{s=1}^{n} a_{ns}^{(k)} R_{ns} \\  \end{pmatrix}=\]
\begin{equation}
\label{lim_A_D}
{\mathop {\lim} \limits_{\lambda \to 0}}\begin{pmatrix}
 \frac{\det \left( {\lambda {\rm {\bf I}} + {\rm {\bf A}}^{ k+1}}
\right)_{1.} \left( {{\rm {\bf a}}_{1.}^{ (k)} }  \right)}{{\det
\left( {\lambda {\rm {\bf I}} + {\rm {\bf A}}^{ k+1} }
\right)}}
 & \ldots & \frac{\det \left( {\lambda {\rm {\bf I}} + {\rm {\bf
A}}^{k+1} } \right)_{n.} \left( {{\rm {\bf
a}}_{n.}^{(k)} }  \right)}{{\det \left( {\lambda {\rm {\bf I}} +
{\rm {\bf
A}}^{ k+1} } \right)}}\\
  \ldots & \ldots & \ldots \\
\frac{\det \left( {\lambda {\rm {\bf I}} + {\rm {\bf A}}^{ k+1}} \right)_{1.} \left( {{\rm {\bf a}}_{n.}^{ (k)} }
\right)}{{\det \left( {\lambda {\rm {\bf I}} + {\rm {\bf A}}^{ k+1}
} \right)}} &
 \ldots &
 \frac{\det \left( {\lambda {\rm {\bf I}} + {\rm {\bf A}}^{ k+1} } \right)_{n.} \left( {{\rm {\bf a}}_{n.}^{ (k)} }
\right)}{{\det \left( {\lambda {\rm {\bf I}} + {\rm {\bf A}}^{k+1}
} \right)}}
\end{pmatrix}
\end{equation}
Taking into account Theorem
\ref{th:char_pol} , we have
\[ \det \left( {\lambda{\rm {\bf I}} + {\rm {\bf A}}^{k+1}
} \right) = \lambda ^{n} + d_{1} \lambda ^{n - 1} +
d_{2} \lambda ^{n - 2} + \ldots + d_{n},
\]
 where $d_{s}={  {\sum\limits_{\alpha \in I_{s,n}} {{\left| {\left( {{\rm {\bf
A}}^{k + 1}} \right)_{\alpha} ^{\alpha }} \right|}}}}$  is a
sum of the principal minors  of ${\rm {\bf A}}^{ k+1} $
of order $s$,  for all $ s=\overline{1,n-1}$, and $d_{n}=\det  {\rm {\bf A}}^{k+1} $.
Since $\rank{\rm {\bf A}}^{ k+1} = r$, then $d_{n} = d_{n - 1} = \ldots = d_{r + 1} = 0$ and
\begin{equation}
\label{eq_detA}\det \left( {\lambda {\rm {\bf I}} + {\rm {\bf A}}^{k+1} } \right)
= \lambda ^{n} + d_{1} \lambda ^{n - 1} + d_{2} \lambda ^{n - 2} + \ldots + d_{r} \lambda
^{n - r}.
 \end{equation}
By Lemma \ref{lem:char_pol_dr}  for  all $i,j=\overline{1,n} $,
\[ \det \left( {\lambda {\rm {\bf I}} + {\rm {\bf
A}}^{k+1}} \right)_{j.} \left( {{\rm {\bf
a}}_{i.}^{ (k)} } \right) = l_{1}^{\left( {ij} \right)} \lambda ^{n
- 1} + l_{2}^{\left( {ij} \right)} \lambda ^{n - 2} + \ldots +
l_{n}^{\left( {ij} \right)}, \]
 where for  all $s=\overline{1,n-1} $,
 \[l_{s}^{\left( {ij}
\right)} = {\sum\limits_{\alpha \in I_{s,n} {\left\{ {j} \right\}}}
{{\left| {\left( {{\rm {\bf A}}_{j\,.}^{k + 1} \left( {{\rm {\bf
a}}_{i.}^{\left( {k} \right)}}  \right)} \right)_{\alpha}
^{\alpha} } \right|}}},\] and \,\,$l_{n}^{\left( {i\,j} \right)} = \det
{\rm {\bf A}}^{ k+1}_{j\,.} \left(
{{\rm {\bf a}}_{i.}^{\left( {k} \right)} } \right)$.

By Lemma \ref{lem:rank_A_row},
$ \rank{\rm {\bf A}}_{j \,.}^{k + 1}
\left( {{\rm {\bf a}}_{i\,.}^{\left( {k} \right)}}  \right) \le r$, so that if
$s>r$, then for all $\alpha \in I_{s,n} {\left\{
{i} \right\}}$ and for all $i,j = \overline {1,n}$,
\[{{\left| {\left( {{\rm {\bf A}}_{j\,.}^{k + 1} \left( {{\rm {\bf
a}}_{i.}^{\left( {k} \right)}}  \right)} \right)_{\alpha}
^{\alpha} } \right|}}= 0.\]

 Therefore if $r + 1 \le s < n$, then for all $
i,j = \overline {1,n}$,
\[l_{s}^{\left( {ij}
\right)} =  {\sum\limits_{\alpha \in I_{s,n} {\left\{ {j} \right\}}}
{{\left| {\left( {{\rm {\bf A}}_{j.}^{k + 1} \left( {{\rm {\bf
a}}_{i.}^{\left( {k} \right)}}  \right)} \right)_{\alpha}
^{\alpha} } \right|}}}= 0,\]
and $l_{n}^{\left( {ij} \right)} =
\det {\rm {\bf A}}^{ k+1}_{j\,.}
\left( {{\rm {\bf a}}_{i.}^{ (k)} } \right) = 0$.
Finally we obtain
\begin{equation}\label{eq_det_A_j}
  \det \left( {\lambda {\rm {\bf I}} + {\rm {\bf A}}^{k+1} } \right)_{j.} \left( {{\rm {\bf a}}_{i.}^{(k)} } \right) =
l_{1}^{\left( {i\,j} \right)} \lambda ^{n - 1} + l_{2}^{\left(
{i\,j} \right)} \lambda ^{n - 2} + \ldots + l_{r}^{\left( {ij}
\right)} \lambda ^{n - r}.
\end{equation}

 By replacing the denominators and the nominators  of the fractions in the entries of the matrix (\ref{lim_A_D}) with  the
 expressions (\ref{eq_detA}) and (\ref{eq_det_A_j}) respectively, finally we obtain
\[
   {\rm {\bf A}}^{D}  = {\mathop {\lim} \limits_{\lambda \to 0}}
\begin{pmatrix}
 {{\frac{{l_{1}^{\left( {11} \right)} \lambda ^{n - 1} + \ldots +
l_{r}^{\left( {11} \right)} \lambda^{n - r}}}{{\lambda ^{n} +
d_{1} \lambda ^{n - 1} + \ldots + d_{r} \lambda ^{n - r}}}}}
 & \ldots & {{\frac{{l_{1}^{\left( {1n} \right)} \lambda ^{n
- 1} + \ldots + l_{r}^{\left( {1n} \right)} \lambda ^{n -
r}}}{{\lambda ^{n} + d_{1} \lambda
^{n - 1} + \ldots + d_{r} \lambda ^{n - r}}}}} \\
  \ldots & \ldots & \ldots \\
  {{\frac{{l_{1}^{\left( {n1} \right)} \lambda ^{n - 1} + \ldots +
l_{r}^{\left( {n1} \right)} \lambda ^{n - r}}}{{\lambda ^{n} +
d_{1} \lambda ^{n - 1} + \ldots + d_{r} \lambda ^{n - r}}}}} &
\ldots & {{\frac{{l_{1}^{\left( {nn} \right)} \lambda ^{n - 1} +
\ldots + l_{r}^{\left( {nn} \right)} \lambda ^{n - r}}}{{\lambda
^{n} + d_{1} \lambda ^{n - 1} + \ldots + d_{r} \lambda ^{n -
r}}}}}
\end{pmatrix}=\]
  \[ = \left( {{\begin{array}{*{20}c}
 {{\frac{{l_{r}^{\left( {11} \right)}} }{{d_{r}} }}} \hfill & {\ldots}
\hfill & {{\frac{{l_{r}^{\left( {1n} \right)}} }{{d_{r}} }}} \hfill \\
 {\ldots}  \hfill & {\ldots}  \hfill & {\ldots}  \hfill \\
 {{\frac{{l_{r}^{\left( {n1} \right)}} }{{d_{r}} }}} \hfill & {\ldots}
\hfill & {{\frac{{l_{r}^{\left( {nn} \right)}} }{{d_{r}} }}} \hfill \\
\end{array}} } \right),
\]
where for all $i,j=\overline{1,n} $,
 \[l_{r}^{\left( {ij}
\right)} = {\sum\limits_{\alpha \in I_{r,n} {\left\{ {j} \right\}}}
{{\left| {\left( {{\rm {\bf A}}_{j.}^{k + 1} \left( {{\rm {\bf
a}}_{i.}^{\left( {k} \right)}}  \right)} \right)_{\alpha}
^{\alpha} } \right|}}},\,\,\,d_{r}={  {\sum\limits_{\alpha \in I_{r,n}} {{\left| {\left( {{\rm {\bf
A}}^{k + 1}} \right)_{\alpha} ^{\alpha }} \right|}}}}. \]
The equation (\ref{eq:dr_repr_col}) can be proved similarly.

This completes the proof.
$\blacksquare$

Using Theorem \ref{theor:dr_repr_row} we evidently can obtain  determinantal representations of the group inverse and  the following determinantal representation of the  identities ${\bf A}^{D}{\bf A}$ and ${\bf A}{\bf A}^{D}$   on  $R  ({\bf A}^{k})$
 \begin{corollary}
 If $Ind{\kern 1pt} {\rm {\bf A}} = 1 $ and $\rank{\rm {\bf
A}}^{2} = \rank{\rm {\bf A}}=r \le n$ for
${\rm {\bf A}}\in {\mathbb C}^{n\times n} $, then the
group inverse  ${\rm {\bf A}}^{g} = \left( {a_{ij}^{ g} } \right)
\in {\rm {\mathbb{C}}}_{}^{n\times n} $ possess the following
determinantal representations:
\begin{equation}\label{eq:AAg}
 a_{ij}^{g} ={\frac{\sum\limits_{\alpha \in I_{r,n} {\left\{ {j} \right\}}} {{\left|
{\left( {{\rm {\bf A}}_{j.}^{2} \left( {{\rm {\bf a}}_{i.}}
\right)} \right)_{\alpha} ^{\alpha} } \right|}}}{{\sum\limits_{\alpha \in I_{r,n}} {{\left| {\left( {{\rm {\bf
A}}^{2}} \right)_{\alpha} ^{\alpha }}  \right|}}}}},
 \end{equation}\[
 a_{ij}^{g} ={\frac{\sum\limits_{\beta \in J_{r,n} {\left\{ {i} \right\}}} {{\left|
{\left( {{\rm {\bf A}}_{.\,i}^{2} \left( {{\rm {\bf a}}_{.j}}
\right)} \right)_{\beta} ^{\beta} } \right|}}}{{\sum\limits_{\beta \in J_{r,n}} {{\left| {\left( {{\rm {\bf
A}}^{2}} \right)_{\beta} ^{\beta }}  \right|}}}}},
\]
 for all $i,j = \overline {1,n}. $
\end{corollary}

\begin{corollary}
If $Ind{\kern 1pt} {\rm {\bf A}} = k $ and $\rank{\rm {\bf A}}^{k
+ 1} = \rank{\rm {\bf A}}^{k}=r \le n$ for
${\rm {\bf A}}\in {\mathbb C}^{n\times n} $, then the matrix ${\bf A}{\bf A}^{D}= (q_{ij})\in {\mathbb C}^{n\times n}$ possess the following determinantal representation
\begin{equation}\label{eq:AAd}
  q_{ij}= \frac{
  {\sum\limits_{\alpha \in I_{r,n} {\left\{ {j} \right\}}}
{{\left| {\left( {{\rm {\bf A}}_{j.}^{k + 1} \left( {{\rm {\bf
a}}_{i.}^{\left( {k+1} \right)}}  \right)} \right)_{\beta}
^{\beta} } \right|}}}  }
{  {\sum\limits_{\alpha \in I_{r,n}} {{\left| {\left( {{\rm {\bf
A}}^{k + 1}} \right)_{\beta} ^{\beta }} \right|}}}},
 \end{equation}
for all $i,j = \overline {1,n}. $
\end{corollary}
\begin{corollary}
If $Ind{\kern 1pt} {\rm {\bf A}} = k $ and $\rank{\rm {\bf A}}^{k
+ 1} = \rank{\rm {\bf A}}^{k}=r \le n$ for
${\rm {\bf A}}\in {\mathbb C}^{n\times n} $, then the matrix ${\bf A}^{D}{\bf A}= (p_{ij})\in {\mathbb C}^{n\times n}$ possess the following determinantal representation
\begin{equation}\label{eq:AdA}
 p_{ij} =\frac{{\sum\limits_{\beta \in J_{r,n} {\left\{ {i} \right\}}} {{\left|
{\left( {{\rm {\bf A}}_{.\,i}^{k+1} \left( {{\rm {\bf
a}}_{.j}}^{(k+1)} \right)} \right)_{\beta} ^{\beta} } \right|}}}}{{\sum\limits_{\beta \in J_{r,n} } {{\left|
{\left( {{\rm {\bf A}}_{.\,i}^{k+1} } \right)_{\beta} ^{\beta} } \right|}}}},
\end{equation}
for all $i,j = \overline {1,n}. $
\end{corollary}

\subsection{Analogues of the classical adjoint matrix for the W-weighted Drazin
inverse}
Cline and Greville \cite{cl} extended the Drazin inverse of square matrix to
rectangular matrix and called it as \textbf{the weighted Drazin inverse} (WDI).
The
 W-weighted Drazin inverse of  ${\bf A}\in {\mathbb{C}}^{m\times n}$  with respect to the matrix ${\bf W}\in {\mathbb{C}}^{n\times m}$  is defined to be the
unique solution ${\bf X}\in {\mathbb{C}}^{m\times n}$  of the following three matrix equations:
\begin{equation}\label{eq:WDr_prop}
\begin{array}{l}
  1)\,\,  ({{\bf A}{\bf W}})^{k+1}{\rm
{\bf X}}{\bf W}=({\rm {\bf
         A}}{\bf W})^{k};\\
  2)\,\,{\rm {\bf X}}{\bf W}{\rm {\bf A}}{\bf W}{\rm {\bf X}}={\rm {\bf X}};\\
  3)\,\, {\rm
{\bf A}}{\bf W}{\rm {\bf X}}={\rm {\bf X}}{\bf W}{\rm {\bf A}},
\end{array}
\end{equation}
where $k= {\rm max}\{Ind({\bf A}{\bf W}), Ind({\bf W}{\bf A})\}$.
It is denoted by ${\rm {\bf X}}={\rm {\bf A}}_{d, W}$.
In particular, when  ${\bf A}\in {\mathbb{C}}^{m\times m}$ and ${\bf W}={\bf I}_{ m}$ ,
then ${\bf A}_{d, W}$ reduce to ${\bf A}^{D}$. If ${\bf A}\in {\mathbb{C}}^{m\times m}$ is non-singular square
matrix and ${\bf W}={\bf I}_{ m}$, then $Ind(A) = 0$ and ${\bf A}_{d, W}={\bf A}^{D}={\bf A}^{-1}$.

The properties of WDI can be found in (e.g.,\cite{wei1,wei2,wei4,ze}). We note the general algebraic structures of  the W-weighted Drazin inverse \cite{wei1}. Let for ${\bf A}\in {\mathbb{C}}^{m\times n}$ and ${\bf W}\in {\mathbb{C}}^{n\times m}$
exist ${\bf L}\in {\mathbb{C}}^{m\times m}$ and ${\bf Q}\in {\mathbb{C}}^{n\times n}$ such that
\[\begin{array}{cc}
            {\bf A}={\bf L}\left(
                   \begin{array}{cc}
                     {\bf A}_{11} & {\bf 0} \\
                     {\bf 0} & {\bf A}_{22} \\
                   \end{array}
                 \right){\bf Q}^{-1}, & {\bf W}={\bf Q}\left(
                   \begin{array}{cc}
                     {\bf W}_{11} & {\bf 0} \\
                     {\bf 0} & {\bf W}_{22} \\
                   \end{array}
                 \right){\bf L}^{-1}
          \end{array}
\]
Then \[{\bf A}_{d, W}={\bf L}\left(
                   \begin{array}{cc}
                    ( {\bf W}_{11}{\bf A}_{11}{\bf W}_{11})^{-1} & {\bf 0} \\
                     {\bf 0} &{\bf 0} \\
                   \end{array}
                 \right){\bf Q}^{-1},\]
where ${\bf L}$, ${\bf L}$, ${\bf A}_{11}$, ${\bf W}_{11}$ are non-singular matrices, and ${\bf A}_{22}$, ${\bf W}_{22}$ are nilpotent matrices.
By  \cite{ca} we have the following  limit representations of the W-weighted
Drazin inverse,
\begin{equation}\label{eq:lim_repr_AW}
{\rm {\bf A}}_{d,W} = {\mathop {\lim} \limits_{\lambda \to 0}}
\left( {\lambda {\rm {\bf I}}_m + ({\rm {\bf A}}{\bf W})^{k + 2}} \right)^{
- 1}({\rm {\bf A}}{\bf W})^{k}{\bf A}
\end{equation}
and
\begin{equation}\label{eq:lim_repr_WA}
{\rm {\bf A}}_{d,W} = {\mathop {\lim} \limits_{\lambda \to 0}}{\bf A}({\bf W}{\rm {\bf A}})^{k}
\left( {\lambda {\rm {\bf I}}_n + ({\bf W}{\rm {\bf A}})^{k + 2}} \right)^{
- 1}
\end{equation}
where  $\lambda \in {\mathbb R} _{ +}  $, and ${\mathbb R} _{ +} $
is a set of the real positive numbers.

Denote ${\bf W}{\bf A}=:{\bf U} $ and ${\bf A}{\bf W}=:{\bf V} $.
 Denote by ${\rm {\bf v}}_{.j}^{(k)} $ and ${\rm {\bf
v}}_{i.}^{(k)} $ the $j$th column  and the $i$th row of  ${\rm
{\bf V}}^{k} $ respectively. Denote by $ \bar{{\bf V}}^{k}:=({\bf A}{\bf W})^{k}{\bf A}\in {\mathbb{C}}^{m\times n}$ and $ \bar{{\bf W}}={\bf W}{\bf A}{\bf W}\in {\mathbb{C}}^{n\times m}$.

\begin{lemma} \label{lem:rank_wcol} If ${\bf A}{\bf W}={\bf
V}=\left(v_{ij} \right)\in {\mathbb{C}}^{m\times m}$ with $ Ind{\kern 1pt} {\rm {\bf V}}=k$, then
\begin{equation}\label{eq:rank_wcol}
 \rank\,\left( {{\rm {\bf V}}^{ k+2} }
\right)_{.\,i} \left( { \bar{{\bf v}}_{.j}^{ (k)} }  \right) \le
\rank\,\left( {{\rm {\bf V}}^{k+2} } \right).
\end{equation}
\end{lemma}

{\textit{Proof.}} We have ${\rm {\bf V}}^{
 k+2}=\bar{{\bf V}}^{ k}\bar{{\bf W}}$.
Let ${\rm {\bf P}}_{i\,s} \left( {-\bar{w}_{j\,s}}  \right)\in {\mathbb
C}^{m\times m} $, $(s \ne i )$, be a matrix with $-\bar{w}_{j\,s} $ in
the $(i, s)$ entry, 1 in all diagonal entries, and 0 in others.
The matrix ${\rm {\bf P}}_{i\,s} \left( {-\bar{w}_{j\,s}}  \right)$, $(s \ne i )$, is a matrix of an elementary transformation. It follows that
\[
\left( {\rm {\bf V}}^{ k+2} \right)_{.\,i} \left(
{ \bar{{\bf v}}_{.\,j}^{ (k)} }  \right) \cdot {\prod\limits_{s \ne
i} {{\rm {\bf P}}_{i\,s} \left( {-\bar{w}_{j\,s}}  \right) =
{\mathop
{\left( {{\begin{array}{*{20}c}
 {{\sum\limits_{s \ne j} {\bar{v}_{1s}^{ (k)}  \bar{w}_{s1}} } } \hfill & {\ldots}  \hfill
& {\bar{v}_{1j}^{ (k)} }  \hfill & {\ldots}  \hfill & {{\sum\limits_{s
\ne
j } {\bar{v}_{1s}^{ (k)}  \bar{w}_{sm}}}} \hfill \\
 {\ldots}  \hfill & {\ldots}  \hfill & {\ldots}  \hfill & {\ldots}  \hfill &
{\ldots}  \hfill \\
 {{\sum\limits_{s \ne j} {\bar{v}_{ms}^{(k)}  \bar{w}_{s1}} } } \hfill & {\ldots}  \hfill
& {\bar{v}_{mj}^{(k)} }  \hfill & {\ldots}  \hfill & {{\sum\limits_{s
\ne
j } {\bar{v}_{ms}^{(k)}  \bar{w}_{sm}}}} \hfill \\
\end{array}} }
\right)}\limits_{i-th}}}}.
\]
 We have the next factorization of the obtained matrix.
\[
{\mathop
{\left( {{\begin{array}{*{20}c}
 {{\sum\limits_{s \ne j} {\bar{v}_{1s}^{ (k)}  \bar{w}_{s1}} } } \hfill & {\ldots}  \hfill
& {\bar{v}_{1j}^{ (k)} }  \hfill & {\ldots}  \hfill & {{\sum\limits_{s
\ne
j } {\bar{v}_{1s}^{ (k)}  \bar{w}_{sm}}}} \hfill \\
 {\ldots}  \hfill & {\ldots}  \hfill & {\ldots}  \hfill & {\ldots}  \hfill &
{\ldots}  \hfill \\
 {{\sum\limits_{s \ne j} {\bar{v}_{ms}^{(k)}  \bar{w}_{s1}} } } \hfill & {\ldots}  \hfill
& {\bar{v}_{mj}^{(k)} }  \hfill & {\ldots}  \hfill & {{\sum\limits_{s
\ne
j } {\bar{v}_{ms}^{(k)}  \bar{w}_{sm}}}} \hfill \\
\end{array}} }
\right)}\limits_{i-th}} =
\]
\[
 = \left( {\begin{array}{*{20}c}
 {\bar{v}_{11}^{ (k)} }  \hfill & {\bar{v}_{12}^{ (k)} }  \hfill & {\ldots}  \hfill &
{\bar{v}_{1n}^{ (k)} }  \hfill \\
 {\bar{v}_{21}^{ (k)} }  \hfill & {\bar{v}_{22}^{ (k)} }  \hfill & {\ldots}  \hfill &
{\bar{v}_{2n}^{ (k)} }  \hfill \\
 {\ldots}  \hfill & {\ldots}  \hfill & {\ldots}  \hfill & {\ldots}  \hfill
\\
 {\bar{v}_{m1}^{ (k)} }  \hfill & {\bar{v}_{m2}^{ (k)} }  \hfill & {\ldots}  \hfill &
{\bar{v}_{mn}^{ (k)} }  \hfill \\
\end{array}}  \right)
{\mathop {\left( {{\begin{array}{*{20}c}
 {\bar{w}_{11}}  \hfill & {\ldots}  \hfill & {0} \hfill & {\ldots}  \hfill &
{\bar{w}_{1m}}  \hfill \\
 {\ldots}  \hfill & {\ldots}  \hfill & {\ldots}  \hfill & {\ldots}  \hfill &
{\ldots}  \hfill \\
 {0} \hfill & {\ldots}  \hfill & {1} \hfill & {\ldots}  \hfill & {0} \hfill
\\
 {\ldots}  \hfill & {\ldots}  \hfill & {\ldots}  \hfill & {\ldots}  \hfill &
{\ldots}  \hfill \\
 {\bar{w}_{n1}}  \hfill & {\ldots}  \hfill & {0} \hfill & {\ldots}  \hfill &
{\bar{w}_{nm}}  \hfill \\
\end{array}} } \right)}\limits_{i-th}} j-th.
\]
 Denote ${\rm {\bf \tilde {W}}}: = {\mathop
{\left( {{\begin{array}{*{20}c}
 {\bar{w}_{11}}  \hfill & {\ldots}  \hfill & {0} \hfill & {\ldots}  \hfill &
{\bar{w}_{1m}}  \hfill \\
 {\ldots}  \hfill & {\ldots}  \hfill & {\ldots}  \hfill & {\ldots}  \hfill &
{\ldots}  \hfill \\
 {0} \hfill & {\ldots}  \hfill & {1} \hfill & {\ldots}  \hfill & {0} \hfill
\\
 {\ldots}  \hfill & {\ldots}  \hfill & {\ldots}  \hfill & {\ldots}  \hfill &
{\ldots}  \hfill \\
 {\bar{w}_{n1}}  \hfill & {\ldots}  \hfill & {0} \hfill & {\ldots}  \hfill &
{\bar{w}_{nm}}  \hfill \\
\end{array}} } \right)}\limits_{i-th}} j-th$. The matrix
${\rm {\bf \tilde {W}}}$ is obtained from $ \bar{{\bf W}}={\bf W}{\bf A}{\bf W}$ by
replacing all entries of the $j$th row  and the $i$th column with
zeroes except for 1 in the $(i, j)$ entry. Since elementary
transformations of a matrix do not change a rank, then $\rank{\rm
{\bf V}}^{ k+2} _{.\,i} \left( { \bar{{\bf v}}_{.j}^{ (k)} }
\right) \le \min {\left\{ {\rank \bar{{\bf V}}^{ k },\rank{\rm {\bf
\tilde {W}}}} \right\}}$. It is obvious that
  \[\begin{array}{c}
      \rank \bar{{\bf V}}^{ k }=\rank \,({\bf A}{\bf W})^{k}{\bf A}\geq  \rank\, ({\bf A}{\bf W})^{k+2}, \\
       \rank {\bf \tilde {W}}\geq \rank \,{\bf W}{\bf A}{\bf W}\geq  \rank\, ({\bf A}{\bf W})^{k+2}.
    \end{array}
     \]
 From this the inequality (\ref{eq:rank_wcol})
follows immediately.
$\blacksquare$

The next lemma is  proved similarly.
\begin{lemma} If ${\bf W}{\bf A}={\bf
U}=\left(u_{ij} \right)\in {\mathbb{C}}^{n\times n}$ with $ Ind{\kern 1pt} {\rm {\bf U}}=k$, then \[ \rank\left( {{\rm
{\bf U}}^{ k+2} } \right)_{i\,.} \left( { \bar{{\bf u}}_{j\,.}^{
(k)} }  \right) \le \rank\left( {{\rm {\bf U}}^{k+2} } \right),\]
where $ \bar{{\bf U}}^{k}:={\bf A}({\bf W}{\bf A})^{k}\in {\mathbb{C}}^{m\times n}$
\end{lemma}

Analogues of the characteristic polynomial are considered in the following two lemmas.
\begin{lemma}\label{lem:char_wcdet}
If ${\bf A}{\bf W}={\bf
V}=\left(v_{ij} \right)\in {\mathbb{C}}^{m\times m}$ with $ Ind{\kern 1pt} {\rm {\bf V}}=k$ and $\lambda \in \mathbb{R}$, then

\begin{equation}
\label{eq:char_wcdet} \left|  \left(\lambda {\rm {\bf I}}_m + {{\rm {\bf V}}^{ k+2} }
\right)_{.\,i} \left( { \bar{{\bf v}}_{.j}^{ (k)} }  \right)\right| = c_{1}^{\left( {ij} \right)} \lambda^{m - 1}
+ c_{2}^{\left( {ij} \right)} \lambda^{m - 2} + \ldots + c_{m}^{\left(
{ij} \right)},
\end{equation}

\noindent where $c_{m}^{\left( {ij} \right)} = {\det}
\left( {{\rm {\bf V}}^{ k+2} }
\right)_{.\,i} \left( { \bar{{\bf v}}_{.j}^{ (k)} }  \right)$  and $c_{s}^{\left( {ij} \right)} =
{\sum\limits_{\beta \in J_{s,\,m} {\left\{ {i} \right\}}}
{\det  \left( {\left( {{\rm {\bf V}}^{ k+2} }
\right)_{.\,i} \left( { \bar{{\bf v}}_{.j}^{ (k)} }  \right)}
\right){\kern 1pt}  _{\beta} ^{\beta} } }$ for all $s = \overline
{1, m- 1} $, $i = \overline {1,m}$, and $j = \overline {1,n}$.
\end{lemma}
{\textit{Proof.}}
Consider the  matrix $\left( {\lambda{\rm {\bf I}} + {\rm {\bf
V}}^{ k+2} } \right)_{. {\kern 1pt} i} ({\rm {\bf
v}}^{(k+2)}_{.{\kern 1pt}  i} ) \in {\rm {\mathbb{C}}}^{m\times
m}$.  Taking into account Theorem
\ref{th:char_pol} we obtain
\begin{equation}
\label{eq:char_wcdet1}  \left|\left( {\lambda{\rm {\bf I}} + {\rm {\bf
V}}^{ k+2} }
\right)_{.{\kern 1pt} i} \left( {{\rm {\bf
v}}^{(k+2)}_{.{\kern 1pt}
{\kern 1pt} i}}  \right)\right| = d_{1} \lambda^{m - 1} + d_{2} \lambda^{m - 2} +
\ldots + d_{m},
\end{equation}
where $d_{s} = {\sum\limits_{\beta \in J_{s,\,m} {\left\{ {i}
\right\}}} {| \left( {{\rm {\bf V}}^{ k+2} } \right) _{\beta} ^{\beta} }| } $ is the sum of all principal
minors of order $s$ that contain the $i$-th column for all $s =
\overline {1, m - 1} $ and $d_{m} = \det \left( {{\rm {\bf V}}^{
k+2} } \right)$.
Since ${\rm {\bf v}}^{(k+2)}_{.
{\kern 1pt} i} = \left( {{\begin{array}{*{20}c}
 {{\sum\limits_{l} {\bar{v}_{1l}^{(k)}  \bar{w}_{li}} } } \hfill \\
 {{\sum\limits_{l} {\bar{v}_{2l}^{(k)}  \bar{w}_{li}} } } \hfill \\
 { \vdots}  \hfill \\
 {{\sum\limits_{l} {\bar{v}_{nl}^{(k)}  \bar{w}_{li}} } } \hfill \\
\end{array}} } \right) = {\sum\limits_{l} {{\rm {\bf \bar{v}}}_{.\,l}^{(k)}  \bar{w}_{li}}
}$, where ${\rm {\bf \bar{v}}}_{.{\kern 1pt}  l}^{(k)}  $ is
the $l$th column-vector  of ${\rm {\bf \bar{V}}}^{k}=({\bf A}{\bf W})^{k}{\bf A}$  and ${\bf W}{\bf A}{\bf W}={\bf \bar{W}}=(\bar{w}_{li})$ for all
$l=\overline{1,n}$, then we have on the one hand
\begin{equation}
\label{eq:char_wcdet2}
\begin{array}{c}
  \left| \left( {\lambda{\rm {\bf I}} + {\rm {\bf V}}^{ k+2} }
\right)_{.{\kern 1pt} i} \left( {\rm {\bf v}}^{(k+2)}_{.
{\kern 1pt} i} \right)\right|
   = \sum\limits_{l} \left|  \left( {\lambda{\rm {\bf I}} + {\rm {\bf
   V}}^{k+2
}} \right)_{.{\kern 1pt} l} \left( {{\rm {\bf \bar{v}}}_{.\,l}^{(k)}  \bar{w}_{li}} \right) \right|= \\ {\sum\limits_{l}
\left| {\left( {\lambda{\rm {\bf I}} + {\rm {\bf V}}^{ k+2} }
\right)_{.{\kern 1pt} i} \left( {{\rm {\bf \bar{v}}}_{.{\kern 1pt}
{\kern 1pt} l}^{(k)} }  \right)  } \right|} \cdot \bar{w}_{li}
\end{array}
\end{equation}
 Having changed the order of summation,  we obtain on the other hand
for all $s = \overline {1, m - 1} $
\begin{equation}
\label{eq:char_wcdet3}
\begin{array}{c}
  d_{s} = {\sum\limits_{\beta \in J_{s,\,m} {\left\{ {i} \right\}}}
{\left| \left( {{\rm {\bf V}}^{ k+2} } \right) {\kern
1pt} _{\beta} ^{\beta} \right|} }  =
  \sum\limits_{\beta \in J_{s,\,m}
{\left\{ {i} \right\}}} \sum\limits_{l} \left|{ \left(
{\left( {{\rm {\bf V}}^{ k+2} {\kern 1pt}} \right)_{.\, i} \left(
{{\rm {\bf \bar{v}}}_{.\, l}^{(k)} \bar{w}_{l\,i}} \right)} \right)}  \,
_{\beta} ^{\beta}\right|=\\
  \sum\limits_{l}  \sum\limits_{\beta \in J_{s,\,m} {\left\{ {i}
\right\}}} \left| \left( {\left( {{\rm {\bf V}}^{ k+2}
}  \right)_{.{\kern 1pt} i} \left( {{\rm {\bf
\bar{v}}}_{.{\kern 1pt}  l}^{(k)} }  \right)} \right){\kern
1pt} _{\beta} ^{\beta}    \right| \cdot \bar{w}_{l{\kern 1pt} i}.
\end{array}
\end{equation}
By substituting (\ref{eq:char_wcdet2}) and (\ref{eq:char_wcdet3}) in (\ref{eq:char_wcdet1}),
and equating factors at $\bar{w} _ {l \, i} $ when $l = j $, we obtain
the equality (\ref{eq:char_wcdet}).
$\blacksquare$

 By analogy can be proved
the following lemma.
\begin{lemma}
If ${\bf W}{\bf A}={\bf
U}=\left(u_{ij} \right)\in {\mathbb{C}}^{n\times n}$ with $ Ind{\kern 1pt} {\rm {\bf U}}=k$ and $\lambda \in \mathbb{R}$, then
\[ \left|  {( \lambda{\rm {\bf I}} + {\rm {\bf U}}^{k+2} )_{j\,.\,} ({\rm {\bf \bar{u}}}_{i.}^{(k)}
)}\right|
   = r_{1}^{\left( {ij} \right)} \lambda^{n
- 1} +r_{2}^{\left( {ij} \right)} \lambda^{n - 2} + \ldots +
r_{n}^{\left( {ij} \right)},
\]

\noindent where  $r_{n}^{\left( {ij} \right)} = \left|
{({\rm {\bf U}}^{k+2} )_{j\,.\,} ({\rm {\bf \bar{u}}}_{i.\,}^{(k)} )}\right|$
and $r_{s}^{\left( {ij} \right)} = \sum\limits_{\alpha \in
I_{s,n} {\left\{ {j} \right\}}} \left| \left( {({\rm
{\bf U}}^{k+2} )_{j\,.\,} ({\rm {\bf \bar{u}}}_{i.\,}^{(k)} )}
\right)\,_{\alpha} ^{\alpha}  \right|$ for all $s = \overline {1,n -
1} $,  $i = \overline {1,m}$, and $j = \overline {1,n}$.
\end{lemma}

\begin{theorem}\label{theor:det_rep_wdraz}
If ${\rm {\bf A}} \in  {\mathbb{C}}^{m\times n}$, ${\rm {\bf W}} \in  {\mathbb{C}}^{n\times m}$ with
$k= {\rm max}\{Ind({\bf A}{\bf W}), Ind({\bf W}{\bf A})\}$ and $\rank({\rm {\bf A}}{\bf W})^{ k} = r$, then the W-weighted Drazin inverse ${\rm {\bf
A}}_{d,W} = \left( {a_{ij}^{d,W} } \right) \in {\rm
{\mathbb{C}}}^{m\times n} $ with respect to ${\bf W}$ possess the following determinantal
representations:
\begin{equation}
\label{eq:dr_rep_wcdet} a_{ij}^{d,W}  = {\frac{\sum\limits_{\beta
\in J_{r,\,m} {\left\{ {i} \right\}}} \left| \left(
{\left( {{\rm {\bf A}}{\bf W}} \right)^{k+2}_{\,. \,i} \left( {{\rm {\bf
\bar{v}}}_{.j}^{ (k)} }  \right)} \right){\kern 1pt} {\kern 1pt} _{\beta}
^{\beta} \right|  }{{{\sum\limits_{\beta \in J_{r,\,\,m}} {{\left|
{\left( { {\bf A}{\bf W}} \right)^{k+2}{\kern 1pt} _{\beta} ^{\beta}
}  \right|}}} }}},
\end{equation}
or
\begin{equation}
\label{eq:dr_rep_wrdet} a_{ij}^{d,W}  = {\frac{\sum\limits_{\alpha
\in I_{r,n} {\left\{ {j} \right\}}} \left| \left(
{({\bf W}{\rm {\bf A}} )^{ k+2}_{j\,.\,} ({\rm {\bf \bar{u}}}_{i.\,}^{ (k)} )}
\right)\,_{\alpha} ^{\alpha} \right| }{{{\sum\limits_{\alpha \in
I_{r,\,n}}  {{\left| {\left({\bf W} {{\rm {\bf A}} } \right)^{k+2}{\kern
1pt}  _{\alpha} ^{\alpha} } \right|}}} }}}.
\end{equation}
where ${\rm {\bf \bar{v}}}_{.j}^{(k)} $ is the $j$th column of  ${\rm
{\bf \bar{V}}}^{k}=({\bf A}{\bf W})^{k}{\bf A} $  for all $j=1,...,m$ and ${\rm {\bf
\bar{u}}}_{i.}^{(k)} $ is the $i$th row of  ${\rm
{\bf \bar{U}}}^{k}={\bf A}({\bf W}{\bf A})^{k} $ for all $i=1,...,n$.

\end{theorem}
{\textit{Proof.}}   At first we shall prove (\ref{eq:dr_rep_wcdet}). By
(\ref{eq:lim_repr_AW}),
\[{\rm {\bf A}}_{d,W} = {\mathop {\lim} \limits_{\lambda \to 0}}
\left( {\lambda {\rm {\bf I}}_m + ({\rm {\bf A}}{\bf W})^{k + 2}} \right)^{
- 1}({\rm {\bf A}}{\bf W})^{k}{\bf A}.\]
Let
\[
\left( {\lambda {\rm {\bf I}}_m + ({\rm {\bf A}}{\bf W})^{k + 2}} \right)^{
- 1}
= {\frac{{1}}{{\det \left( {\lambda {\rm {\bf I}}_m + ({\rm {\bf A}}{\bf W})^{k + 2}} \right)}}}\left( {{\begin{array}{*{20}c}
 {L_{11}}  \hfill & {L_{21}}  \hfill & {\ldots}  \hfill & {L_{m1}}  \hfill
\\
 {L_{12}}  \hfill & {L_{22}}  \hfill & {\ldots}  \hfill & {L_{m2}}  \hfill
\\
 {\ldots}  \hfill & {\ldots}  \hfill & {\ldots}  \hfill & {\ldots}  \hfill
\\
 {L_{1m}}  \hfill & {L_{2m}}  \hfill & {\ldots}  \hfill & {L_{mm}}  \hfill
\\
\end{array}} } \right),
\]
\noindent where $L_{ij} $ is a left $ij$-th cofactor  of a matrix
${\lambda {\rm {\bf I}}_m + ({\rm {\bf A}}{\bf W})^{k + 2}}$. Then we have
\[\begin{array}{l}
   \left( {\lambda {\rm {\bf I}}_m + ({\rm {\bf A}}{\bf W})^{k + 2}} \right)^{
- 1}({\rm {\bf A}}{\bf W})^{k}{\bf A}  = \\
  ={\frac{{1}}{{\det \left( {\lambda {\rm {\bf I}}_m + ({\rm {\bf A}}{\bf W})^{k + 2}} \right)}}}\left(
{{\begin{array}{*{20}c}
 {{\sum\limits_{s = 1}^{m} {L_{s1} \bar{v}_{s1}^{ (k)} } } } \hfill &
{{\sum\limits_{s = 1}^{m} {L_{s1} \bar{v}_{s2}^{ (k)} } } } \hfill &
{\ldots}
\hfill & {{\sum\limits_{s = 1}^{m} {L_{s1} \bar{v}_{sn}^{ (k)} } } } \hfill \\
 {{\sum\limits_{s = 1}^{m} {L_{s2} \bar{v}_{s1}^{ (k)} } } } \hfill &
{{\sum\limits_{s = 1}^{m} {L_{s2} \bar{v}_{s2}^{ (k)} } } } \hfill &
{\ldots}
\hfill & {{\sum\limits_{s = 1}^{m} {L_{s2} \bar{v}_{sn}^{ (k)} } } } \hfill \\
 {\ldots}  \hfill & {\ldots}  \hfill & {\ldots}  \hfill & {\ldots}  \hfill
\\
 {{\sum\limits_{s = 1}^{m} {L_{sm} \bar{v}_{s1}^{ (k)} } } } \hfill &
{{\sum\limits_{s = 1}^{m} {L_{sm} \bar{v}_{s2}^{ (k)} } } } \hfill &
{\ldots}
\hfill & {{\sum\limits_{s = 1}^{m} {L_{sm} \bar{v}_{sn}^{ (k)} } } } \hfill \\
\end{array}} } \right).
\end{array}
\]
By (\ref{eq:lim_repr_AW}), we obtain
\begin{equation}
\label{eq:A_d1} {\rm {\bf A}}_{d,W}  = {\mathop {\lim}
\limits_{\lambda \to 0}} \left( {{\begin{array}{*{20}c}
 {{\frac{\left| \left( {\lambda {\rm {\bf I}}_m + ({\rm {\bf A}}{\bf W})^{k + 2}} \right)_{.1} \left( {{\rm {\bf \bar{v}}}_{.1}^{ (k)} }
\right)\right|}{{\left| \left( {\lambda {\rm {\bf I}}_m + ({\rm {\bf A}}{\bf W})^{k + 2}} \right)\right|}}}} \hfill & {\ldots}  \hfill & {{\frac{{\left| \left( {\lambda {\rm {\bf I}}_m + ({\rm {\bf A}}{\bf W})^{k + 2}} \right)_{.1} \left( {{\rm {\bf \bar{v}}}_{.n}^{(k)} } \right)\right|}}{{\left|
\left( {\lambda {\rm {\bf I}}_m + ({\rm {\bf A}}{\bf W})^{k + 2}} \right)\right|}}}} \hfill \\
 {\ldots}  \hfill & {\ldots}  \hfill & {\ldots}  \hfill \\
 {{\frac{{\left| \left( {\lambda {\rm {\bf I}}_m + ({\rm {\bf A}}{\bf W})^{k + 2}} \right)_{.n} \left( {{\rm {\bf \bar{v}}}_{.1}^{ (k)} }
\right)\right|}}{{\left| \left( {\lambda {\rm {\bf I}}_m + ({\rm {\bf A}}{\bf W})^{k + 2}} \right)\right|}}}} \hfill & {\ldots}  \hfill & {{\frac{{\left| \left( {\lambda {\rm {\bf I}}_m + ({\rm {\bf A}}{\bf W})^{k + 2}} \right)_{.m} \left( {{\rm {\bf \bar{v}}}_{.n}^{(k)} } \right)\right|}}{{\left|
\left( {\lambda {\rm {\bf I}}_m + ({\rm {\bf A}}{\bf W})^{k + 2}} \right)\right|}}}} \hfill \\
\end{array}} } \right).
\end{equation}
By Theorem \ref{th:char_pol} we have
 \[\left| \left( {\lambda {\rm {\bf I}}_m + ({\rm {\bf A}}{\bf W})^{k + 2}} \right)\right| = \lambda ^{m} + d_{1}
\lambda ^{m - 1} + d_{2} \lambda ^{m - 2} + \ldots + d_{m},\]
where
$d_{s} = {\sum\limits_{\beta \in J_{s,\,m}}  {{\left| {\left( {\lambda {\rm {\bf I}}_m + ({\rm {\bf A}}{\bf W})^{k + 2}} \right){\kern 1pt} {\kern 1pt} _{\beta}
^{\beta} } \right|}}} $ is a sum of principal minors of $({\rm {\bf A}}{\bf W})^{k + 2}$ of order $s$ for all  $s = \overline {1,m - 1} $ and
$d_{m} = \left|( {\rm {\bf A}}{\bf W})^{k + 2}\right|$.

Since $\rank({\rm {\bf A}}{\bf W})^{k + 2}=\rank({\rm {\bf A}}{\bf W})^{k + 1}=\rank({\rm {\bf A}}{\bf W})^{k} = r$,
then $d_{m} = d_{m - 1} = \ldots = d_{r + 1} = 0$. It follows that
$\det \left( {\lambda {\rm {\bf I}}_m + ({\rm {\bf A}}{\bf W})^{k + 2}} \right)
= \lambda ^{m} + d_{1} \lambda ^{m - 1} + d_{2}\lambda ^{m - 2} +
\ldots + d_{r} \lambda ^{m - r}$.

By Lemma \ref{lem:char_wcdet}
 \[\left| \left(\lambda {\rm {\bf I}}_m + {({\rm {\bf A}}{\bf W})^{ k+2} }
\right)_{.\,i} \left( { \bar{{\bf v}}_{.j}^{ (k)} }  \right)\right| = c_{1}^{\left( {ij} \right)} \lambda^{m - 1}
+ c_{2}^{\left( {ij} \right)} \lambda^{m - 2} + \ldots + c_{m}^{\left(
{ij} \right)} \]
for  $i = \overline {1,m} $ and $j = \overline {1,n} $, where $c_{s}^{\left( {ij}
\right)} = \sum\limits_{\beta \in J_{s,\,m} {\left\{ {i}
\right\}}} {\left| \left( {({\rm {\bf A}}{\bf W})^{ k+2} _{.\,i}
\left( {{\rm {\bf \bar{v}}}_{.j}^{ (k)} } \right)} \right)
{\kern 1pt} _{\beta }^{\beta}  \right|} $ for all $s = \overline {1,m -
1} $ and $c_{m}^{\left( {ij} \right)} = \left| ({\rm {\bf A}}{\bf W})^{ k+2} _{.i} \left( {{\rm {\bf \bar{v}}}_{.j}^{ (k)}
} \right)\right|$.

We shall prove that $c_{k}^{\left( {ij} \right)} = 0$, when $k \ge r +
1$ for   $i = \overline {1,m} $ and $j = \overline {1,n} $.
By Lemma \ref{lem:rank_wcol}
$\left( {({\rm {\bf A}}{\bf W})^{ k+2} _{.\,i}
\left( {{\rm {\bf \bar{v}}}_{.j}^{ (k)} } \right)} \right) \le r$, then the matrix
$\left( {({\rm {\bf A}}{\bf W})^{ k+2} _{.\,i}
\left( {{\rm {\bf \bar{v}}}_{.j}^{ (k)} } \right)} \right)$ has no more $r$ linearly
independent columns.

Consider $\left( {({\rm {\bf A}}{\bf W})^{ k+2} _{.\,i}
\left( {{\rm {\bf \bar{v}}}_{.j}^{ (k)} } \right)} \right) {\kern 1pt} _{\beta}
^{\beta}  $, when $\beta \in J_{s,m} {\left\{ {i} \right\}}$. It
is a principal submatrix of  $\left( {({\rm {\bf A}}{\bf W})^{ k+2} _{.\,i}
\left( {{\rm {\bf \bar{v}}}_{.j}^{ (k)} } \right)} \right)$ of order
$s \ge r + 1$. Deleting both its $i$-th row and column, we obtain
a principal submatrix of order $s - 1$ of  $({\rm {\bf A}}{\bf W})^{ k+2}$.
 We denote it by ${\rm {\bf M}}$. The following cases are
possible.
\begin{itemize}
  \item Let $s = r + 1$ and $\det {\rm {\bf M}} \ne 0$. In this case all
columns of $ {\rm {\bf M}} $ are right-linearly independent. The
addition of all of them on one coordinate to columns of
 $\left( {({\rm {\bf A}}{\bf W})^{ k+2} _{.\,i}
\left( {{\rm {\bf \bar{v}}}_{.j}^{ (k)} } \right)} \right)
{\kern 1pt} _{\beta} ^{\beta}$ keeps their right-linear
independence. Hence, they are basis in a matrix $\left( {({\rm {\bf A}}{\bf W})^{ k+2} _{.\,i}
\left( {{\rm {\bf \bar{v}}}_{.j}^{ (k)} } \right)} \right) {\kern 1pt} _{\beta}
^{\beta} $, and  the $i$-th
column is the right linear combination of its basis columns. From
this, $\left| \left( {({\rm {\bf A}}{\bf W})^{ k+2} _{.\,i}
\left( {{\rm {\bf \bar{v}}}_{.j}^{ (k)} } \right)} \right) {\kern 1pt}
_{\beta} ^{\beta} \right| = 0$, when $\beta \in J_{s,n} {\left\{ {i}
\right\}}$ and $s= r + 1$.
  \item If $s = r + 1$ and $\det {\rm {\bf M}} = 0$, than $p$, ($p \leq r$),
columns are basis in  ${\rm {\bf M}}$ and in  $\left( {({\rm {\bf A}}{\bf W})^{ k+2} _{.\,i}
\left( {{\rm {\bf \bar{v}}}_{.j}^{ (k)} } \right)} \right) {\kern 1pt} _{\beta} ^{\beta} $.
Then  $\left| \left( {({\rm {\bf A}}{\bf W})^{ k+2} _{.\,i}
\left( {{\rm {\bf \bar{v}}}_{.j}^{ (k)} } \right)} \right){\kern 1pt} _{\beta}
^{\beta} \right| = 0$ as well.
  \item If $s > r + 1$, then  $\det {\rm {\bf M}} = 0$ and
$p$, ($p < r$), columns are basis in the  both matrices ${\rm {\bf
M}}$ and $\left( {({\rm {\bf A}}{\bf W})^{ k+2} _{.\,i}
\left( {{\rm {\bf \bar{v}}}_{.j}^{ (k)} } \right)} \right)
{\kern 1pt} _{\beta} ^{\beta}  $. Therefore, $\left| \left( {({\rm {\bf A}}{\bf W})^{ k+2} _{.\,i}
\left( {{\rm {\bf \bar{v}}}_{.j}^{ (k)} } \right)} \right) {\kern 1pt} _{\beta} ^{\beta} \right| = 0.$
\end{itemize}
Thus in all cases we have $\left| \left( {({\rm {\bf A}}{\bf W})^{ k+2} _{.\,i}
\left( {{\rm {\bf \bar{v}}}_{.j}^{ (k)} } \right)} \right) {\kern 1pt} _{\beta} ^{\beta}\right|  = 0$,
when $\beta \in J_{s,m} {\left\{ {i} \right\}}$ and $r + 1 \le s <
m$. From here if $r + 1 \le s < m$, then
\[c_{s}^{\left( {ij} \right)} = \sum\limits_{\beta \in J_{s,\,m}
{\left\{ {i} \right\}}} {\left| \left( {({\rm {\bf A}}{\bf W})^{ k+2} _{.\,i}
\left( {{\rm {\bf \bar{v}}}_{.j}^{ (k)} } \right)} \right) {\kern 1pt} _{\beta} ^{\beta}\right| } =0,\] and
$c_{m}^{\left( {ij} \right)} = \det \left( {({\rm {\bf A}}{\bf W})^{ k+2} _{.\,i}
\left( {{\rm {\bf \bar{v}}}_{.j}^{ (k)} } \right)} \right) = 0$ for  $i = \overline {1,m} $ and $j = \overline {1,n} $.

Hence, $\left|\left( {\lambda {\rm {\bf I}} + ({\rm {\bf A}}{\bf W})^{ k+2}} \right)_{.\,i} \left( {{\rm {\bf \bar{v}}}_{.\,j}^{ (k)} }
\right)\right| =c_{1}^{\left( {ij} \right)} \lambda ^{m - 1} +
 \ldots +
c_{r}^{\left( {ij} \right)} \lambda ^{m - r}$ for  $i = \overline {1,m} $ and $j = \overline {1,n} $. By substituting these values in the matrix from
(\ref{eq:A_d1}), we obtain

\[\begin{array}{c}
  {\rm {\bf A}}_{ d,W}  = {\mathop {\lim} \limits_{\lambda \to 0}} \left(
{{\begin{array}{*{20}c}
 {{\frac{{c_{1}^{\left( {11} \right)} \lambda ^{m - 1} + \ldots +
c_{r}^{\left( {11} \right)} \lambda ^{m - r}}}{{\lambda ^{m} + d_{1}
\lambda ^{m - 1} + \ldots + d_{r} \lambda ^{m - r}}}}} \hfill &
{\ldots}  \hfill & {{\frac{{c_{1}^{\left( {1n} \right)} \lambda ^{m
- 1} + \ldots + c_{r}^{\left( {1n} \right)} \lambda ^{m -
r}}}{{\lambda ^{m} + d_{1} \lambda
^{m - 1} + \ldots + d_{r} \lambda ^{m - r}}}}} \hfill \\
 {\ldots}  \hfill & {\ldots}  \hfill & {\ldots}  \hfill \\
 {{\frac{{c_{1}^{\left( {m1} \right)} \lambda ^{m - 1} + \ldots +
c_{r}^{\left( {m1} \right)} \lambda ^{m - r}}}{{\lambda ^{m} + d_{1}
\lambda ^{m - 1} + \ldots + d_{r} \lambda ^{m - r}}}}} \hfill &
{\ldots}  \hfill & {{\frac{{c_{1}^{\left( {mn} \right)} \lambda ^{m
- 1} + \ldots + c_{r}^{\left( {mn} \right)} \lambda ^{m -
r}}}{{\lambda ^{m} + d_{1} \lambda
^{m - 1} + \ldots + d_{r} \lambda ^{m - r}}}}} \hfill \\
\end{array}} } \right) =\\
  \left( {{\begin{array}{*{20}c}
 {{\frac{{c_{r}^{\left( {11} \right)}} }{{d_{r}} }}} \hfill & {\ldots}
\hfill & {{\frac{{c_{r}^{\left( {1n} \right)}} }{{d_{r}} }}} \hfill \\
 {\ldots}  \hfill & {\ldots}  \hfill & {\ldots}  \hfill \\
 {{\frac{{c_{r}^{\left( {m1} \right)}} }{{d_{r}} }}} \hfill & {\ldots}
\hfill & {{\frac{{c_{r}^{\left( {mn} \right)}} }{{d_{r}} }}} \hfill \\
\end{array}} } \right).
\end{array}
\]
where $c_{r}^{\left( {ij} \right)} = \sum\limits_{\beta \in
J_{r,\,m} {\left\{ {i} \right\}}} {\left| \left( {\left(
{{\rm {\bf A}}^{k+1}} \right)_{\,.\,i} \left( {{\rm {\bf
a}}_{.j}^{ (k)} }  \right)} \right){\kern 1pt}  _{\beta}
^{\beta} \right|}  $ and $d_{r} = {\sum\limits_{\beta \in J_{r,\,m}}
{{\left| {\left(
{{\rm {\bf A}}^{k+1}} \right){\kern
1pt}  _{\beta} ^{\beta} } \right|}}} $. Thus, we have
obtained the determinantal representation of ${\rm {\bf A}}_{d,W}  $ by (\ref{eq:dr_rep_wcdet}).

By analogy can be proved  (\ref{eq:dr_rep_wrdet}).
$\blacksquare$

\section{ Cramer's rules for generalized inverse solutions of systems of  linear
equations}
An obvious consequence of a determinantal representation of the inverse matrix by the classical adjoint matrix is the Cramer rule.
As we know, Cramer rule gives an explicit expression for the solution of nonsingular linear equations. In \cite{ro}, Robinson gave an elegant proof of Cramer rule which aroused great interest in finding determinantal formulas
for solutions of some restricted linear equations both consistent and  nonconsistent. It has been widely discussed by Robinson \cite{ro},
Ben-Israel \cite{ben2}, Verghese \cite{ver}, Werner \cite{wer}, Chen \cite{ch}, Ji \cite{ji} ,Wang \cite{waG1}, Wei \cite{wei4}.

In this section we demonstrate that the obtained analogues of the adjoint matrix for the generalized inverse matrices enable us to obtain
 natural analogues of Cramer's rule for generalized inverse solutions of systems of  linear equations.

\subsection{Cramer's rule for the least
squares solution with the minimum norm}
\begin{definition} Suppose in a complex system of  linear
equations:
\begin{equation}
\label{eq:Ax=b} {\rm {\bf A}} \cdot {\rm {\bf x}} = {\rm {\bf y}}
\end{equation}
the coefficient matrix ${\rm {\bf A}} \in {\mathbb C}^{m\times
n}_{r}$ and a column of constants ${\rm {\bf y}} = \left( {y_{1}
,\ldots ,y_{m} } \right)^{T}\in {\mathbb C}^{m}$. The least
squares solution with the minimum norm of  (\ref{eq:Ax=b})
 is  the vector ${\rm {\bf x}}^{0}\in\mathbb{C}^{n} $
  satisfying
 \[{\left\| {{\rm
{\bf x}}^{0}} \right\|} = {\mathop {\min} \limits_{{\rm {\bf
\tilde {x}}} \in\mathbb{C}^{n}} } {\left\{ {{\left\| {{\rm {\bf
\tilde {x}}}} \right\|}\,\vert\, {\left\| {{\rm {\bf A}} \cdot
{\rm {\bf \tilde {x}}} - {\rm {\bf y}}} \right\|} = {\mathop
{\min} \limits_{{\rm {\bf x}} \in \mathbb{C}^{n}} }{\left\| {{\rm
{\bf A}} \cdot {\rm {\bf x}} - {\rm {\bf y}}} \right\|}}
\right\}},\] where $\mathbb{C}^{n}$ is an $n$-dimension complex
vector space.
\end{definition}
If  the equation (\ref{eq:Ax=b}) has no precision solutions, then
${\rm {\bf x}}^{0}$ is its optimal approximation.

The following important proposition is well-known.

 \begin{theorem}
\cite{ho}  The vector ${\rm {\bf x}} = {\rm {\bf A}}^{ +} {\rm
{\bf y}}$
 is the least squares
solution with the minimum norm of the system (\ref{eq:Ax=b}).
\end{theorem}
\begin{theorem}\label{kyrc13}
 The following
statements are true for the system of  linear equations
(\ref{eq:Ax=b}).
\begin{itemize}
\item [ i)] If $\rank{\rm {\bf A}} = n$,
 then the components of the least squares
solution with the minimum norm ${\rm {\bf x}}^{0}=\left( {{x^{0}_{1}} ,\ldots ,x^{0}_{n} }
\right)^{T}$ are obtained by the formula
\begin{equation}
\label{eq:x01} x_{j}^{0} = {\frac{{\det  ({\rm {\bf A}}^{ *} {\rm
{\bf A}})_{.\,j} \left( {{\rm {\bf f}}} \right)}}{{\det {\rm {\bf
A}}^{ *} {\rm {\bf A}}}}}, \quad \left( {\forall j = \overline
{1,n}} \right),
\end{equation}
 where  ${\rm {\bf f}} = {\rm {\bf A}}^{ *} {\rm {\bf y}}$.

 \item [ii)] If $\rank{\rm {\bf A}} = r \le m < n$, then
\begin{equation}
\label{eq:x02} x_{j}^{0} = {\frac{{\sum\limits_{\beta \in J_{r,n}
{\left\{ {j} \right\}}} {{\left| {\left( {({\rm {\bf A}}^{ *} {\rm
{\bf A}})_{.\,j} ({\rm {\bf f}}  )} \right)_{\beta} ^{\beta} }
\right|}}}}{{d_{r} \left( {{\rm {\bf A}}^{ *} {\rm {\bf A}}}
\right)}}}, \quad \left( {\forall j = \overline {1,n}} \right).
\end{equation}
\end{itemize}
\end{theorem}

{\textit{Proof.}} i) If $\rank{\rm {\bf A}} = n$, then  we can
represent ${\rm {\bf A}}^{ +} $ by (\ref{eq:A*A_ful}). By multiplying
${\rm {\bf A}}^{ +}$  into ${\rm {\bf y}}$ we get (\ref{eq:x01}).

ii) If $\rank{\rm {\bf A}} = k \le m < n$, then    ${\rm {\bf
A}}^{ +} $ can be  represented by  (\ref{eq:det_repr_A*A}). By multiplying
${\rm {\bf A}}^{ +}$  into ${\rm {\bf y}}$ the least squares
solution with the minimum norm of the linear system (\ref{eq:Ax=b}) is given by components
as in (\ref{eq:x02}). $\blacksquare$

Using (\ref{eq3}) and (\ref{eq:AA*_ful}), we can obtain another
representation of the Cramer rule for the least squares solution with the minimum norm
of a
 linear system.
\begin{theorem}
 The following statements are
true for a system of  linear equations written in the form ${\rm
{\bf x}}\cdot {\rm {\bf A}} = {\rm {\bf y}}$.
\begin{itemize}
\item [ i)] If $\rank{\rm {\bf A}} = m$,
 then the components of the least squares
solution ${\rm {\bf x}}^{0}={\rm {\bf y}}{\rm {\bf A}}^{+}$ are
obtained by the formula
\[ x_{i}^{0} = {\frac{{\det  ({\rm {\bf A}}{\rm {\bf A}}^{ *} )_{i\,.}
\left( {{\rm {\bf g}}} \right)}}{{\det {\rm {\bf A}}{\rm {\bf
A}}^{ *} }}}, \quad \left( {\forall i = \overline {1,m}} \right),
\]
 where  ${\rm {\bf g}} = {\rm {\bf y}}{\rm {\bf A}}^{ *} $.

 \item [ii)] If $\rank{\rm {\bf A}} = r \le  n< m$, then
\[ x_{i}^{0} = {\frac{{\sum\limits_{\alpha \in I_{r,m} {\left\{ {i} \right\}}}
{{\left| {\left( {({\rm {\bf A}}{\rm {\bf A}}^{ *} )_{\,i\,.}
({\rm {\bf g}}  )} \right)_{\alpha} ^{\alpha} } \right|}}}}{{d_{r}
\left( {{\rm {\bf A}}{\rm {\bf A}}^{ *} } \right)}}}, \quad \left(
{\forall i = \overline {1,m}} \right).
\]
\end{itemize}
\end{theorem}
{\textit{Proof}}. The proof of this theorem is analogous to that
of Theorem \ref{kyrc13}.
\begin{remark}
The obtained formulas of the Cramer rule for the least squares
solution differ from  similar formulas in \cite{ben2,wer,ch,
ji,waG1}. They give  a  closer approximation to the Cramer
rule for consistent nonsingular systems of linear equations.
\end{remark}

\subsection{Cramer's rule for the Drazin inverse solution}

 In some situations,
however, people pay more attention to the Drazin inverse solution
of  singular linear systems \cite{waG,chen,si,wei3}.

Consider a general
system of linear equations (\ref{eq:Ax=b}), where ${\rm {\bf A}} \in
{\mathbb C}^{n\times n}$ and ${\rm {\bf x}}$, ${\rm {\bf y}}$ are
 vectors in ${\mathbb C}^{n}$. $R({\rm {\bf A}})$ denotes the
range of ${\rm {\bf A}}$ and $N({\rm {\bf A}})$ denotes the null
space of ${\rm {\bf A}}$.

The characteristic of the Drazin inverse
solution ${\rm {\bf A}}^{D}\rm {\bf y}$ is given in \cite{wei} by
the following theorem.

\begin{theorem}\label{kyrc14}
Let ${\rm {\bf A}} \in {\mathbb C}^{n\times n}$ with $Ind(A) = k$.
Then ${\rm {\bf A}}^{D}{\rm {\bf y}}$ is both the unique
solution in $R({\rm {\bf A}}^{k})$ of
\begin{equation}\label{eq15}
{\rm {\bf A}}^{k+1}\rm {\bf x} = {\rm {\bf A}}^{k}{\rm {\bf y}},
\end{equation}
and the unique minimal ${\rm {\bf P}}$-norm least squares solution
of (\ref{eq:Ax=b}).
\end{theorem}

\begin{remark}
 The ${\rm {\bf P}}$-norm is defined as $\|{\rm
{\bf x}}\|_{{\rm {\bf P}}} = \|{\rm {\bf P}}^{-1}{\rm {\bf x}}\|$
for ${\rm {\bf x}}\in {\mathbb C}^{n}$, where ${\rm {\bf P}}$ is a
nonsingular matrix that transforms ${\rm {\bf A}}$ into its Jordan
canonical form (\ref{eq:jord_dr}).
\end{remark}

In other words, the the Drazin inverse solution  ${\bf x}={\rm {\bf A}}^{D}{\bf y}$ is the unique solution of the problem: for a given ${\bf A}$ and a given vector ${\bf y} \in R({\rm {\bf A}}^{k})$, find a  vector ${\bf x} \in R({\rm {\bf A}}^{k})$ satisfying ${\rm {\bf A}}{\bf x}={\bf y}$, $Ind \, {\bf A}=k$.

In general, unlike ${\rm {\bf A}}^{+} {\rm {\bf y}}$, the Drazin inverse solution ${\rm {\bf A}}^{D} {\rm {\bf y}}$ is not a true
solution of a singular system (\ref{eq:Ax=b}), even if the system is consistent. However, Theorem
\ref{kyrc14} means that ${\rm {\bf A}}^{D} {\rm {\bf y}}$ is the unique minimal P-norm least squares solution of (\ref{eq:Ax=b}).

 The determinantal representation of the  ${\rm {\bf P}}$-norm least
squares solution a system of linear equations (\ref{eq:Ax=b}) by the determinantal representation of the Drazin inverse (\ref{eq:det_dr_st}) is introduced in \cite{st3}.

We obtain  Cramer's rule for the  ${\rm {\bf P}}$-norm least
squares solution (the Drazin inverse solution) of (\ref{eq:Ax=b}) in the following theorem.
\begin{theorem}
Let ${\rm {\bf A}} \in {\mathbb C}^{n\times n}$ with $Ind({\rm
{\bf A}}) = k$ and $\rank{\rm {\bf
A}}^{k + 1} = \rank{\rm {\bf A}}^{k}=r$. Then the unique minimal ${\rm {\bf P}}$-norm least
squares solution ${\rm {\bf \widehat{x}}}=(\widehat{x}_{1},
\ldots,\widehat{x}_{n})^{T}$ of the system (\ref{eq:Ax=b}) is given
by
\begin{equation}
\label{eq16} \widehat{x}_{i} = {\frac{{{\sum\limits_{\beta \in
J_{r,\,n} {\left\{ {i} \right\}}} {{\left| {\left( {{\rm {\bf
A}}_{.{\kern 1pt} i}^{k + 1} \left( {{\rm {\bf f}}} \right)}
\right)_{\beta} ^{\beta} }  \right|}} }}}{{{\sum\limits_{\beta \in
J_{r,n}} {{\left| {\left( {{\rm {\bf A}}^{k + 1}} \right)_{\beta}
^{\beta} }  \right|}}} }}} \quad \forall i = \overline {1,n},
\end{equation}
where ${\rm {\bf f}} = {\rm {\bf A}}^{k}{\rm {\bf y}}.$
\end{theorem}
{\textit{Proof}}. Representing  the Drazin inverse by
(\ref{eq:dr_repr_col}) and by virtue of Theorem \ref{kyrc14}, we have
\[
{\rm {\bf \widehat{x}}} =\begin{pmatrix}
  \widehat{x}_{1} \\
  \ldots \\
  \widehat{x}_{n}
\end{pmatrix}={\rm {\bf A}}^{D}{\rm {\bf y}}=\frac{1}{d_{r}
\left( {{\rm {\bf A}}^{k+1} } \right)}\begin{pmatrix}
  \sum\limits_{s = 1}^{n} {d_{1s} y_{s}}   \\
  \ldots \\
 \sum\limits_{s = 1}^{n} {d_{ns} y_{s}}
\end{pmatrix}.
\]
Therefore,
\[
\widehat{x}_{i} = {\frac{{1}}{{d_{r} \left( {{\rm {\bf A}}^{k +
1}} \right)}}}{\sum\limits_{s = 1\,\,}^{n} {{\sum\limits_{\beta
\in \, J_{r,n} {\left\{ {i} \right\}}} {{\left| {\left( {{\rm {\bf
A}}_{.{\kern 1pt} i}^{k + 1} \left( {{\rm {\bf a}}_{.\,s}^{\left(
{k} \right)}}  \right)} \right)_{\beta} ^{\beta} }  \right|}}} }}
\cdot y_{s} =
\]
\[
 = {\frac{{1}}{{d_{r} \left( {{\rm {\bf A}}^{k + 1}}
\right)}}}{\sum\limits_{\beta \in J_{r,\,n} {\left\{ {i}
\right\}}} {\,{\sum\limits_{\,s = 1\,\,}^{n} {{\left| {\left(
{{\rm {\bf A}}_{.\,i}^{k + 1} \left( {{\rm {\bf a}}_{.\,s}^{\left(
{k} \right)}}  \right)} \right)_{\beta} ^{\beta} }  \right|} \cdot
y_{s}} } }}=
\]
\[
 = {\frac{{1}}{{d_{r} \left( {{\rm {\bf A}}^{k + 1}}
\right)}}}{\sum\limits_{\beta \in J_{r,\,n} {\left\{ {i}
\right\}}} {\;{\sum\limits_{\,s = 1\,\,}^{n} {{\left| {\left(
{{\rm {\bf A}}_{.\,i}^{k + 1} \left( {{\rm {\bf a}}_{.\,s}^{\left(
{k} \right)} \cdot y_{s}}  \right)} \right)_{\beta} ^{\beta} }
\right|}}} }}.
\]
From this (\ref{eq16}) follows immediately. $\blacksquare$

If we shall present a system of linear equations as,
\begin{equation}\label{eq:xA=b}
\rm {\bf x}{\rm {\bf A}} = {\rm {\bf y}},
\end{equation}
where ${\rm {\bf A}} \in {\mathbb C}^{n\times n}$ with $Ind(A) = k$ and $\rank{\rm {\bf
A}}^{k + 1} = \rank{\rm {\bf A}}^{k}=r$, then by using the determinantal representation of the Drazin inverse
(\ref{eq:dr_repr_row}) we shall obtain the following analog of Cramer's rule for the  Drazin inverse solution of (\ref{eq:xA=b}):
\[ \widehat{x}_{i} = {\frac{{{\sum\limits_{\alpha \in
I_{r,\,n} {\left\{ {i} \right\}}} {{\left| {\left( {{\rm {\bf
A}}_{i{\kern 1pt} .}^{k + 1} \left( {{\rm {\bf g}}} \right)}
\right)_{\alpha} ^{\alpha} }  \right|}} }}}{{{\sum\limits_{\alpha \in
I_{r,n}} {{\left| {\left( {{\rm {\bf A}}^{k + 1}} \right)_{\alpha}
^{\alpha} }  \right|}}} }}}, \quad \forall i = \overline {1,n},
\]
where ${\rm {\bf g}} = {\rm {\bf y}}{\rm {\bf A}}^{k}.$

\subsection{Cramer's rule for the W-weighted Drazin inverse solution}
Consider restricted linear equations
\begin{equation}\label{eq:WAWx=y}
{\bf W} {\rm {\bf A}}{\bf W}{\rm {\bf x}} = {\rm {\bf y}},
\end{equation}
where ${\rm {\bf
A}}\in {\mathbb{C}}^{m\times n}$, ${\rm {\bf
W}}\in {\mathbb{C}}^{n\times m}$, $k_{1}=Ind({\rm {\bf
A}}{\rm {\bf
W}})$, $k_{2}=Ind({\rm {\bf
W}}{\rm {\bf
A}})$ with ${\bf y}\in ({\rm {\bf
W}}{\rm {\bf
A}})^{k_{2}}$ and $\rank({\rm {\bf
W}}{\rm {\bf
A}})^{k_{2}}=\rank({\rm {\bf
A}}{\rm {\bf
W}})^{k_{1}}=r$.

In \cite{wei4}, Wei has showed that there exists an unique solution ${\bf
A}_{d,W}{\bf y}$ of the linear equations (\ref{eq:WAWx=y}) and given a Cramer rule for
the W-weighted Drazin inverse solution of (\ref{eq:WAWx=y}) by the following theorem.
\begin{theorem}\label{theor:wei_Wdr}
Let ${\rm {\bf A}}$, ${\bf W}$ be the same as in (\ref{eq:WAWx=y}). Suppose that ${\rm {\bf
U}}\in {\mathbb{C}}^{n\times (n-r)}_{n-r}$ and ${\rm {\bf
V}}^{*}\in {\mathbb{C}}^{m\times (m-r)}_{m-r}$ be matrices whose columns form bases for $N(({\rm {\bf
W}}{\rm {\bf
A}})^{k_{2}})$ and $N(({\rm {\bf
A}}{\rm {\bf
W}})^{k_{1}})$, respectively. Then the unique W-weighted Drazin inverse solution ${\bf x}=(x_{1},...,x_{m}$ of (\ref{eq:WAWx=y}) satisfies
\[
x_{i}=\det \left(
             \begin{array}{cc}
               {\bf W} {\rm {\bf A}}{\bf W}(i\rightarrow {\bf y}) & {\rm {\bf
U}} \\
               {\rm {\bf
V}}(i\rightarrow {\bf 0}) & {\bf 0} \\
             \end{array}
           \right)\left/ \det \left(
             \begin{array}{cc}
               {\bf W} {\rm {\bf A}}{\bf W} & {\rm {\bf
U}} \\
               {\rm {\bf
V}} & {\bf 0} \\
             \end{array}
           \right)\right.
\]
where $i=1,2,..,m$.

\end{theorem}

Let  $k= {\rm max}\{{k_{1}}, {k_{2}}\}$. Denote ${\bf f}= ({\rm {\bf A}}{\bf W})^{ k}{\bf A}\cdot{\bf y}$. Then by Theorem
\ref{theor:det_rep_wdraz} using the  determinantal
representation (\ref{eq:dr_rep_wcdet}) of the W-weighted Drazin inverse ${\rm {\bf
A}}_{d,W}$, we evidently obtain the following Cramer's rule the W-weighted Drazin inverse solution of (\ref{eq:WAWx=y}),

\begin{equation}
\label{eq:cr_rul_wcdet} x_{i}  = {\frac{\sum\limits_{\beta
\in J_{r,\,m} {\left\{ {i} \right\}}} \left| \left(
{\left( {{\rm {\bf A}}{\bf W}} \right)^{k+2}_{\,. \,i} \left( {{\rm {\bf
f}}}   \right)} \right){\kern 1pt} {\kern 1pt} _{\beta}
^{\beta} \right|  }{{{\sum\limits_{\beta \in J_{r,\,\,m}} {{\left|
{\left( { {\bf A}{\bf W}} \right)^{k+2}{\kern 1pt} _{\beta} ^{\beta}
}  \right|}}} }}},
\end{equation}
where $i=1,2,..,m$.
\begin{remark}
Note that for (\ref{eq:cr_rul_wcdet}) unlike Theorem
\ref{theor:wei_Wdr}, we do not need auxiliary matrices ${\rm {\bf
U}}$ and ${\rm {\bf
V}}$.
\end{remark}

\subsection{Examples}
1. Let us consider the system of linear equations.
\begin{equation} \label{eq17}
\left\{ \begin{array}{c}
  2x_1-5x_3+4x_4=1, \\
   7x_1-4x_2-9x_3+1.5x_4=2, \\
  3x_1-4x_2+7x_3-6.5x_4=3, \\
  x_1-4x_2+12x_3-10.5x_4=1.
\end{array}\right.\end{equation}
The coefficient matrix of the system is the matrix  ${\rm \bf
A}=\begin{pmatrix}
  2 & 0 & -5 & 4 \\
  7 & -4 & -9 & 1.5 \\
  3 & -4 & 7 & -6.5 \\
  1 & -4 & 12 & -10.5
\end{pmatrix}$. We calculate the rank of ${\rm {\bf A}}$ which is equal to 3, and we have
\[{\rm \bf A}^{\ast}=\begin{pmatrix}
  2 & 7 & 3 & 1 \\
  0 & -4 & -4 & -4 \\
  -5 & -9 & 7 & 12 \\
  4 & 1.5 & -6.5 & -10.5
\end{pmatrix}, {\rm \bf A}^{*}{\rm \bf A}=\begin{pmatrix}
  63 & -44 & -40 & -11.5 \\
-44 & 48 & -40 & 62 \\
  -40 & -40 & 299 & -205 \\
  -11.5 & 62 & -205 & 170.75
\end{pmatrix}.\]
At first we  obtain  entries of ${\rm \bf A}^{+}$ by (\ref{eq:A*A_ful}):
\[\begin{array}{c}
  d_{3}({\rm \bf A}^{*}{\rm \bf A})=\left|\begin{array}{ccc}
  63 & -44 & -40 \\
  -44 & 48 & -40 \\
  -40 & -40 & 299
\end{array}\right|+\left|\begin{array}{ccc}
  63 & -44 & -11.5 \\
  -44 & 48 & 62 \\
  -11.5 & 62 & 170.75
\end{array}\right|+ \\
  +\left|\begin{array}{ccc}
  63 & -40 & -11.5 \\
  -40 & 299 & -205 \\
  -11.5 & -205 & 170.75
\end{array}\right|+\left|\begin{array}{ccc}
  48 & -40 & 62 \\
  -40 & 299 & -205 \\
  62 & -205 & 170.75
\end{array}\right|=102060,
\end{array}
\]
\[\begin{array}{l}
   l_{11}=\left|\begin{array}{ccc}
  2 & -44 & -40 \\
  0 & 48 & -40 \\
  -5 & -40 & 299
\end{array}\right|+\left|\begin{array}{ccc}
  2 & -44 & -11.5 \\
  0 & 48 & 62 \\
  4 & 62 & 170.75
\end{array}\right|+\left|\begin{array}{ccc}
  2 & -40 & -11.5 \\
  -5 & 299 & -205 \\
  4 & -205 & 170.75
\end{array}\right|=\\
  =25779,
\end{array}
\] and so forth.
Continuing in the same way, we get
\[
{\rm \bf A}^{+}=\frac{1}{102060}\begin{pmatrix}
  25779 & -4905 & 20742 & -5037 \\
  -3840 & -2880 & -4800 & -960 \\
  28350 & -17010 & 22680 & -5670 \\
  39558 & -18810 & 26484 & -13074
\end{pmatrix}.\]
Now we  obtain the least squares solution of the system
(\ref{eq17}) by the matrix method.
 \[{\rm \bf x}^{0}=
   \begin{pmatrix}
  x^0_{1} \\
  x^0_{2} \\
  x^0_{3} \\
  x^0_{4}
\end{pmatrix}=\frac{1}{102060}\begin{pmatrix}
  25779 & -4905 & 20742 & -5037 \\
  -3840 & -2880 & -4800 & -960 \\
  28350 & -17010 & 22680 & -5670 \\
  39558 & -18810 & 26484 & -13074
\end{pmatrix}\cdot\begin{pmatrix}
 1\\
  2 \\
  3 \\
  1
\end{pmatrix}=\]

 \[ =\frac{1}{102060}\begin{pmatrix}
  73158 \\
  -24960 \\
  56700 \\
  68316
\end{pmatrix}=\begin{pmatrix}
  \frac{12193}{17010}\\
  -\frac{416}{1071} \\
 \frac{5}{9} \\
  \frac{5693}{8505}
\end{pmatrix}
\]
Next we get the least squares solution with minimum norm of the system (\ref{eq17})
by the Cramer rule (\ref{eq:x02}), where
\[{\rm \bf f}=\begin{pmatrix}
  2 & 7 & 3 & 1 \\
  0 & -4 & -4 & -4 \\
  -5 & -9 & 7 & 12 \\
  4 & 1.5 & -6.5 & -10.5
\end{pmatrix}\cdot\begin{pmatrix}
 1\\
  2 \\
  3 \\
  1
\end{pmatrix}=\begin{pmatrix}
 26\\
  -24 \\
  10 \\
  -23
\end{pmatrix}.
\]
Thus we have
\[
  x^0_1=\frac{1}{102060}\left(
  \left|\begin{array}{ccc}
  26 & -44 & -40 \\
  -24 & 48 & -40 \\
  10 & -40 & 299
\end{array}\right|+
\left|\begin{array}{ccc}
  26 & -44 & -11.5 \\
  -24 & 48 & 62 \\
  -23 & 62 & 170.75
\end{array}\right|+\right.\]
\[\left.+\left|\begin{array}{ccc}
  26 & -40 & -11.5 \\
  10 & 299 & -205 \\
  23 & -205 & 170.75
\end{array}\right|\right)=
\frac{73158}{102060}=\frac{12193}{17010};\]
\[  x^0_2=\frac{1}{102060}\left(\left|\begin{array}{ccc}
  63 & 26 & -40 \\
  -44 & -24 & -40 \\
  -40 & 10 & 299
\end{array}\right|+\left|\begin{array}{ccc}
  63 & 26 & -11.5 \\
  -44 & -24 & 62 \\
  -11.5 & -23 & 170.75
\end{array}\right|+\right. \]\[
 \left.+\left|\begin{array}{ccc}
  -24 & -40 & 62 \\
  10 & 299 & -205 \\
  -23 & -205 & 170.75
\end{array}\right|\right)=\frac{-24960}{102060}=-\frac{416}{1071};
\]
\[
  x^0_3=\frac{1}{102060}\left(\left|\begin{array}{ccc}
  63 & -44 & 26 \\
  -44 & 48 & -24 \\
  -40 & -40 & 10
\end{array}\right|
  +\left|\begin{array}{ccc}
  63 & 26 & -11.5 \\
  -40 & 10 & -205 \\
  -11.5 & -23 & 170.75
\end{array}\right|+\right.\]\[\left.+\left|\begin{array}{ccc}
  48 & -24 & 62 \\
  -40 & 10 & -205 \\
  62 & -23 & 170.75
\end{array}\right|\right)=\frac{56700}{102060}=\frac{5}{9};
\]
\[
  x^0_4=\frac{1}{102060}\left(\left|\begin{array}{ccc}
  63 & -44 & 26 \\
  -44 & 48 & -24 \\
  -11.5 & 62 & -23
\end{array}\right| \\
  +\left|\begin{array}{ccc}
  63 & -40 & 26 \\
  -40 & 299 & 10 \\
  -11.5 & -205 & -23
\end{array}\right|+\right.\]
\[\left.+\left|\begin{array}{ccc}
  48 & -40 & -24 \\
  -40 & 299 & 10 \\
  62 & -205 &-23
\end{array}\right|\right)=\frac{68316}{102060}=\frac{5693}{8505}.
\]

2. Let us consider the following system of linear equations.
\begin{equation} \label{eq18}
\left\{ \begin{array}{c}
x_1-x_2+x_3+x_4=1, \\
  x_2-x_3+x_4=2, \\
  x_1-x_2+x_3+2x_4=3, \\
  x_1-x_2+x_3+x_4=1.
\end{array}\right.\end{equation}
The coefficient matrix of the system is the matrix ${\rm \bf
A}=\begin{pmatrix}
  1 & -1 & 1 & 1 \\
  0 & 1 & -1 & 1 \\
  1 & -1 & 1 & 2 \\
  1 & -1 & 1 & 1
\end{pmatrix}$. It is easy to verify the following:
\[{\rm \bf
A}^{2}=\begin{pmatrix}
  3 & -4 & 4 & 3 \\
  0 & 1 & -1 & 0 \\
  4 & -5 & 5 & 4 \\
  3 & -4 & 4 & 3
\end{pmatrix},\,\,{\rm \bf
A}^{3}=\begin{pmatrix}
  10 & -14 & 14 & 10 \\
  -1 & 2 & -2 & -1 \\
  13 & -18 & 18 & 13 \\
  10 & -14 & 14 & 10
\end{pmatrix},\] and $\rank{\rm \bf A}=3$,
 $\rank{\rm {\bf A}}^{2}=\rank{\rm \bf A}^{3}=2$.
 This implies $k=Ind({\rm {\bf A}})=2$.
We  obtain  entries of ${\rm \bf A}^{D}$ by   (\ref{eq:dr_repr_col}).
\[\begin{array}{c}
   d_{2}({\rm \bf A}^{3})=\left|\begin{array}{cc}
  10 & -14 \\
  -1 & 2
\end{array}\right|+\left|\begin{array}{cc}
  10 & 14 \\
  13 & 18
\end{array}\right|+\left|\begin{array}{cc}
  10 & 10 \\
  10 & 10
\end{array}\right|\\
   +\left|\begin{array}{cc}
 2 & -2 \\
  -18 & 18
\end{array}\right|+\left|\begin{array}{cc}
  2 & -1 \\
  -14 & 10
\end{array}\right|+\left|\begin{array}{cc}
  18 & 13 \\
  14 & 10
\end{array}\right|=8,
\end{array}
\]
\[   d_{11}=\left|\begin{array}{cc}
 3 & -14 \\
  0 & 2
\end{array}\right|+\left|\begin{array}{cc}
  3 & 14 \\
  4 & 18
\end{array}\right|+\left|\begin{array}{cc}
  3 & 10 \\
  3 & 10
\end{array}\right|=4,
\] and so forth.

Continuing in the same way, we get ${\rm \bf
A}^{D}=\begin{pmatrix}
  0.5 & 0.5 & -0.5 & 0.5 \\
 1.75 & 2.5 & -2.5 & 1.75 \\
  1.25 & 1.5 & -1.5 & 1.25 \\
  0.5 & 0.5 & -0.5 & 0.5
\end{pmatrix}.$
Now we obtain the Drazin inverse solution ${\rm {\bf
\widehat{x}}}$ of the system (\ref{eq18}) by the Cramer rule
(\ref{eq16}), where
\[{\rm \bf g}={\rm \bf
A}^{2}{\rm {\bf y}}=\begin{pmatrix}
  3 & -4 & 4 & 3 \\
  0 & 1 & -1 & 0 \\
  4 & -5 & 5 & 4 \\
  3 & -4 & 4 & 3
\end{pmatrix}\cdot\begin{pmatrix}
 1\\
  2 \\
  3 \\
  1
\end{pmatrix}=\begin{pmatrix}
 10\\
  -1 \\
  13 \\
  10
\end{pmatrix}.
\]
Thus we have
\[
  \widehat{x}_1=\frac{1}{8}\left(
  \left|\begin{array}{cc}
 10 & -14 \\
  -1 & 2
\end{array}\right|+\left|\begin{array}{cc}
  10 & 14 \\
  13 & 18
\end{array}\right|+\left|\begin{array}{cc}
  10 & 10 \\
  10 & 10
\end{array}\right|\right)=\frac{1}{2},
\]
\[
  \widehat{x}_2=\frac{1}{8}\left(
  \left|\begin{array}{cc}
 10 &10 \\
  -1 & -1
\end{array}\right|+\left|\begin{array}{cc}
  -1 & -2 \\
  13 & 18
\end{array}\right|+\left|\begin{array}{cc}
 -1 & -1 \\
  10 & 10
\end{array}\right|\right)=1,
\]
\[
  \widehat{x}_3=\frac{1}{8}\left(
  \left|\begin{array}{cc}
 10 &10 \\
  13 & 13
\end{array}\right|+\left|\begin{array}{cc}
  2 & -1 \\
  -18 & 13
\end{array}\right|+\left|\begin{array}{cc}
13 & 13 \\
  10 & 10
\end{array}\right|\right)=1,
\]
\[
  \widehat{x}_4=\frac{1}{8}\left(
  \left|\begin{array}{cc}
 10 & 10 \\
  10 & 10
\end{array}\right|+\left|\begin{array}{cc}
  2 & -1 \\
  -14 & 10
\end{array}\right|+\left|\begin{array}{cc}
  18 & 13 \\
  14 & 10
\end{array}\right|\right)=\frac{1}{2}.\]

\section{Cramer's rule of the generalized inverse solutions of some  matrix
equations}
Matrix equation is one of the important study fields of linear
algebra. Linear matrix equations, such as
\begin{equation}\label{eq:AX=C}
 {\rm {\bf A}}{\rm {\bf X}} = {\rm {\bf C}},
\end{equation}
\begin{equation}\label{eq:XB=D}
 {\rm {\bf X}}{\rm {\bf B}} = {\rm {\bf D}},
\end{equation}
 and \begin{equation}\label{eq_s:AXB=D}
 {\rm {\bf A}}{\rm {\bf X}}{\rm {\bf B}} = {\rm {\bf
D}},
\end{equation}
play an important role in linear system theory therefore a large
number of papers have presented several methods for solving these
matrix equations \cite{da,de,do,peng,vet}. In \cite{kh}, Khatri and
Mitra studied the Hermitian solutions to the  matrix equations
(\ref{eq:AX=B}) and (\ref{eq_s:AXB=D}) over the complex field and
the system of the equations  (\ref{eq:AX=B}) and (\ref{eq:XB=D}).
Wang, in \cite{wa1,wa2}, and Li and Wu, in \cite{li}   studied the
bisymmetric, symmetric and skew-antisymmetric least squares
solution to this system over the  quaternion skew field. Extreme
ranks of real matrices in least squares solution of the equation
(\ref{eq_s:AXB=D}) was investigated in \cite{liu} over the complex
field and in \cite{wa3} over the  quaternion skew field.

As we know, the Cramer rule gives an explicit expression for the
solution of nonsingular linear equations. The  Robinson's result  ( \cite{ro})
 aroused
great interest in finding determinantal representations of least squares solution as analogs of the
Cramer rule for the matrix equations (for example,
\cite{wang,qi,ya}).
 The Cramer rule for solutions of the
restricted matrix equations (\ref{eq:AX=B}), (\ref{eq:XB=D}) and
(\ref{eq_s:AXB=D}) was established in \cite{wag,gu,so1}.

In this section, we  obtain analogs of the
Cramer rule for generalized inverse solutions of the matrix equations
(\ref{eq:AX=B}), (\ref{eq:XB=D}) and (\ref{eq:AXB=D}) without any
restriction, namely for the minimum norm least squares solution  and Drazin inverse solution.

 We shall show  numerical examples to illustrate the main results as well.

\subsection{Cramer's rule for the minimum norm least squares solution of some  matrix
equations}
\begin{definition}
Consider a  matrix equation
\begin{equation}
\label{eq:AX=B}
{\bf A}{\bf X}={\bf B},
\end{equation} where ${\rm {\bf
A}}\in {\mathbb{C}}^{m\times n},{\rm {\bf B}} \in
{\mathbb{C}}^{m\times s} $ are given, ${\rm {\bf X}} \in
{\mathbb{C}}^{n\times s}$ is unknown. Suppose
\[S_{1}=\{{\rm {\bf X}}|{\rm {\bf X}} \in
{\mathbb{C}}^{n\times s}, \|{\rm {\bf A}}{\rm {\bf X}} - {\rm {\bf
B}}\|=\min \}.\] Then matrices ${\rm {\bf X}} \in
{\mathbb{C}}^{n\times s}$ such that ${\rm {\bf X}}\in S_{1}$ are
called least squares solutions of the matrix equation
(\ref{eq:AX=B}).
 If
${\rm {\bf X}}_{LS}={\min}_ {{\rm {\bf X}}\in S_{1}}\|{\rm {\bf
X}}\|$, then ${\rm {\bf X}}_{LS}$ is called the minimum norm least squares
solution of (\ref{eq:AX=B}).
\end{definition}
If  the equation (\ref{eq:AX=B}) has no precision solutions, then
${\rm {\bf X}}_{LS}$ is its optimal approximation.

The following important proposition is well-known.

\begin{lemma}(\cite{ji})\label{theor:LS} The least squares solutions of (\ref{eq:AX=B}) are
\[{\rm {\bf X}} = {\rm {\bf A}}^{+}{\rm {\bf B}} + ({\rm {\bf
I}}_{n} - {\rm {\bf A}}^{+} {\rm {\bf A}}){\rm {\bf C}},\] where
${\rm {\bf A}}\in {\mathbb{C}}^{m\times n},{\rm {\bf B}} \in
{\mathbb{C}}^{m\times s} $ are given, and ${\rm {\bf C}}\in
{\mathbb{C}}^{n\times s}$ is an arbitrary
 matrix. The least squares minimum norm solution is ${\rm {\bf
X}}_{LS}={\rm {\bf A}}^{+}{\rm {\bf B}}$.
\end{lemma}
 We denote ${\rm {\bf A}}^{ \ast}{\rm {\bf B}}=:\hat{{\rm
{\bf B}}}= (\hat{b}_{ij})\in {\mathbb{C}}^{n\times s}$.
\begin{theorem}\label{theor:ls_AX}
\begin{enumerate}
\item[(i)] If ${\rm rank}\,{\rm {\bf A}} = r \le m < n$, then for the minimum norm least squares solution ${\rm {\bf X}}_{LS}=(x_{ij})\in
{\mathbb{C}}^{n\times s}$ for all
$i = \overline {1,n} $, $j = \overline {1,s} $ we have
\begin{equation}
\label{eq:ls_AX} x_{ij} = {\frac{{{\sum\limits_{\beta \in
J_{r,\,n} {\left\{ {i} \right\}}} {{\left| \left( {\left( {{\rm
{\bf A}}^{ *} {\rm {\bf A}}} \right)_{\,.\,i} \left( {\hat{{\rm
{\bf b}}}_{.j}} \right)} \right){\kern 1pt} {\kern 1pt} _{\beta}
^{\beta} \right|}} } }}{{{\sum\limits_{\beta \in J_{r,\,\,n}}
{{\left| {\left( {{\rm {\bf A}}^{ *} {\rm {\bf A}}} \right){\kern
1pt} {\kern 1pt} _{\beta} ^{\beta} }  \right|}}} }}}.
\end{equation}
\item[(ii)] If ${\rm rank}\,{\rm {\bf A}} = n$, then   for all $i = \overline {1,n} $, $j =
\overline {1,s} $ we have
\begin{equation}
\label{eq:ls_AX_full} x_{i\,j} = {\frac{{\det ({\rm {\bf A}}^{ *}
{\rm {\bf A}})_{.\,i} \left( {\hat{{\rm {\bf b}}}_{.j}}
\right)}}{{ \det ({\rm {\bf A}}^{ *} {\rm {\bf A}})}}}
\end{equation}
\noindent where $\hat{{\rm {\bf b}}}_{.j}$ is the $j$th column of
$\hat{{\rm {\bf B}}}$ for all $j = \overline {1,s} $.
\end{enumerate}
\end{theorem}

{\textit{Proof}}. i) If ${\rm rank}\,{\rm {\bf A}} = r \le m < n$, then
by Theorem  \ref{theor:det_repr_MP} we can represent the matrix
${\rm {\bf A}}^{ +} $ by (\ref{eq:det_repr_A*A}). Therefore, we
obtain for all $i = \overline {1,n} $, $j = \overline {1,s} $
\[
x_{ij}  =\sum_{k=1}^{m}
a_{ik}^{+}b_{kj}=\sum_{k=1}^{m}{\frac{{{\sum\limits_{\beta \in
J_{r,\,n} {\left\{ {i} \right\}}} {{\left| \left( {\left( {{\rm
{\bf A}}^{ *} {\rm {\bf A}}} \right)_{\,. \,i} \left( {{\rm {\bf
a}}_{.k}^{ *} }  \right)} \right) {\kern 1pt} _{\beta}
^{\beta}\right|} } } }}{{{\sum\limits_{\beta \in J_{r,\,\,n}}
{{\left| {\left( {{\rm {\bf A}}^{ *} {\rm {\bf A}}} \right){\kern
1pt} _{\beta} ^{\beta} }  \right|}}} }}}\cdot b_{kj}=
\]
\[
{\frac{{{\sum\limits_{\beta \in J_{r,\,n} {\left\{ {i}
\right\}}}\sum_{k=1}^{m} {{\left| \left( {\left( {{\rm {\bf A}}^{
*} {\rm {\bf A}}} \right)_{\,. \,i} \left( {{\rm {\bf a}}_{.k}^{
*} } \right)} \right) {\kern 1pt} _{\beta} ^{\beta}\right|} } }
}\cdot b_{kj}} {{{\sum\limits_{\beta \in J_{r,\,\,n}} {{\left|
{\left( {{\rm {\bf A}}^{ *} {\rm {\bf A}}} \right){\kern 1pt}
_{\beta} ^{\beta} } \right|}}} }}}.
\]
Since ${\sum\limits_{k} {{\rm {\bf a}}_{.\,k}^{ *}  b_{kj}} }=
\left( {{\begin{array}{*{20}c}
 {{\sum\limits_{k} {a_{1k}^{ *}  b_{kj}} } } \hfill \\
 {{\sum\limits_{k} {a_{2k}^{ *}  b_{kj}} } } \hfill \\
 { \vdots}  \hfill \\
 {{\sum\limits_{k} {a_{nk}^{ *}  b_{kj}} } } \hfill \\
\end{array}} } \right) = \hat{{\rm {\bf b}}}_{.j}$, then it follows (\ref{eq:ls_AX}).

\noindent (ii) The proof of this case is similarly  to that of  (i) by using Corollary \ref{cor:A+A-1}. $\blacksquare$
\begin{definition}
Consider a  matrix equation
\begin{equation}\label{eq:XA=B}
 {\rm {\bf X}}{\rm {\bf A}} = {\rm {\bf B}},
\end{equation}
\noindent  where ${\rm {\bf A}}\in {\mathbb{C}}^{m\times n},{\rm
{\bf B}} \in {\mathbb{C}}^{s\times n} $ are given, ${\rm {\bf X}}
\in {\mathbb{C}}^{s\times m}$ is unknown. Suppose
\[S_{2}=\{{\rm {\bf X}}|\,{\rm {\bf X}} \in
{\mathbb{C}}^{s\times m}, \|{\rm {\bf X}}{\rm {\bf A}} - {\rm {\bf
B}}\|=\min \}.\] Then matrices ${\rm {\bf X}} \in
{\mathbb{C}}^{s\times m}$ such that ${\rm {\bf X}}\in S_{2}$ are
called least squares solutions of the matrix equation
(\ref{eq:XA=B}).
 If
${\rm {\bf X}}_{LS}={\min}_ {{\rm {\bf X}}\in S_{2}}\|{\rm {\bf
X}}\|$, then ${\rm {\bf X}}_{LS}$ is called the minimum norm least squares
solution of (\ref{eq:XA=B}).
\end{definition}
The following  lemma can be obtained by analogy to Lemma
\ref{theor:LS}.

\begin{lemma} The least squares solutions of (\ref{eq:XA=B}) are
\[{\rm {\bf X}} = {\rm {\bf B}}{\rm {\bf A}}^{+} + {\rm {\bf C}}({\rm {\bf
I}}_{m} -  {\rm {\bf A}}{\rm {\bf A}}^{+}),\]  where ${\rm {\bf
A}}\in {\mathbb{C}}^{m\times n},{\rm {\bf B}} \in
{\mathbb{C}}^{s\times n} $ are given, and ${\rm {\bf C}}\in
{\mathbb{C}}^{s\times m}$ is an arbitrary matrix. The minimum norm
least squares solution is ${\rm {\bf X}}_{LS}={\rm {\bf B}}{\rm
{\bf A}}^{+}$.
\end{lemma}
 We denote  ${\rm {\bf B}}{\rm {\bf A}}^{
\ast}=:\check{{\rm {\bf B}}}= (\check{b}_{ij})\in
{\mathbb{C}}^{s\times m}$.
\begin{theorem}\label{theor:ls_XA}
\begin{enumerate}
\item[(i)] If ${\rm rank}\,{\rm {\bf A}} = r \le n < m$, then for the minimum norm least squares solution ${\rm {\bf X}}_{LS}=(x_{ij})\in
{\mathbb{C}}^{s\times m}$ for all $i = \overline {1,s} $, $j =
\overline {1,m} $ we have
\begin{equation}
\label{eq:ls_XA} x_{ij} = {\frac{{{\sum\limits_{\alpha \in I_{r,m}
{\left\{ {j} \right\}}} {{\left| \left( {\left( {{\rm {\bf A}}{\rm
{\bf A}}^{ *} } \right)_{\,j\,.} \left( {\check{{\rm {\bf
b}}}_{i\,.}} \right)} \right)\,_{\alpha} ^{\alpha}\right|} }
}}}{{{\sum\limits_{\alpha \in I_{r,\,m}}  {{\left| {\left( {{\rm
{\bf A}}{\rm {\bf A}}^{ *} } \right) {\kern 1pt} _{\alpha}
^{\alpha} } \right|}}} }}}.
\end{equation}
\item[(ii)] If ${\rm rank}\,{\rm {\bf A}} = m$, then for all $i = \overline {1,s} $, $j =
\overline {1,m} $ we have
\begin{equation}
\label{eq:ls_XA_full} x_{i\,j} = {\frac{{\det ({\rm {\bf A}}{\rm
{\bf A}}^{ *} )_{j.\,} \left( {\check{{\rm {\bf b}}}_{i\,.}}
\right)}}{{\det ({\rm {\bf A}}{\rm {\bf A}}^{ *} )}}}
\end{equation}
\noindent where $\check{{\rm {\bf b}}}_{i.}$ is the $i$th row of
$\check{{\rm {\bf B}}}$ for all $i = \overline {1,s} $.
\end{enumerate}
\end{theorem}
{\textit{Proof}}. (i) If ${\rm rank}\,{\rm {\bf A}} = r \le n <m $, then
by Theorem  \ref{theor:det_repr_MP} we can represent the matrix
${\rm {\bf A}}^{ +} $ by (\ref{eq:det_repr_AA*}). Therefore, for
all $i = \overline {1,s} $, $j = \overline {1,m} $ we obtain
\[
x_{ij}  =\sum_{k=1}^{n}b_{ik} a_{kj}^{+}=\sum_{k=1}^{n}b_{ik}\cdot
{\frac{{{\sum\limits_{\alpha \in I_{r,m} {\left\{ {j} \right\}}}
{{\left| \left( {\left( {{\rm {\bf A}}{\rm {\bf A}}^{ *} }
\right)_{\,j\,.} \left( {\bf a}_{k\,.}^{*} \right)}
\right)\,_{\alpha} ^{\alpha}\right|} } }}}{{{\sum\limits_{\alpha
\in I_{r,\,m}}  {{\left| {\left( {{\rm {\bf A}}{\rm {\bf A}}^{ *}
} \right) {\kern 1pt} _{\alpha} ^{\alpha} } \right|}}} }}} =
\]
\[
{\frac{{{\sum_{k=1}^{n}b_{ik} \sum\limits_{\alpha \in I_{r,m}
{\left\{ {j} \right\}}} {{\left| \left( {\left( {{\rm {\bf A}}{\rm
{\bf A}}^{ *} } \right)_{\,j\,.} \left( {\bf a}_{k\,.}^{*}
\right)} \right)\,_{\alpha} ^{\alpha}\right|} } } }}
{{{\sum\limits_{\alpha \in I_{r,\,m}}  {{\left| {\left( {{\rm {\bf
A}}{\rm {\bf A}}^{ *} } \right) {\kern 1pt} _{\alpha} ^{\alpha} }
\right|}}} }}}
\]
Since for all $i = \overline {1,s} $
 \[{\sum\limits_{k} {{  b_{ik}\rm {\bf a}}_{k\,.}^{ *}}
}=\begin{pmatrix}
  \sum\limits_{k} {b_{ik}a_{k1}^{ *} } & \sum\limits_{k} {b_{ik}a_{k2}^{ *} } & \cdots & \sum\limits_{k} {b_{ik}a_{km}^{ *} }
\end{pmatrix}
 = \check{{\rm {\bf b}}}_{i.},\]
then it follows (\ref{eq:ls_XA}).

\noindent (ii) The proof of this case is similarly  to that of  (i) by using Corollary \ref{cor:A+A-1}. $\blacksquare$

\begin{definition}
Consider a  matrix equation
\begin{equation}
\label{eq:AXB=D}
{\bf A}{\bf X}{\bf B}={\bf D},
\end{equation}   where $ {\rm {\bf
A}}\in{\rm {\mathbb{C}}}^{m \times n}_{r_{1}},{\rm {\bf B}}\in{\rm
{\mathbb{C}}}^{p \times q}_{r_{2}}, {\rm {\bf D}}\in{\rm
{\mathbb{C}}}^{m \times q}$ are given, $ {\rm {\bf X}}\in{\rm
{\mathbb{C}}}^{n \times p}$ is unknown. Suppose
\[S_{3}=\{{\rm {\bf X}}|\,{\rm {\bf X}} \in
{\mathbb{C}}^{n \times p}, \|{\rm {\bf A}}{\rm {\bf X}}{\rm {\bf
B}} - {\rm {\bf D}}\|=\min \}.\] Then matrices ${\rm {\bf X}} \in
{\mathbb{C}}^{n \times p}$ such that ${\rm {\bf X}}\in S_{3}$ are
called least squares solutions of the matrix equation
(\ref{eq:AXB=D}).
 If
${\rm {\bf X}}_{LS}={\min}_ {{\rm {\bf X}}\in S_{3}}\|{\rm {\bf
X}}\|$, then ${\rm {\bf X}}_{LS}$ is called the minimum norm least squares
solution of (\ref{eq:AXB=D}).
\end{definition}
The following important proposition is well-known.
\begin{lemma}(\cite{wer})\label{theor:LS_AXB} The least squares solutions of (\ref{eq:AXB=D}) are
\[{\rm {\bf X}} = {\rm {\bf A}}^{+}{\rm {\bf
D}}{\rm {\bf B}}^{+} + ({\rm {\bf I}}_{n} - {\rm {\bf A}}^{+} {\rm
{\bf A}}){\rm {\bf V}}+{\rm {\bf W}}({\rm {\bf I}}_{p} -{\rm{\bf
B}}{\rm{\bf B}}^{+}),\]  where $ {\rm {\bf A}}\in{\rm
{\mathbb{C}}}^{m \times n}_{r_{1}},{\rm {\bf B}}\in{\rm
{\mathbb{C}}}^{p \times q}_{r_{2}}, {\rm {\bf D}}\in{\rm
{\mathbb{C}}}^{m \times q}$ are given, and $\{{\rm {\bf V}},{\rm
{\bf W}}\}\subset {\mathbb{C}}^{n\times p}$  are  arbitrary
quaternion matrices. The minimum norm least squares solution is
${\rm {\bf X}}_{LS}={\rm {\bf A}}^{+}{\rm {\bf D}}{\rm {\bf
B}}^{+}$.
\end{lemma}
 We denote  ${\rm {\bf \widetilde{D}}}= {\rm {\bf
A}}^\ast{\rm {\bf D}}{\rm {\bf B}}^\ast$.
\begin{theorem}\label{theor:AXB=D}
\begin{enumerate}
\item[(i)]  If ${\rm rank}\,{\rm {\bf A}} = r_{1} <  n$ and ${\rm rank}\,{\rm {\bf B}} = r_{2} < p$, then for the minimum norm least squares solution ${\rm {\bf X}}_{LS}=(x_{ij})\in
{\mathbb{C}}^{n\times p}$  of (\ref{eq:AXB=D}) we have
\begin{equation}\label{eq:d^B}
x_{ij} = {\frac{{{\sum\limits_{\beta \in J_{r_{1},\,n} {\left\{
{i} \right\}}} { \left| {\left( {{\rm {\bf A}}^{ *} {\rm {\bf A}}}
\right)_{\,.\,i} \left( {{{\rm {\bf d}}}\,_{.\,j}^{{\rm {\bf B}}}}
\right)\, _{\beta} ^{\beta}} \right| } } }}{{{\sum\limits_{\beta
\in J_{r_{1},n}} {{\left| {\left( {{\rm {\bf A}}^{ *} {\rm {\bf
A}}} \right)_{\beta} ^{\beta} } \right|}} \sum\limits_{\alpha \in
I_{r_{2},p}}{{\left| {\left( {{\rm {\bf B}}{\rm {\bf B}}^{ *} }
\right) _{\alpha} ^{\alpha} } \right|}}} }}},
\end{equation}
or
\begin{equation}\label{eq:d^A}
 x_{ij}={\frac{{{\sum\limits_{\alpha
\in I_{r_{2},p} {\left\{ {j} \right\}}} { \left| {\left( {{\rm
{\bf B}}{\rm {\bf B}}^{ *} } \right)_{\,j\,.} \left( {{{\rm {\bf
d}}}\,_{i\,.}^{{\rm {\bf A}}}} \right)\,_{\alpha} ^{\alpha}}
\right| } }}}{{{\sum\limits_{\beta \in J_{r_{1},n}} {{\left|
{\left( {{\rm {\bf A}}^{ *} {\rm {\bf A}}} \right) _{\beta}
^{\beta} } \right|}}\sum\limits_{\alpha \in I_{r_{2},p}} {{\left|
{\left( {{\rm {\bf B}}{\rm {\bf B}}^{ *} } \right) _{\alpha}
^{\alpha} } \right|}}} }}},
\end{equation}
where
  \begin{equation}
\label{eq:def_d^B_m}
   {{\rm {\bf d}}_{.\,j}^{{\rm {\bf B}}}}=\left[
\sum\limits_{\alpha \in I_{r_{2},p} {\left\{ {j} \right\}}} {
\left| {\left( {{\rm {\bf B}}{\rm {\bf B}}^{ *} } \right)_{j.}
\left( {\tilde{{\rm {\bf d}}}_{1.}} \right)\,_{\alpha} ^{\alpha}}
\right|},...,\sum\limits_{\alpha \in I_{r_{2},p} {\left\{ {j}
\right\}}} { \left| {\left( {{\rm {\bf B}}{\rm {\bf B}}^{ *} }
\right)_{j.} \left( {\tilde{{\rm {\bf d}}}_{n.}}
\right)\,_{\alpha} ^{\alpha}} \right|} \right]^{T},
\end{equation}
  \begin{equation}
\label{eq:def_d^A_m}
  {{\rm {\bf d}}_{i\,.}^{{\rm {\bf A}}}}=\left[
\sum\limits_{\beta \in J_{r_{1},n} {\left\{ {i} \right\}}} {
\left| {\left( {{\rm {\bf A}}^{ *}{\rm {\bf A}} } \right)_{.i}
\left( {\tilde{{\rm {\bf d}}}_{.1}} \right)\,_{\beta} ^{\beta}}
\right|},...,\sum\limits_{\alpha \in I_{r_{1},n} {\left\{ {i}
\right\}}} { \left| {\left( {{\rm {\bf A}}^{ *}{\rm {\bf A}} }
\right)_{.i} \left( {\tilde{{\rm {\bf d}}}_{.\,p}}
\right)\,_{\beta} ^{\beta}} \right|} \right]
\end{equation}
 are the column-vector and the row-vector, respectively. $\tilde{{\rm {\bf d}}}_{i\,.}$ is the  $i$th row of ${\rm {\bf
\widetilde{D}}}$ for all $i = \overline {1,n} $, and $\tilde{{\rm
{\bf d}}}_{.\,j}$ is the $j$th column of ${\rm {\bf
\widetilde{D}}}$ for all $j = \overline {1,p} $.
\item[(ii)]If ${\rm rank}\,{\rm {\bf A}} = n$ and ${\rm rank}\,{\rm {\bf B}} = p$, then for the
least squares solution ${\rm {\bf X}}_{LS}=(x_{ij})\in
{\mathbb{C}}^{n\times p}$  of (\ref{eq:AXB=D}) we have for all $i
= \overline {1,n} $, $j = \overline {1,p} $,
\begin{equation}
\label{eq:AXB_cdetA*A} x_{i\,j} = {\frac{{\det \left( {({\rm {\bf
A}}^{ *} {\rm {\bf A}})_{.\,i\,} \left( {{\rm {\bf d}}_{.j}^{{\rm
{\bf B}}}} \right)}\right)}}{{\det ({\rm {\bf A}}^{ *} {\rm {\bf
A}})\cdot \det ({\rm {\bf B}}{\rm {\bf B}}^{ *} )}}},
\end{equation}
or
\begin{equation}
\label{eq:AXB_rdetBB*} x_{i\,j} = {\frac{{\det\left( { ({\rm {\bf
B}}{\rm {\bf B}}^{ *} )_{j.\,} \left( {{\rm {\bf d}}_{i\,.}^{{\rm
{\bf A}}}} \right)}\right)}}{{\det ({\rm {\bf A}}^{ *} {\rm {\bf
A}})\cdot \det ({\rm {\bf B}}{\rm {\bf B}}^{ *} )}}},
\end{equation}
 \noindent where
  \begin{equation}
\label{eq:def_d^B}{\rm {\bf d}}_{.j}^{{\rm {\bf B}}} : = \left[
{\det\left( ({\rm {\bf B}}{\rm {\bf B}}^{ *} )_{j.\,} \left(
{\tilde{{\rm {\bf d}}}_{1\,.}^{}} \right)\right),\ldots
,\det\left( ({\rm {\bf B}}{\rm {\bf B}}^{ *} )_{j.\,} \left(
{\tilde{{\rm {\bf d}}}_{n\,.}^{}} \right)\right)} \right]^{T},
\end{equation}
   \begin{equation}
\label{eq:def_d^A}{\rm {\bf d}}_{i\,.}^{{\rm {\bf A}}} : = \left[
{\det  \left(({\rm {\bf A}}^{ *} {\rm {\bf A}})_{.\,i\,} \left(
{\tilde{{\rm {\bf d}}}_{.1}} \right)\right),\ldots ,\det \left(
({\rm {\bf A}}^{ *} {\rm {\bf A}})_{.\,i\,} \left( {\tilde{{\rm
{\bf d}}}_{.p}} \right)\right)} \right]
\end{equation}
are respectively the column-vector and the row-vector.

\item[(iii)] If ${\rm rank}\,{\rm {\bf A}} = n$ and ${\rm rank}\,{\rm {\bf B}} = r_{2} < p$, then for the
least squares solution ${\rm {\bf X}}_{LS}=(x_{ij})\in
{\mathbb{C}}^{n\times p}$  of (\ref{eq:AXB=D}) we have
\begin{equation}
\label{eq:AXB_detA*A_d^B} x_{ij}={\frac{{\det \left( {({\rm {\bf
A}}^{ *} {\rm {\bf A}})_{.\,i\,} \left( {{\rm {\bf d}}_{.j}^{{\rm
{\bf B}}}} \right)}\right)}}{{\det ({\rm {\bf A}}^{ *} {\rm {\bf
A}})\sum\limits_{\alpha \in I_{r_{2},p}}{{\left| {\left( {{\rm
{\bf B}}{\rm {\bf B}}^{ *} } \right) _{\alpha} ^{\alpha} }
\right|}}} }},
\end{equation}
or
\begin{equation}\label{AXB_detA*A_d^A}
x_{ij} = {\frac{{{\sum\limits_{\alpha \in I_{r_{2},p} {\left\{ {j}
\right\}}} { \left| {\left( {{\rm {\bf B}}{\rm {\bf B}}^{ *} }
\right)_{\,j\,.} \left( {{{\rm {\bf d}}}\,_{i\,.}^{{\rm {\bf A}}}}
\right)\,_{\alpha} ^{\alpha}} \right| } }}}{{\det ({\rm {\bf A}}^{
*} {\rm {\bf A}})\sum\limits_{\alpha \in I_{r_{2},p}}{{\left|
{\left( {{\rm {\bf B}}{\rm {\bf B}}^{ *} } \right) _{\alpha}
^{\alpha} } \right|}}} }},
\end{equation}
where $ {{\rm {\bf d}}_{.\,j}^{{\rm {\bf B}}}}$ is
(\ref{eq:def_d^B_m})
   and
 ${\rm {\bf
d}}_{i\,.}^{{\rm {\bf A}}} $ is (\ref{eq:def_d^A}).

\item[(iiii)]
If ${\rm rank}\,{\rm {\bf A}} = r_{1} <  m$ and
${\rm rank}\,{\rm {\bf B}} =  p$, then for the least squares
solution ${\rm {\bf X}}_{LS}=(x_{ij})\in {\mathbb{C}}^{n\times p}$
of (\ref{eq:AXB=D}) we have
\begin{equation}
\label{eq:AXB_detBB*_d^A} x_{i\,j} = {\frac{{\det\left( { ({\rm
{\bf B}}{\rm {\bf B}}^{ *} )_{j.\,} \left( {{\rm {\bf
d}}_{i\,.}^{{\rm {\bf A}}}} \right)}\right)}}{{{\sum\limits_{\beta
\in J_{r_{1},n}} {{\left| {\left( {{\rm {\bf A}}^{ *} {\rm {\bf
A}}} \right) _{\beta} ^{\beta} } \right|}}\cdot \det ({\rm {\bf
B}}{\rm {\bf B}}^{ *} )}}}},
\end{equation}
or
\begin{equation} \label{eq:AXB_detBB*_d^B} x_{i\,j}=
{\frac{{{\sum\limits_{\beta \in J_{r_{1},\,n} {\left\{ {i}
\right\}}} { \left| {\left( {{\rm {\bf A}}^{ *} {\rm {\bf A}}}
\right)_{\,.\,i} \left( {{{\rm {\bf d}}}\,_{.\,j}^{{\rm {\bf B}}}}
\right)\, _{\beta} ^{\beta}} \right| } } }}{{{\sum\limits_{\beta
\in J_{r_{1},n}} {{\left| {\left( {{\rm {\bf A}}^{ *} {\rm {\bf
A}}} \right)_{\beta} ^{\beta} } \right|}}\det ({\rm {\bf B}}{\rm
{\bf B}}^{ *} )}}}},
\end{equation}
 \noindent where $ {{\rm {\bf d}}_{.\,j}^{{\rm {\bf B}}}}$ is
(\ref{eq:def_d^B})
   and
 ${\rm {\bf
d}}_{i\,.}^{{\rm {\bf A}}} $ is (\ref{eq:def_d^A_m}).
\end{enumerate}
\end{theorem}
{\textit{Proof}}. (i) If ${\rm {\bf A}} \in {\rm {\mathbb{C}}}_{r_{1}}^{m\times n}
$, ${\rm {\bf B}} \in {\rm {\mathbb{C}}}_{r_{2}}^{p\times q} $ and
$ r_{1} <  n$, $ r_{2} < p$, then by Theorem
\ref{theor:det_repr_MP} the Moore-Penrose inverses ${\rm {\bf
A}}^{ +} = \left( {a_{ij}^{ +} } \right) \in {\rm
{\mathbb{C}}}_{}^{n\times m} $ and ${\rm {\bf B}}^{ +} = \left(
{b_{ij}^{ +} } \right) \in {\rm {\mathbb{C}}}^{q\times p} $
possess the following determinantal representations respectively,
\[
 a_{ij}^{ +}  = {\frac{{{\sum\limits_{\beta
\in J_{r_{1},\,n} {\left\{ {i} \right\}}} { \left| {\left( {{\rm
{\bf A}}^{ *} {\rm {\bf A}}} \right)_{\,. \,i} \left( {{\rm {\bf
a}}_{.j}^{ *} }  \right){\kern 1pt}  _{\beta} ^{\beta}} \right| }
} }}{{{\sum\limits_{\beta \in J_{r_{1},\,n}} {{\left| {\left(
{{\rm {\bf A}}^{ *} {\rm {\bf A}}} \right){\kern 1pt} _{\beta}
^{\beta} }  \right|}}} }}},
\]
\begin{equation}\label{eq:b+}
 b_{ij}^{ +}  =
{\frac{{{\sum\limits_{\alpha \in I_{r_{2},p} {\left\{ {j}
\right\}}} { \left| {({\rm {\bf B}}{\rm {\bf B}}^{ *} )_{j\,.\,}
({\rm {\bf b}}_{i.\,}^{ *} )\,_{\alpha} ^{\alpha}} \right| }
}}}{{{\sum\limits_{\alpha \in I_{r_{2},p}}  {{\left| {\left( {{\rm
{\bf B}}{\rm {\bf B}}^{ *} } \right){\kern 1pt} _{\alpha}
^{\alpha} } \right|}}} }}}.
\end{equation}
Since by Theorem \ref{theor:LS_AXB}  ${\rm {\bf X}}_{LS}={\rm {\bf
A}}^{+}{\rm {\bf D}}{\rm {\bf B}}^{+}$, then an entry of ${\rm
{\bf X}}_{LS}=(x_{ij})$ is

\begin{equation}
\label{eq:sum+} x_{ij} = {{\sum\limits_{s = 1}^{q} {\left(
{{\sum\limits_{k = 1}^{m} {{a}_{ik}^{+} d_{ks}} } } \right)}}
{b}_{sj}^{+}}.
\end{equation}
Denote  by $\hat{{\rm {\bf d}}_{.s}}$ the $s$th column of ${\rm
{\bf A}}^{ \ast}{\rm {\bf D}}=:\hat{{\rm {\bf D}}}=
(\hat{d}_{ij})\in {\mathbb{C}}^{n\times q}$ for all $s=\overline
{1,q}$. It follows from ${\sum\limits_{k} { {\rm {\bf a}}_{.\,k}^{
*}}d_{ks} }=\hat{{\rm {\bf d}}_{.\,s}}$ that
\[
\sum\limits_{k = 1}^{m} {{a}_{ik}^{+} d_{ks}}=\sum\limits_{k =
1}^{m}{\frac{{{\sum\limits_{\beta \in J_{r_{1},\,n} {\left\{ {i}
\right\}}} { \left| {\left( {{\rm {\bf A}}^{ *} {\rm {\bf A}}}
\right)_{\,. \,i} \left( {{\rm {\bf a}}_{.k}^{ *} } \right) {\kern
1pt} _{\beta} ^{\beta}} \right| } } }}{{{\sum\limits_{\beta \in
J_{r_{1},\,n}} {{\left| {\left( {{\rm {\bf A}}^{ *} {\rm {\bf A}}}
\right){\kern 1pt} _{\beta} ^{\beta} }  \right|}}} }}}\cdot
d_{ks}=
\]
\begin{equation}\label{eq:sum_cdet}
{\frac{{{\sum\limits_{\beta \in J_{r_{1},\,n} {\left\{ {i}
\right\}}}\sum\limits_{k = 1}^{m} { \left| {\left( {{\rm {\bf
A}}^{ *} {\rm {\bf A}}} \right)_{\,. \,i} \left( {{\rm {\bf
a}}_{.k}^{ *} } \right) {\kern 1pt} _{\beta} ^{\beta}} \right| } }
}\cdot d_{ks}}{{{\sum\limits_{\beta \in J_{r_{1},\,n}} {{\left|
{\left( {{\rm {\bf A}}^{ *} {\rm {\bf A}}} \right){\kern 1pt}
_{\beta} ^{\beta} }  \right|}}} }}}={\frac{{{\sum\limits_{\beta
\in J_{r_{1},\,n} {\left\{ {i} \right\}}} { \left| {\left( {{\rm
{\bf A}}^{ *} {\rm {\bf A}}} \right)_{\,. \,i} \left( \hat{{\rm
{\bf d}}_{.\,s}} \right) {\kern 1pt} _{\beta} ^{\beta}} \right| }
} }}{{{\sum\limits_{\beta \in J_{r_{1},\,n}} {{\left| {\left(
{{\rm {\bf A}}^{ *} {\rm {\bf A}}} \right){\kern 1pt} _{\beta}
^{\beta} }  \right|}}} }}}
\end{equation}
Suppose ${\rm {\bf e}}_{s.}$ and ${\rm {\bf e}}_{.\,s}$ are
respectively the unit row-vector and the unit column-vector whose
components are $0$, except the $s$th components, which are $1$.
Substituting  (\ref{eq:sum_cdet}) and (\ref{eq:b+}) in
(\ref{eq:sum+}), we obtain
\[
x_{ij} =\sum\limits_{s = 1}^{q}{\frac{{{\sum\limits_{\beta \in
J_{r_{1},\,n} {\left\{ {i} \right\}}} { \left| {\left( {{\rm {\bf
A}}^{ *} {\rm {\bf A}}} \right)_{\,. \,i} \left( \hat{{\rm {\bf
d}}_{.\,s}} \right) {\kern 1pt} _{\beta} ^{\beta}} \right| } }
}}{{{\sum\limits_{\beta \in J_{r_{1},\,n}} {{\left| {\left( {{\rm
{\bf A}}^{ *} {\rm {\bf A}}} \right){\kern 1pt} _{\beta} ^{\beta}
}  \right|}}} }}}{\frac{{{\sum\limits_{\alpha \in I_{r_{2},p}
{\left\{ {j} \right\}}} { \left| {({\rm {\bf B}}{\rm {\bf B}}^{ *}
)_{j\,.\,} ({\rm {\bf b}}_{s.\,}^{ *} )\,_{\alpha} ^{\alpha}}
\right| } }}}{{{\sum\limits_{\alpha \in I_{r_{2},p}}  {{\left|
{\left( {{\rm {\bf B}}{\rm {\bf B}}^{ *} } \right){\kern 1pt}
_{\alpha} ^{\alpha} } \right|}}} }}}.
\]
Since
\begin{equation}\label{eq:prop}
\hat{{\rm {\bf d}}_{.\,s}}=\sum\limits_{l = 1}^{n}{\rm {\bf
e}}_{.\,l}\hat{ d_{ls}},\,\,\,\,  {\rm {\bf b}}_{s.\,}^{
*}=\sum\limits_{t = 1}^{p}b_{st}^{*}{\rm {\bf e}}_{t.},\,\,\,\,
\sum\limits_{s=1}^{q}\hat{d_{ls}}b_{st}^{*}=\widetilde{d}_{lt},
\end{equation}
then we have
\[
x_{ij} = \]
\[{\frac{{ \sum\limits_{s = 1}^{q}\sum\limits_{t =
1}^{p} \sum\limits_{l = 1}^{n} {\sum\limits_{\beta \in
J_{r_{1},\,n} {\left\{ {i} \right\}}} { \left| {\left( {{\rm {\bf
A}}^{ *} {\rm {\bf A}}} \right)_{\,. \,i} \left( {\rm {\bf
e}}_{.\,l} \right) {\kern 1pt} _{\beta} ^{\beta}} \right| } }
}\hat{ d_{ls}}b_{st}^{*}{\sum\limits_{\alpha \in I_{r_{2},p}
{\left\{ {j} \right\}}} { \left| {({\rm {\bf B}}{\rm {\bf B}}^{ *}
)_{j\,.\,} ({\rm {\bf e}}_{t.} )\,_{\alpha} ^{\alpha}} \right| } }
}{{{\sum\limits_{\beta \in J_{r_{1},\,n}} {{\left| {\left( {{\rm
{\bf A}}^{ *} {\rm {\bf A}}} \right){\kern 1pt} _{\beta} ^{\beta}
}  \right|}}} }{{\sum\limits_{\alpha \in I_{r_{2},p}} {{\left|
{\left( {{\rm {\bf B}}{\rm {\bf B}}^{ *} } \right){\kern 1pt}
_{\alpha} ^{\alpha} } \right|}}} }}    }=
\]
\begin{equation}\label{eq:x_ij}
{\frac{{ \sum\limits_{t = 1}^{p} \sum\limits_{l = 1}^{n}
{\sum\limits_{\beta \in J_{r_{1},\,n} {\left\{ {i} \right\}}}
{\left| {\left( {{\rm {\bf A}}^{ *} {\rm {\bf A}}} \right)_{\,.
\,i} \left( {\rm {\bf e}}_{.\,l} \right) {\kern 1pt} _{\beta}
^{\beta}} \right| } } }\,\,\widetilde{d}_{lt}{\sum\limits_{\alpha
\in I_{r_{2},p} {\left\{ {j} \right\}}} { \left| {({\rm {\bf
B}}{\rm {\bf B}}^{ *} )_{j\,.\,} ({\rm {\bf e}}_{t.} )\,_{\alpha}
^{\alpha}} \right| } } }{{{\sum\limits_{\beta \in J_{r_{1},\,n}}
{{\left| {\left( {{\rm {\bf A}}^{ *} {\rm {\bf A}}} \right){\kern
1pt} _{\beta} ^{\beta} }  \right|}}} }{{\sum\limits_{\alpha \in
I_{r_{2},p}} {{\left| {\left( {{\rm {\bf B}}{\rm {\bf B}}^{ *} }
\right){\kern 1pt} _{\alpha} ^{\alpha} } \right|}}} }}    }.
\end{equation}
Denote by
\[
 d^{{\rm {\bf A}}}_{it}:= \]
\[
{\sum\limits_{\beta \in J_{r_{1},\,n} {\left\{ {i} \right\}}} {
\left| {\left( {{\rm {\bf A}}^{ *} {\rm {\bf A}}} \right)_{\,.
\,i} \left( \widetilde{{\rm {\bf d}}}_{.\,t} \right){\kern 1pt}
_{\beta} ^{\beta}} \right| } }= \sum\limits_{l = 1}^{n}
{\sum\limits_{\beta \in J_{r_{1},\,n} {\left\{ {i} \right\}}} {
\left| {\left( {{\rm {\bf A}}^{ *} {\rm {\bf A}}} \right)_{\,.
\,i} \left( {\rm {\bf e}}_{.\,l} \right){\kern 1pt}
 _{\beta} ^{\beta}} \right| } } \widetilde{d}_{lt}
\]
the $t$th component  of a row-vector ${\rm {\bf d}}^{{\rm {\bf
A}}}_{i\,.}= (d^{{\rm {\bf A}}}_{i1},...,d^{{\rm {\bf A}}}_{ip})$
for all $t=\overline {1,p}$. Substituting it in (\ref{eq:x_ij}),
we have
\[x_{ij} ={\frac{{ \sum\limits_{t = 1}^{p}
 d^{{\rm {\bf A}}}_{it}
}{\sum\limits_{\alpha \in I_{r_{2},p} {\left\{ {j} \right\}}} {
\left| {({\rm {\bf B}}{\rm {\bf B}}^{ *} )_{j\,.\,} ({\rm {\bf
e}}_{t.} )\,_{\alpha} ^{\alpha}} \right| } }
}{{{\sum\limits_{\beta \in J_{r_{1},\,n}} {{\left| {\left( {{\rm
{\bf A}}^{ *} {\rm {\bf A}}} \right){\kern 1pt} _{\beta} ^{\beta}
}  \right|}}} }{{\sum\limits_{\alpha \in I_{r_{2},p}} {{\left|
{\left( {{\rm {\bf B}}{\rm {\bf B}}^{ *} } \right){\kern 1pt}
_{\alpha} ^{\alpha} } \right|}}} }}    }.
\]
Since $\sum\limits_{t = 1}^{p}
 d^{{\rm {\bf A}}}_{it}{\rm {\bf e}}_{t.}={\rm {\bf
d}}^{{\rm {\bf A}}}_{i\,.}$, then it follows (\ref{eq:d^A}).

If we denote by
\begin{equation}\label{eq:d^B_den}
 d^{{\rm {\bf B}}}_{lj}:=
\sum\limits_{t = 1}^{p}\widetilde{d}_{lt}{\sum\limits_{\alpha \in
I_{r_{2},p} {\left\{ {j} \right\}}} { \left| {({\rm {\bf B}}{\rm
{\bf B}}^{ *} )_{j\,.\,} ({\rm {\bf e}}_{t.} )\,_{\alpha}
^{\alpha}} \right| } }={\sum\limits_{\alpha \in I_{r_{2},p}
{\left\{ {j} \right\}}} {\left| {({\rm {\bf B}}{\rm {\bf B}}^{ *}
)_{j\,.\,} (\widetilde{{\rm {\bf d}}}_{l.} )\,_{\alpha} ^{\alpha}}
\right| } }
\end{equation}
the $l$th component  of a column-vector ${\rm {\bf d}}^{{\rm {\bf
B}}}_{.\,j}= (d^{{\rm {\bf B}}}_{1j},...,d^{{\rm {\bf
B}}}_{jn})^{T}$ for all $l=\overline {1,n}$ and substitute it in
(\ref{eq:x_ij}), we obtain
\[x_{ij} ={\frac{{  \sum\limits_{l = 1}^{n}
{\sum\limits_{\beta \in J_{r_{1},\,n} {\left\{ {i} \right\}}} {
\left| {\left( {{\rm {\bf A}}^{ *} {\rm {\bf A}}} \right)_{\,.
\,i} \left( {\rm {\bf e}}_{.\,l} \right){\kern 1pt} _{\beta}
^{\beta}} \right| } } }\,\,d^{{\rm {\bf B}}}_{lj}
}{{{\sum\limits_{\beta \in J_{r_{1},\,n}} {{\left| {\left( {{\rm
{\bf A}}^{ *} {\rm {\bf A}}} \right){\kern 1pt} _{\beta} ^{\beta}
}  \right|}}} }{{\sum\limits_{\alpha \in I_{r_{2},p}} {{\left|
{\left( {{\rm {\bf B}}{\rm {\bf B}}^{ *} } \right){\kern 1pt}
_{\alpha} ^{\alpha} } \right|}}} }}    }.
\]
Since $\sum\limits_{l = 1}^{n}{\rm {\bf e}}_{.l}
 d^{{\rm {\bf B}}}_{lj}={\rm {\bf
d}}^{{\rm {\bf B}}}_{.\,j}$, then it follows (\ref{eq:d^B}).

(ii)  If ${\rm rank}\,{\rm {\bf A}} = n$ and ${\rm rank}\,{\rm {\bf B}} = p$,
  then by Corollary \ref{cor:A+A-1} ${\rm {\bf
A}}^{ +}  = \left( {{\rm {\bf A}}^{ *} {\rm {\bf A}}} \right)^{ -
1}{\rm {\bf A}}^{ *} $ and ${\rm {\bf B}}^{ +}  = {\rm {\bf B}}^{
* }\left( {{\rm {\bf B}}{\rm {\bf B}}^{ *} } \right)^{ - 1}$.
 Therefore, we obtain
\[\begin{array}{c}
{\rm {\bf X}}_{LS} = ({\rm {\bf A}}^{ \ast}{\rm {\bf A}})^{ -
1}{\rm {\bf A}}^{ \ast}{\rm {\bf D}}{\rm {\bf B}}^{ *} \left(
{{\rm {\bf
B}}{\rm {\bf B}}^{ *} } \right)^{ - 1}=\\
 = \begin{pmatrix}
   x_{11} & x_{12} & \ldots & x_{1p} \\
   x_{21} & x_{22} & \ldots & x_{2p} \\
   \ldots & \ldots & \ldots & \ldots \\
   x_{n1} & x_{n2} & \ldots & x_{np} \
 \end{pmatrix}
 = {\frac{{1}}{{\det ({\rm {\bf A}}^{ \ast}{\rm {\bf A}})}}}\begin{pmatrix}
  {L} _{11}^{{\rm {\bf A}}} & {L} _{21}^{{\rm {\bf A}}}&
   \ldots & {L} _{n1}^{{\rm {\bf A}}} \\
  {L} _{12}^{{\rm {\bf A}}} & {L} _{22}^{{\rm {\bf A}}} &
  \ldots & {L} _{n2}^{{\rm {\bf A}}} \\
  \ldots & \ldots & \ldots & \ldots \\
 {L} _{1n}^{{\rm {\bf A}}} & {L} _{2n}^{{\rm {\bf A}}}
  & \ldots & {L} _{nn}^{{\rm {\bf A}}}
\end{pmatrix}\times\\
  \times\begin{pmatrix}
    \tilde{d}_{11} & \tilde{d}_{12} & \ldots & \tilde{d}_{1m} \\
    \tilde{d}_{21} & \tilde{d}_{22} & \ldots & \tilde{d}_{2m} \\
    \ldots & \ldots & \ldots & \ldots \\
    \tilde{d}_{n1} & \tilde{d}_{n2} & \ldots & \tilde{d}_{nm} \
  \end{pmatrix}

{\frac{{1}}{{\det \left( {{\rm {\bf B}}{\rm {\bf B}}^{ *} }
\right)}}}
\begin{pmatrix}
 {R} _{\, 11}^{{\rm {\bf B}}} & {R} _{\, 21}^{{\rm {\bf B}}}
 &\ldots & {R} _{\, p1}^{{\rm {\bf B}}} \\
 {R} _{\, 12}^{{\rm {\bf B}}} & {R} _{\, 22}^{{\rm {\bf B}}} &\ldots &
 {R} _{\, p2}^{{\rm {\bf B}}} \\
 \ldots  & \ldots & \ldots & \ldots \\
 {R} _{\, 1p}^{{\rm {\bf B}}} & {R} _{\, 2p}^{{\rm {\bf B}}} &
 \ldots & {R} _{\, pp}^{{\rm {\bf B}}}
\end{pmatrix},
\end{array}
\]
\noindent where $\tilde{d}_{ij}$ is  $ij$th entry of the matrix
${\rm {\bf \widetilde{D}}}$, ${L}_{ij}^{{\rm {\bf A}}} $ is the
$ij$th
 cofactor of $({\rm {\bf A}}^{ \ast}{\rm {\bf A}})$ for all $i,j =  \overline {1,n} $ and ${R}_{i{\kern 1pt} j}^{{\rm {\bf
B}}} $ is the  $ij$th cofactor of $\left( {{\rm {\bf B}}{\rm {\bf
B}}^{ *} } \right)$ for all $i,j = \overline {1,p} $. This implies
\begin{equation}
\label{eq:sum}x_{ij} = {\frac{{{\sum\limits_{k = 1}^{n} {{
L}_{ki}^{{\rm {\bf A}}}} }\left( {{\sum\limits_{s = 1}^{p}
{\tilde{d}_{\,ks}} } {R}_{js}^{{\rm {\bf B}}}} \right)}}{{\det
({\rm {\bf A}}^{ *} {\rm {\bf A}})\cdot\det ({\rm {\bf B}}{\rm
{\bf B}}^{ *} )}}},
\end{equation}
 \noindent for all $i = \overline {1,n} $, $j =
\overline {1,p} $. We obtain the sum in parentheses  and denote it
as follows
\[{\sum\limits_{s
= 1}^{p} {\tilde{d}_{k\,s}} } {R}_{j\,s}^{{\rm {\bf B}}} = \det
({\rm {\bf B}}{\rm {\bf B}}^{ *} )_{j.\,} \left( {\tilde{{\rm {\bf
d}}}_{k\,.}} \right):=d_{k\,j}^{{\rm {\bf B}}},\] where
$\tilde{{\rm {\bf d}}}_{k\,.} $ is the $k$th row-vector of
$\tilde{{\rm {\bf D}}}$ for all $k = \overline {1,n} $. Suppose
${\rm {\bf d}}_{.\,j}^{{\rm {\bf B}}} : = \left(d_{1\,j}^{{\rm
{\bf B}}},\ldots ,d_{n\,j}^{{\rm {\bf B}}} \right)^{T}$
 is the column-vector for all
$j = \overline {1,p} $. Reducing the sum ${\sum\limits_{k = 1}^{n}
{{ L}_{ki}^{{\rm {\bf A}}}} }d_{k\,j}^{{\rm {\bf B}}} $, we obtain
an analog of Cramer's rule for (\ref{eq:AXB=D}) by
(\ref{eq:AXB_cdetA*A}).

Interchanging the order of summation in (\ref{eq:sum}), we have
\[
 x_{ij} = {\frac{{{\sum\limits_{s = 1}^{p} {\left(
{{\sum\limits_{k = 1}^{n} {{L}_{ki}^{{\rm {\bf A}}}
\tilde{d}_{\,ks}} } } \right)}} {R}_{js}^{{\rm {\bf B}}}} }{{\det
({\rm {\bf A}}^{ *} {\rm {\bf A}})\cdot\det ({\rm {\bf B}}{\rm
{\bf B}}^{ *} )}}}.\]
We obtain the sum in parentheses  and denote
it as follows
\[{\sum\limits_{k = 1}^{n}
{{L}_{ki}^{{\rm {\bf A}}} \tilde{d}_{k\,s}} }  = \det ({\rm {\bf
A}}^{ *} {\rm {\bf A}})_{.\,i\,} \left( {\tilde{{\rm {\bf
d}}}_{.\,s}} \right)=:d_{i\,s}^{{\rm {\bf A}}},\] \noindent where
$\tilde{{\rm {\bf d}}}_{.\,s} $ is the $s$th column-vector of
$\tilde{{\rm {\bf D}}}$ for all $s = \overline {1,p} $. Suppose
${{{\rm {\bf d}}}_{i\,.}^{{\rm {\bf A}}}} : = \left(
d_{i\,1}^{{\rm {\bf A}}},\ldots , d_{i\,p}^{{\rm {\bf A}}}
\right)$ is the row-vector for all $i = \overline {1,n} $.
Reducing the sum ${{\sum\limits_{s = 1}^{n} d_{i\,s}^{{\rm {\bf
A}}}} {R}_{js}^{{\rm {\bf B}}}} $, we obtain another analog of
Cramer's rule for the least squares solutions of (\ref{eq:AXB=D})
by (\ref{eq:AXB_rdetBB*}).

(iii) If ${\rm {\bf A}} \in {\rm {\mathbb{C}}}_{r_{1}}^{m\times n}
$, ${\rm {\bf B}} \in {\rm {\mathbb{C}}}_{r_{2}}^{p\times q} $ and
$ r_{1} =n$, $ r_{2} < p$, then by Remark \ref{rem:det_repr_MP_A*A}
and   Theorem \ref{theor:det_repr_MP} the Moore-Penrose inverses
${\rm {\bf A}}^{ +} = \left( {a_{ij}^{ +} } \right) \in {\rm
{\mathbb{C}}}_{}^{n\times m} $ and ${\rm {\bf B}}^{ +} = \left(
{b_{ij}^{ +} } \right) \in {\rm {\mathbb{C}}}^{q\times p} $
possess the following determinantal representations respectively,
\[
 a_{ij}^{ +}  = {\frac{\det {\left( {{\rm
{\bf A}}^{ *} {\rm {\bf A}}} \right)_{\,. \,i} \left( {{\rm {\bf
a}}_{.j}^{ *} }  \right)}   } {\det {\left( {{\rm {\bf A}}^{ *}
{\rm {\bf A}}} \right)}}},
\]
\begin{equation}\label{eq:b+2}
 b_{ij}^{ +}  =
{\frac{{{\sum\limits_{\alpha \in I_{r_{2},p} {\left\{ {j}
\right\}}} { \left| {({\rm {\bf B}}{\rm {\bf B}}^{ *} )_{j\,.\,}
({\rm {\bf b}}_{i.\,}^{ *} )\,_{\alpha} ^{\alpha}} \right| }
}}}{{{\sum\limits_{\alpha \in I_{r_{2},p}}  {{\left| {\left( {{\rm
{\bf B}}{\rm {\bf B}}^{ *} } \right){\kern 1pt} _{\alpha}
^{\alpha} } \right|}}} }}}.
\end{equation}
Since by Theorem \ref{theor:LS_AXB}  ${\rm {\bf X}}_{LS}={\rm {\bf
A}}^{+}{\rm {\bf D}}{\rm {\bf B}}^{+}$, then an entry of ${\rm
{\bf X}}_{LS}=(x_{ij})$ is (\ref{eq:sum+}). Denote  by $\hat{{\rm
{\bf d}}_{.s}}$ the $s$th column of ${\rm {\bf A}}^{ \ast}{\rm
{\bf D}}=:\hat{{\rm {\bf D}}}= (\hat{d}_{ij})\in
{\mathbb{C}}^{n\times q}$ for all $s=\overline {1,q}$. It follows
from ${\sum\limits_{k} { {\rm {\bf a}}_{.\,k}^{ *}}d_{ks}
}=\hat{{\rm {\bf d}}_{.\,s}}$ that
\begin{equation}\label{eq:sum_det}
\sum\limits_{k = 1}^{m} {{a}_{ik}^{+} d_{ks}}=\sum\limits_{k =
1}^{m}{\frac{\det {\left( {{\rm {\bf A}}^{ *} {\rm {\bf A}}}
\right)_{\,. \,i} \left( {{\rm {\bf a}}_{.k}^{ *} }  \right)}   }
{\det {\left( {{\rm {\bf A}}^{ *} {\rm {\bf A}}} \right)}}}\cdot
d_{ks}={\frac{\det {\left( {{\rm {\bf A}}^{ *} {\rm {\bf A}}}
\right)_{\,. \,i} \left( {\hat{{\rm {\bf d}}_{.\,s}} }  \right)} }
{\det {\left( {{\rm {\bf A}}^{ *} {\rm {\bf A}}} \right)}}}
\end{equation}
Substituting  (\ref{eq:sum_det}) and (\ref{eq:b+2}) in
(\ref{eq:sum+}), and using (\ref{eq:prop}) we have
\[
x_{ij} =\sum\limits_{s = 1}^{q}{\frac{\det {\left( {{\rm {\bf
A}}^{ *} {\rm {\bf A}}} \right)_{\,. \,i} \left( {\hat{{\rm {\bf
d}}_{.\,s}} }  \right)} } {\det {\left( {{\rm {\bf A}}^{ *} {\rm
{\bf A}}} \right)}}}{\frac{{{\sum\limits_{\alpha \in I_{r_{2},p}
{\left\{ {j} \right\}}} { \left| {({\rm {\bf B}}{\rm {\bf B}}^{ *}
)_{j\,.\,} ({\rm {\bf b}}_{s.\,}^{ *} )\,_{\alpha} ^{\alpha}}
\right| } }}}{{{\sum\limits_{\alpha \in I_{r_{2},p}} {{\left|
{\left( {{\rm {\bf B}}{\rm {\bf B}}^{ *} } \right){\kern 1pt}
_{\alpha} ^{\alpha} } \right|}}} }}}=
\]
\[{\frac{{ \sum\limits_{s = 1}^{q}\sum\limits_{t =
1}^{p} \sum\limits_{l = 1}^{n} \det {\left( {{\rm {\bf A}}^{ *}
{\rm {\bf A}}} \right)_{\,. \,i} \left( {{\rm {\bf e}}_{.\,l} }
\right)} }\hat{ d_{ls}}b_{st}^{*}{\sum\limits_{\alpha \in
I_{r_{2},p} {\left\{ {j} \right\}}} { \left| {({\rm {\bf B}}{\rm
{\bf B}}^{ *} )_{j\,.\,} ({\rm {\bf e}}_{t.} )\,_{\alpha}
^{\alpha}} \right| } } }{{{\det {\left( {{\rm {\bf A}}^{ *} {\rm
{\bf A}}} \right)}} }{{\sum\limits_{\alpha \in I_{r_{2},p}}
{{\left| {\left( {{\rm {\bf B}}{\rm {\bf B}}^{ *} } \right){\kern
1pt} _{\alpha} ^{\alpha} } \right|}}} }}    }=
\]
\begin{equation}\label{eq:x_ij_2}
{\frac{{ \sum\limits_{t = 1}^{p} \sum\limits_{l = 1}^{n} \det
{\left( {{\rm {\bf A}}^{ *} {\rm {\bf A}}} \right)_{\,. \,i}
\left( {{\rm {\bf e}}_{.\,l} } \right)}
}\,\,\widetilde{d}_{lt}{\sum\limits_{\alpha \in I_{r_{2},p}
{\left\{ {j} \right\}}} { \left| {({\rm {\bf B}}{\rm {\bf B}}^{ *}
)_{j\,.\,} ({\rm {\bf e}}_{t.} )\,_{\alpha} ^{\alpha}} \right| } }
}{{{\det {\left( {{\rm {\bf A}}^{ *} {\rm {\bf A}}} \right)}}
}{{\sum\limits_{\alpha \in I_{r_{2},p}} {{\left| {\left( {{\rm
{\bf B}}{\rm {\bf B}}^{ *} } \right){\kern 1pt} _{\alpha}
^{\alpha} } \right|}}} }}    }.
\end{equation}
If we substitute (\ref{eq:d^B_den}) in (\ref{eq:x_ij_2}), then we
get
\[x_{ij} ={\frac{{  \sum\limits_{l = 1}^{n}
\det {\left( {{\rm {\bf A}}^{ *} {\rm {\bf A}}} \right)_{\,. \,i}
\left( {{\rm {\bf e}}_{.\,l} } \right)} }\,\,d^{{\rm {\bf
B}}}_{lj} }{{{\det {\left( {{\rm {\bf A}}^{ *} {\rm {\bf A}}}
\right)}} }{{\sum\limits_{\alpha \in I_{r_{2},p}} {{\left| {\left(
{{\rm {\bf B}}{\rm {\bf B}}^{ *} } \right){\kern 1pt} _{\alpha}
^{\alpha} } \right|}}} }}    }.
\]
Since again $\sum\limits_{l = 1}^{n}{\rm {\bf e}}_{.l}
 d^{{\rm {\bf B}}}_{lj}={\rm {\bf
d}}^{{\rm {\bf B}}}_{.\,j}$, then it follows
(\ref{eq:AXB_detA*A_d^B}), where $ {{\rm {\bf d}}_{.\,j}^{{\rm
{\bf B}}}}$ is (\ref{eq:def_d^B_m}).

If we denote by
\[
 d^{{\rm {\bf A}}}_{it}:= \]
\[
\sum\limits_{l = 1}^{n} \det {\left( {{\rm {\bf A}}^{ *} {\rm {\bf
A}}} \right)_{\,. \,i} \left( {\widetilde{{\rm {\bf d}}}_{.\,t} }
\right)} = \sum\limits_{l = 1}^{n} \det {\left( {{\rm {\bf A}}^{
*} {\rm {\bf A}}} \right)_{\,. \,i} \left( {{\rm {\bf e}}_{.\,l} }
\right)} \,\,\widetilde{d}_{lt}
\]
the $t$th component  of a row-vector ${\rm {\bf d}}^{{\rm {\bf
A}}}_{i\,.}= (d^{{\rm {\bf A}}}_{i1},...,d^{{\rm {\bf A}}}_{ip})$
for all $t=\overline {1,p}$ and substitute it in
(\ref{eq:x_ij_2}), we obtain
\[x_{ij} ={\frac{{ \sum\limits_{t = 1}^{p}
 d^{{\rm {\bf A}}}_{it}
}{\sum\limits_{\alpha \in I_{r_{2},p} {\left\{ {j} \right\}}} {
\left| {({\rm {\bf B}}{\rm {\bf B}}^{ *} )_{j\,.\,} ({\rm {\bf
e}}_{t.} )\,_{\alpha} ^{\alpha}} \right| } } }{{{\det {\left(
{{\rm {\bf A}}^{ *} {\rm {\bf A}}} \right)}}
}{{\sum\limits_{\alpha \in I_{r_{2},p}} {{\left| {\left( {{\rm
{\bf B}}{\rm {\bf B}}^{ *} } \right){\kern 1pt} _{\alpha}
^{\alpha} } \right|}}} }}    }.
\]
Since again $\sum\limits_{t = 1}^{p}
 d^{{\rm {\bf A}}}_{it}{\rm {\bf e}}_{t.}={\rm {\bf
d}}^{{\rm {\bf A}}}_{i\,.}$, then it follows
(\ref{AXB_detA*A_d^A}), where ${\rm {\bf d}}^{{\rm {\bf
A}}}_{i\,.}$ is (\ref{eq:def_d^A}).

(iiii) The proof is similar to the proof of (iii).
$\blacksquare$

\subsection{Cramer's rule of the Drazin inverse solutions of some  matrix
equations}
Consider a  matrix equation
\begin{equation}\label{eq:AX_d=B}
 {\rm {\bf A}}{\rm {\bf X}} = {\rm {\bf B}},
\end{equation}
where ${\rm {\bf
A}}\in  {\mathbb{C}}^{n\times n} $ with $Ind{\kern 1pt} {\rm {\bf A}}= k $, ${\rm {\bf B}}\in  {\mathbb{C}}^{n\times m}$ are given and ${\rm {\bf X}} \in
{\mathbb{C}}^{n\times m}$ is unknown.
\begin{theorem}(\cite{xu}, Theorem 1) If  the range space $R({\bf B})\subset R({\bf A}^{k})$,
then the matrix equation (\ref{eq:AX_d=B}) with constrain
$R({\bf X})\subset R({\bf A}^{k})$ has a unique solution
\[ {\bf X}={\bf A}^{D}{\bf B}.\]
\end{theorem}
 We denote ${\rm {\bf A}}^{ k}{\rm {\bf B}}=:\hat{{\rm
{\bf B}}}= (\hat{b}_{ij})\in {\mathbb{C}}^{n\times m}$.
\begin{theorem}
If $\rank{\rm {\bf
A}}^{k + 1} = \rank{\rm {\bf A}}^{k}=r \le n$ for ${\rm {\bf A}}\in {\mathbb C}^{n\times n} $, then for Drazin inverse solution $ {\bf X}={\bf A}^{D}{\bf B}= (x_{ij})\in {\mathbb{C}}^{n\times m}$ of  (\ref{eq:AX_d=B}) we have for all $i=\overline {1,n} $, $j=\overline {1,m} $,
\begin{equation}
\label{eq:dr_AX} x_{ij} = {\frac{{{\sum\limits_{\beta \in
J_{r,\,n} {\left\{ {i} \right\}}} {{\left| \left( {{\rm
{\bf A}}^{ k+1}_{\,.\,i} \left( {\hat{{\rm
{\bf b}}}_{.j}} \right)} \right) {\kern 1pt} _{\beta}
^{\beta} \right|}} } }}{{{\sum\limits_{\beta \in J_{r,\,\,n}}
{{\left| {\left( {{\rm {\bf A}}^{ k+1} } \right){\kern
1pt} {\kern 1pt} _{\beta} ^{\beta} }  \right|}}} }}}.
\end{equation}
\end{theorem}

{\textit{Proof}}.  By Theorem  \ref{theor:dr_repr_row} we can represent the matrix
${\rm {\bf A}}^{ D} $ by (\ref{eq:dr_repr_col}). Therefore, we
obtain for all $i=\overline {1,n} $, $j=\overline {1,m} $,
\[
x_{ij}  =\sum_{s=1}^{n}
a_{is}^{D}b_{sj}=\sum_{s=1}^{n}{\frac{{{\sum\limits_{\beta \in
J_{r,\,n} {\left\{ {i} \right\}}} {{\left| \left( {{\rm
{\bf A}}^{ k+1}_{\,. \,i} \left( {{\rm {\bf
a}}_{\,.s}^{ (k)} }  \right)} \right) {\kern 1pt} _{\beta}
^{\beta}\right|} } } }}{{{\sum\limits_{\beta \in J_{r,\,n}}
{{\left| {\left( {{\rm {\bf A}}^{ k+1} } \right){\kern
1pt} _{\beta} ^{\beta} }  \right|}}} }}}\cdot b_{sj}=
\]
\[
{\frac{{{\sum\limits_{\beta \in J_{r,\,n} {\left\{ {i}
\right\}}}\sum_{s=1}^{n} {{\left| \left( {{\rm {\bf A}}^{
k+1} _{\,. \,i} \left( {{\rm {\bf a}}_{\,.s}^{
(k)} } \right)} \right) {\kern 1pt} _{\beta} ^{\beta}\right|} } }
}\cdot b_{sj}} {{{\sum\limits_{\beta \in J_{r,\,n}} {{\left|
{\left( {{\rm {\bf A}}^{ k+1} } \right){\kern 1pt}
_{\beta} ^{\beta} } \right|}}} }}}.
\]
Since ${\sum\limits_{s} {{\rm {\bf a}}_{\,.s}^{ (k)}  b_{sj}} }=
\left( {{\begin{array}{*{20}c}
 {{\sum\limits_{s} {a_{1s}^{(k)}  b_{sj}} } } \hfill \\
 {{\sum\limits_{s} {a_{2s}^{ (k)}  b_{sj}} } } \hfill \\
 { \vdots}  \hfill \\
 {{\sum\limits_{s} {a_{ns}^{(k)}  b_{sj}} } } \hfill \\
\end{array}} } \right) = \hat{{\rm {\bf b}}}_{.j}$, then it follows (\ref{eq:dr_AX}).
$\blacksquare$

Consider a  matrix equation
\begin{equation}\label{eq:XA_d=B}
 {\rm {\bf X}}{\rm {\bf A}} = {\rm {\bf B}},
\end{equation}
 where ${\rm {\bf
A}}\in{\mathbb{C}}^{m\times m}$ with  $Ind{\kern 1pt} {\rm {\bf A}} = k $, ${\rm {\bf B}}\in {\mathbb{C}}^{n\times m} $ are given and  ${\rm {\bf X}} \in
{\mathbb{C}}^{n\times m}$ is unknown.
\begin{theorem}(\cite{xu}, Theorem 2) If the null space $N({\bf B})\supset N({\bf A}^{k})$, then
the matrix equation (\ref{eq:XA_d=B}) with constrain
$N({\bf X})\supset N({\bf A}^{k})$ has a unique solution
\[ {\bf X}={\bf B}{\bf A}^{D}.\]
\end{theorem}

 We denote  ${\rm {\bf B}}{\rm {\bf A}}^{
k}=:\check{{\rm {\bf B}}}= (\check{b}_{ij})\in
{\mathbb{C}}^{n\times m}$.
\begin{theorem}
 If $\rank{\rm {\bf
A}}^{k + 1} = \rank{\rm {\bf A}}^{k}=r \le m$ for ${\rm {\bf A}}\in {\mathbb C}^{m\times m} $, then for Drazin inverse solution $ {\bf X}={\bf B}{\bf A}^{D}= (x_{ij})\in {\mathbb C}^{n\times m}$ of  (\ref{eq:XA_d=B}),  we have for all $i=\overline {1,n} $, $j=\overline {1,m} $,
\begin{equation}
\label{eq:dr_XA} x_{ij} = {\frac{{{\sum\limits_{\alpha \in I_{r,m}
{\left\{ {j} \right\}}} {{\left| \left( {{\rm
{\bf A}}^{ k+1} _{\,j\,.} \left( {\check{{\rm {\bf
b}}}_{i\,.}} \right)} \right)\,_{\alpha} ^{\alpha}\right|} }
}}}{{{\sum\limits_{\alpha \in I_{r,\,m}}  {{\left| {\left( {{\rm {\bf A}}^{ k+1} } \right) {\kern 1pt} _{\alpha}
^{\alpha} } \right|}}} }}}.
\end{equation}
\end{theorem}
{\textit{Proof}}.  By Theorem  \ref{theor:dr_repr_row} we can represent the matrix
${\rm {\bf A}}^{ D} $ by (\ref{eq:dr_repr_row}). Therefore, for
all $i=\overline {1,n} $, $j=\overline {1,m} $,  we obtain
\[
x_{ij}  =\sum_{s=1}^{m}b_{is} a_{sj}^{D}=\sum_{s=1}^{m}b_{is}\cdot
{\frac{{{\sum\limits_{\alpha \in I_{r,m} {\left\{ {j} \right\}}}
{{\left| \left( {{\rm {\bf A}}^{ k+1}_{\,j\,.} \left( {\bf a}_{s\,.}^{(k)} \right)}
\right)\,_{\alpha} ^{\alpha}\right|} } }}}{{{\sum\limits_{\alpha
\in I_{r,\,m}}  {{\left| {\left( {{\rm {\bf A}}^{ k+1}
} \right) {\kern 1pt} _{\alpha} ^{\alpha} } \right|}}} }}} =
\]
\[
{\frac{{{\sum_{s=1}^{m}b_{ik} \sum\limits_{\alpha \in I_{r,m}
{\left\{ {j} \right\}}} {{\left| \left( {{\rm
{\bf A}}^{k+1}_{\,j\,.} \left( {\bf a}_{s\,.}^{(k)}
\right)} \right)\,_{\alpha} ^{\alpha}\right|} } } }}
{{{\sum\limits_{\alpha \in I_{r,\,m}}  {{\left| {\left( {{\rm {\bf A}}^{ k+1} } \right) {\kern 1pt} _{\alpha} ^{\alpha} }
\right|}}} }}}
\]
Since for all $i = \overline {1,n} $
 \[{\sum\limits_{s} {{  b_{is}\rm {\bf a}}_{s\,.}^{ (k)}}
}=\begin{pmatrix}
  \sum\limits_{s} {b_{is}a_{s1}^{(k)} } & \sum\limits_{s} {b_{is}a_{s2}^{ (k)} } & \cdots & \sum\limits_{s} {b_{is}a_{sm}^{ (k)} }
\end{pmatrix}
 = \check{{\rm {\bf b}}}_{i.},\]
then it follows (\ref{eq:dr_XA}).
$\blacksquare$

 Consider a  matrix equation
 \begin{equation}\label{eq:AXB_d=D}
 {\rm {\bf A}}{\rm {\bf X}}{\rm {\bf B}} = {\rm {\bf
D}},
\end{equation}
where  ${\rm {\bf
A}}\in{\mathbb{C}}^{n\times n}$ with  $Ind{\kern 1pt} {\rm {\bf A}} = k_{1} $,  ${\rm {\bf B}}\in {\mathbb{C}}^{m\times m} $ with  $Ind{\kern 1pt} {\rm {\bf B}} = k_{2} $ and $ {\rm {\bf D}}\in{\rm
{\mathbb{C}}}^{n \times m}$ are given, and  $ {\rm {\bf X}}\in{\rm
{\mathbb{C}}}^{n \times m}$ is unknown.
\begin{theorem}(\cite{xu}, Theorem 3) If $R({\bf D})\subset R({\bf A}^{k_{1}})$ and $N({\bf D})\supset N({\bf B}^{k_{2}})$, $k=max\{k_{1}, k_{2}\}$, then the matrix equation (\ref{eq:AXB_d=D}) with constrain
$R({\bf X})\subset R({\bf A}^{k})$ and $N({\bf X})\supset N({\bf B}^{k})$ has a unique solution
\[{\rm {\bf X}}={\rm {\bf A}}^{D}{\rm {\bf D}}{\rm {\bf B}}^{
D}.\]
\end{theorem}

We denote  ${\rm {\bf A}}^{
 k_{1} }{\rm {\bf D}}{\rm {\bf B}}^{
 k_{2} }=:\widetilde{{\rm {\bf D}}}= (\widetilde{d}_{ij})\in
{\mathbb{C}}^{n\times m}$.
\begin{theorem}\label{theor:AXB_d=D}
 If  $\rank{\rm {\bf
A}}^{k_{1} + 1} = \rank{\rm {\bf A}}^{k_{1}}=r_{1} \le n$ for ${\rm {\bf A}}\in {\mathbb C}^{n\times n} $, and  $\rank{\rm {\bf
B}}^{k_{2} + 1} = \rank{\rm {\bf B}}^{k_{2}}=r_{2} \le m$ for ${\rm {\bf B}}\in {\mathbb C}^{m\times m} $, then for the
Drazin inverse solution ${\rm {\bf X}}={\rm {\bf A}}^{D}{\rm {\bf D}}{\rm {\bf B}}^{
D}=:(x_{ij})\in
{\mathbb{C}}^{n\times m}$  of (\ref{eq:AXB_d=D}) we have
\begin{equation}\label{eq:d^B_d}
x_{ij} = {\frac{{{\sum\limits_{\beta \in J_{r_{1},\,n} {\left\{
{i} \right\}}} { \left| {{\rm {\bf A}}^{ k_{1}+1} _{\,.\,i} \left( {{{\rm {\bf d}}}\,_{.\,j}^{{\rm {\bf B}}}}
\right)\, _{\beta} ^{\beta}} \right| } } }}{{{\sum\limits_{\beta
\in J_{r_{1},n}} {{\left| {\left( {{\rm {\bf A}}^{ k_{1}+1} } \right)_{\beta} ^{\beta} } \right|}} \sum\limits_{\alpha \in
I_{r_{2},m}}{{\left| {\left( {{\rm {\bf B}}^{k_{2}+1} }
\right) _{\alpha} ^{\alpha} } \right|}}} }}},
\end{equation}
or
\begin{equation}\label{eq:d^A_d}
 x_{ij}={\frac{{{\sum\limits_{\alpha
\in I_{r_{2},m} {\left\{ {j} \right\}}} { \left| {{\rm {\bf B}}^{k_{2}+1} _{j\,.} \left( {{{\rm {\bf
d}}}\,_{i\,.}^{{\rm {\bf A}}}} \right)\,_{\alpha} ^{\alpha}}
\right| } }}}{{{\sum\limits_{\beta \in J_{r_{1},n}} {{\left|
{\left( {{\rm {\bf A}}^{k_{1}+1} } \right) _{\beta}
^{\beta} } \right|}}\sum\limits_{\alpha \in I_{r_{2},m}} {{\left|
{\left( {{\rm {\bf B}}^{k_{2}+1} } \right) _{\alpha}
^{\alpha} } \right|}}} }}},
\end{equation}
where
  \begin{equation}
\label{eq:def_d^B_m_d}
   {{\rm {\bf d}}_{.\,j}^{{\rm {\bf B}}}}=\left[
\sum\limits_{\alpha \in I_{r_{2},m} {\left\{ {j} \right\}}} {
\left| {{\rm {\bf B}}^{ k_{2}+1} _{j.}
\left( {\widetilde{{\rm {\bf d}}}_{1.}} \right)\,_{\alpha} ^{\alpha}}
\right|},...,\sum\limits_{\alpha \in I_{r_{2},m} {\left\{ {j}
\right\}}} { \left| {{\rm {\bf B}}^{k_{2}+1}_{j.} \left( {\widetilde{{\rm {\bf d}}}_{n.}}
\right)\,_{\alpha} ^{\alpha}} \right|} \right]^{T},
\end{equation}
  \[
  {{\rm {\bf d}}_{i\,.}^{{\rm {\bf A}}}}=\left[
\sum\limits_{\beta \in J_{r_{1},n} {\left\{ {i} \right\}}} {
\left| {{\rm {\bf A}}^{k_{1}+1}_{.i}
\left( {\widetilde{{\rm {\bf d}}}_{.1}} \right)\,_{\beta} ^{\beta}}
\right|},...,\sum\limits_{\alpha \in I_{r_{1},n} {\left\{ {i}
\right\}}} { \left| {{\rm {\bf A}}^{k_{1}+1}_{.i} \left( {\widetilde{{\rm {\bf d}}}_{.\,m}}
\right)\,_{\beta} ^{\beta}} \right|} \right]
\]
 are the column-vector and the row-vector. ${\widetilde{{\rm {\bf d}}}_{i.}}$  and ${\widetilde{{\rm {\bf d}}}_{.j}}$ are respectively the $i$th row  and the $j$th column of $\widetilde{{\rm {\bf D}}}$ for all $i=\overline {1,n}$, $j=\overline {1,m}$.

\end{theorem}
{\textit{Proof}}.
By (\ref{eq:dr_repr_col}) and
(\ref{eq:dr_repr_row}) the Drazin inverses ${\rm {\bf
A}}^{D} = \left( {a_{ij}^{D} } \right) \in {\rm
{\mathbb{C}}}_{}^{n\times n} $ and ${\rm {\bf B}}^{D} = \left(
{b_{ij}^{D} } \right) \in {\rm {\mathbb{C}}}^{m\times m} $
possess the following determinantal representations, respectively,
\[
 a_{ij}^{D}  = {\frac{{{\sum\limits_{\beta
\in J_{r_{1},\,n} {\left\{ {i} \right\}}} { \left| {{\rm
{\bf A}}^{k_{1}+1} _{\,. \,i} \left( {{\rm {\bf
a}}_{.j}^{ (k_{1})} }  \right){\kern 1pt}  _{\beta} ^{\beta}} \right| }
} }}{{{\sum\limits_{\beta \in J_{r_{1},\,n}} {{\left| {\left(
{{\rm {\bf A}}^{{k_{1}+1}} } \right){\kern 1pt} _{\beta}
^{\beta} }  \right|}}} }}},
\]
\begin{equation}\label{eq:b_d}
 b_{ij}^{D}  =
{\frac{{{\sum\limits_{\alpha \in I_{r_{2},m} {\left\{ {j}
\right\}}} { \left| {{\rm {\bf B}}^{k_{2}+1} _{j\,.\,}
({\rm {\bf b}}_{i.\,}^{ (k_{2})} )\,_{\alpha} ^{\alpha}} \right| }
}}}{{{\sum\limits_{\alpha \in I_{r_{2},m}}  {{\left| {\left( {{\rm {\bf B}}^{k_{2}+1} } \right){\kern 1pt} _{\alpha}
^{\alpha} } \right|}}} }}}.
\end{equation}
Then an entry of the Drazin inverse solution ${\rm {\bf X}}={\rm {\bf A}}^{D}{\rm {\bf D}}{\rm {\bf B}}^{
D}=:(x_{ij})\in
{\mathbb{C}}^{n\times m}$ is

\begin{equation}
\label{eq:sum_d} x_{ij} = {{\sum\limits_{s = 1}^{m} {\left(
{{\sum\limits_{t = 1}^{n} {{a}_{it}^{D} d_{ts}} } } \right)}}
{b}_{sj}^{D}}.
\end{equation}
Denote  by $\hat{{\rm {\bf d}}_{.s}}$ the $s$th column of ${\rm
{\bf A}}^{ k}{\rm {\bf D}}=:\hat{{\rm {\bf D}}}=
(\hat{d}_{ij})\in {\mathbb{C}}^{n\times m}$ for all $s=\overline
{1,m}$. It follows from ${\sum\limits_{t} { {\rm {\bf a}}_{.\,t}^{
D}}d_{ts} }=\hat{{\rm {\bf d}}_{.\,s}}$ that
\[
\sum\limits_{t = 1}^{n} {{a}_{it}^{D} d_{ts}}=\sum\limits_{t =
1}^{n}{\frac{{{\sum\limits_{\beta \in J_{r_{1},\,n} {\left\{ {i}
\right\}}} { \left| {{\rm {\bf A}}^{k_{1}+1} _{\,. \,i} \left( {{\rm {\bf a}}_{.t}^{
(k_{1})} } \right) {\kern
1pt} _{\beta} ^{\beta}} \right| } } }}{{{\sum\limits_{\beta \in
J_{r_{1},\,n}} {{\left| {\left( {{\rm {\bf A}}^{k_{1}+1} }
\right){\kern 1pt} _{\beta} ^{\beta} }  \right|}}} }}}\cdot
d_{ts}=
\]
\begin{equation}\label{eq:sum_d_cdet}
{\frac{{{\sum\limits_{\beta \in J_{r_{1},\,n} {\left\{ {i}
\right\}}}\sum\limits_{t = 1}^{n} { \left| {{\rm {\bf
A}}^{k_{1}+1}_{\,. \,i} \left( {{\rm {\bf
a}}_{.t}^{(k_{1})} } \right) {\kern 1pt} _{\beta} ^{\beta}} \right| } }
}\cdot d_{ts}}{{{\sum\limits_{\beta \in J_{r_{1},\,n}} {{\left|
{\left( {{\rm {\bf A}}^{k_{1}+1}} \right){\kern 1pt}
_{\beta} ^{\beta} }  \right|}}} }}}={\frac{{{\sum\limits_{\beta
\in J_{r_{1},\,n} {\left\{ {i} \right\}}} { \left| {{\rm
{\bf A}}^{k_{1}+1}_{\,. \,i} \left( \hat{{\rm
{\bf d}}_{.\,s}} \right) {\kern 1pt} _{\beta} ^{\beta}} \right| }
} }}{{{\sum\limits_{\beta \in J_{r_{1},\,n}} {{\left| {\left(
{{\rm {\bf A}}^{k_{1}+1}} \right){\kern 1pt} _{\beta}
^{\beta} }  \right|}}} }}}
\end{equation}

Substituting  (\ref{eq:sum_d_cdet}) and (\ref{eq:b_d}) in
(\ref{eq:sum_d}), we obtain
\[
x_{ij} =\sum\limits_{s = 1}^{m}{\frac{{{\sum\limits_{\beta \in
J_{r_{1},\,n} {\left\{ {i} \right\}}} { \left| {{\rm {\bf
A}}^{k_{1}+1} _{\,. \,i} \left( \hat{{\rm {\bf
d}}_{.\,s}} \right) {\kern 1pt} _{\beta} ^{\beta}} \right| } }
}}{{{\sum\limits_{\beta \in J_{r_{1},\,n}} {{\left| {\left( {{\rm
{\bf A}}^{k_{1}+1}} \right){\kern 1pt} _{\beta} ^{\beta}
}  \right|}}} }}}{\frac{{{\sum\limits_{\alpha \in I_{r_{2},m}
{\left\{ {j} \right\}}} { \left| {{\rm {\bf B}}^{k_{2}+1}
_{j\,.\,} ({\rm {\bf b}}_{s.\,}^{ (k_{2})} )\,_{\alpha} ^{\alpha}}
\right| } }}}{{{\sum\limits_{\alpha \in I_{r_{2},m}}  {{\left|
{\left( {{\rm {\bf B}}^{k_{2}+1} } \right){\kern 1pt}
_{\alpha} ^{\alpha} } \right|}}} }}}.
\]
Suppose ${\rm {\bf e}}_{s.}$ and ${\rm {\bf e}}_{.\,s}$ are
respectively the unit row-vector and the unit column-vector whose
components are $0$, except the $s$th components, which are $1$. Since
\[
\hat{{\rm {\bf d}}_{.\,s}}=\sum\limits_{l = 1}^{n}{\rm {\bf
e}}_{.\,l}\hat{ d_{ls}},\,\,\,\,  {\rm {\bf b}}_{s.\,}^{
(k_{2})}=\sum\limits_{t = 1}^{m}b_{st}^{(k_{2})}{\rm {\bf e}}_{t.},\,\,\,\,
\sum\limits_{s=1}^{m}\hat{d_{ls}}b_{st}^{(k_{2})}=\widetilde{d}_{lt},
\]
then we have
\[
x_{ij} = \]
\[{\frac{{ \sum\limits_{s = 1}^{m}\sum\limits_{t =
1}^{m} \sum\limits_{l = 1}^{n} {\sum\limits_{\beta \in
J_{r_{1},\,n} {\left\{ {i} \right\}}} { \left| {{\rm {\bf
A}}^{k_{1}+1}_{\,. \,i} \left( {\rm {\bf
e}}_{.\,l} \right) {\kern 1pt} _{\beta} ^{\beta}} \right| } }
}\hat{ d_{ls}}b_{st}^{(k_{2})}{\sum\limits_{\alpha \in I_{r_{2},m}
{\left\{ {j} \right\}}} { \left| {{\rm {\bf B}}^{k_{2}+1}
_{j\,.\,} ({\rm {\bf e}}_{t.} )\,_{\alpha} ^{\alpha}} \right| } }
}{{{\sum\limits_{\beta \in J_{r_{1},\,n}} {{\left| {\left( {{\rm
{\bf A}}^{k_{1}+1}} \right){\kern 1pt} _{\beta} ^{\beta}
}  \right|}}} }{{\sum\limits_{\alpha \in I_{r_{2},m}} {{\left|
{\left( {{\rm {\bf B}}^{k_{2}+1} } \right){\kern 1pt}
_{\alpha} ^{\alpha} } \right|}}} }}    }=
\]
\begin{equation}\label{eq:x_d_ij}
{\frac{{ \sum\limits_{t = 1}^{m} \sum\limits_{l = 1}^{n}
{\sum\limits_{\beta \in J_{r_{1},\,n} {\left\{ {i} \right\}}}
{\left| {{\rm {\bf A}}^{k_{1}+1}_{\,.
\,i} \left( {\rm {\bf e}}_{.\,l} \right) {\kern 1pt} _{\beta}
^{\beta}} \right| } } }\,\,\widetilde{d}_{lt}{\sum\limits_{\alpha
\in I_{r_{2},m} {\left\{ {j} \right\}}} { \left| {{\rm {\bf B}}^{k_{2}+1}_{j\,.\,} ({\rm {\bf e}}_{t.} )\,_{\alpha}
^{\alpha}} \right| } } }{{{\sum\limits_{\beta \in J_{r_{1},\,n}}
{{\left| {\left( {{\rm {\bf A}}^{k_{1}+1}} \right){\kern
1pt} _{\beta} ^{\beta} }  \right|}}} }{{\sum\limits_{\alpha \in
I_{r_{2},m}} {{\left| {\left( {{\rm {\bf B}}^{k_{2}+1} }
\right){\kern 1pt} _{\alpha} ^{\alpha} } \right|}}} }}    }.
\end{equation}
Denote by
\[
 d^{{\rm {\bf A}}}_{it}:= \]
\[
{\sum\limits_{\beta \in J_{r_{1},\,n} {\left\{ {i} \right\}}} {
\left| {{\rm {\bf A}}^{k_{1}+1}_{\,.
\,i} \left( \widetilde{{\rm {\bf d}}}_{.\,t} \right){\kern 1pt}
_{\beta} ^{\beta}} \right| } }= \sum\limits_{l = 1}^{n}
{\sum\limits_{\beta \in J_{r_{1},\,n} {\left\{ {i} \right\}}} {
\left| {{\rm {\bf A}}^{k_{1}+1}_{\,.
\,i} \left( {\rm {\bf e}}_{.\,l} \right){\kern 1pt}
 _{\beta} ^{\beta}} \right| } } \widetilde{d}_{lt}
\]
the $t$th component  of a row-vector ${\rm {\bf d}}^{{\rm {\bf
A}}}_{i\,.}= (d^{{\rm {\bf A}}}_{i1},...,d^{{\rm {\bf A}}}_{im})$
for all $t=\overline {1,m}$. Substituting it in (\ref{eq:x_d_ij}),
we obtain
\[x_{ij} ={\frac{{ \sum\limits_{t = 1}^{m}
 d^{{\rm {\bf A}}}_{it}
}{\sum\limits_{\alpha \in I_{r_{2},m} {\left\{ {j} \right\}}} {
\left| {{\rm {\bf B}}^{k_{2}+1} _{j\,.\,} ({\rm {\bf
e}}_{t.} )\,_{\alpha} ^{\alpha}} \right| } }
}{{{\sum\limits_{\beta \in J_{r_{1},\,n}} {{\left| {\left( {{\rm
{\bf A}}^{k_{1}+1}} \right){\kern 1pt} _{\beta} ^{\beta}
}  \right|}}} }{{\sum\limits_{\alpha \in I_{r_{2},m}} {{\left|
{\left( {{\rm {\bf B}}^{k_{2}+1} } \right){\kern 1pt}
_{\alpha} ^{\alpha} } \right|}}} }}    }.
\]
Since $\sum\limits_{t = 1}^{m}
 d^{{\rm {\bf A}}}_{it}{\rm {\bf e}}_{t.}={\rm {\bf
d}}^{{\rm {\bf A}}}_{i\,.}$, then it follows (\ref{eq:d^A_d}).

If we denote by
\[
 d^{{\rm {\bf B}}}_{lj}:=
\sum\limits_{t = 1}^{m}\widetilde{d}_{lt}{\sum\limits_{\alpha \in
I_{r_{2},m} {\left\{ {j} \right\}}} { \left| {{\rm
{\bf B}}^{k_{2}+1} _{j\,.\,} ({\rm {\bf e}}_{t.} )\,_{\alpha}
^{\alpha}} \right| } }={\sum\limits_{\alpha \in I_{r_{2},m}
{\left\{ {j} \right\}}} {\left| {{\rm {\bf B}}^{k_{2}+1}
_{j\,.\,} (\widetilde{{\rm {\bf d}}}_{l.} )\,_{\alpha} ^{\alpha}}
\right| } }
\]
the $l$th component  of a column-vector ${\rm {\bf d}}^{{\rm {\bf
B}}}_{.\,j}= (d^{{\rm {\bf B}}}_{1j},...,d^{{\rm {\bf
B}}}_{jn})^{T}$ for all $l=\overline {1,n}$ and substitute it in
(\ref{eq:x_d_ij}), we obtain
\[x_{ij} ={\frac{{  \sum\limits_{l = 1}^{n}
{\sum\limits_{\beta \in J_{r_{1},\,n} {\left\{ {i} \right\}}} {
\left| {{\rm {\bf A}}^{k_{1}+1}_{\,.
\,i} \left( {\rm {\bf e}}_{.\,l} \right){\kern 1pt} _{\beta}
^{\beta}} \right| } } }\,\,d^{{\rm {\bf B}}}_{lj}
}{{{\sum\limits_{\beta \in J_{r_{1},\,n}} {{\left| {\left( {{\rm
{\bf A}}^{k_{1}+1}} \right){\kern 1pt} _{\beta} ^{\beta}
}  \right|}}} }{{\sum\limits_{\alpha \in I_{r_{2},m}} {{\left|
{\left( {{\rm {\bf B}}^{k_{2}+1} } \right){\kern 1pt}
_{\alpha} ^{\alpha} } \right|}}} }}    }.
\]
Since $\sum\limits_{l = 1}^{n}{\rm {\bf e}}_{.l}
 d^{{\rm {\bf B}}}_{lj}={\rm {\bf
d}}^{{\rm {\bf B}}}_{.\,j}$, then it follows (\ref{eq:d^B_d}).
$\blacksquare$

\subsection{Examples}
In this subsection, we give an example to illustrate  results obtained in the section.

1. Let
us consider the matrix equation
\begin{equation}\label{eq_ex:AXB=D}
 {\rm {\bf A}}{\rm {\bf X}}{\rm {\bf B}} = {\rm {\bf
D}},
\end{equation}
where
\[{\bf A}=\begin{pmatrix}
  1 & i & i \\
  i & -1 & -1 \\
  0 & 1 & 0 \\
  -1 & 0 & -i
\end{pmatrix},\,\, {\bf B}=\begin{pmatrix}
  i & 1 & -i \\
  -1 & i & 1
\end{pmatrix},\,\, {\bf D}=\begin{pmatrix}
  1 & i & 1 \\
  i & 0 & 1 \\
  1 & i & 0 \\
  0 & 1 & i
\end{pmatrix}.\]
Since ${\rm rank}\,{\rm {\bf A}} = 2$ and ${\rm rank}\,{\rm {\bf
B}} = 1$, then we have the case (ii) of Theorem \ref{theor:AXB=D}.
We shall find the least squares solution of (\ref{eq_ex:AXB=D}) by
(\ref{eq:d^B}). Then we have

\[{\bf A}^{*}{\bf A}=\begin{pmatrix}
  3 & 2i & 3i \\
  -2i & 3 & 2 \\
  -3i & 2 & 3
\end{pmatrix},\,\, {\bf B}{\bf B}^{*}=\begin{pmatrix}
  3 & -3i \\
  3i & 3
\end{pmatrix},\,\, {\rm {\bf \widetilde{D}}}= {\rm {\bf
A}}^\ast{\rm {\bf D}}{\rm {\bf B}}^{\ast}=\begin{pmatrix}
  1 & -i \\
  -i & -1 \\
  -i & -1
\end{pmatrix},\] and $
{{{\sum\limits_{\alpha \in I_{1,\,2}}  {{\left| {\left( {{\rm {\bf
B}}{\rm {\bf B}}^{ *} } \right) {\kern 1pt} _{\alpha} ^{\alpha} }
\right|}}} }}=3+3=6,$
\[ {{{\sum\limits_{\beta \in J_{2,\,3}}
{{\left| {\left( {{\rm {\bf A}}^{ *}{\rm {\bf A}} } \right) {\kern
1pt} _{\beta} ^{\beta} } \right|}}} }}=\det\begin{pmatrix}
  3 & 2i \\
  -2i & 3
\end{pmatrix}+\det\begin{pmatrix}
  3 & 2 \\
  2 & 3
\end{pmatrix}+\det\begin{pmatrix}
  3 & 3i \\
  -3i & 3
\end{pmatrix}=10.\]
By (\ref{eq:def_d^B}), we can get \[{\bf d}_{.1}^{{\bf
B}}=\begin{pmatrix}
  1 \\
  -i \\
  -i
\end{pmatrix},\,\,\,\,\,{\bf d}_{.2}^{{\bf  B}}=\begin{pmatrix}
  -i \\
  -1 \\
  -1
\end{pmatrix}.\]
Since $\left( {{\rm {\bf A}}^{ *} {\rm {\bf A}}} \right)_{\,.\,1}
\left( {{{\rm {\bf d}}}\,_{.\,1}^{{\rm {\bf B}}}} \right)=
\begin{pmatrix}
  1 & 2i & 3i \\
  -i & 3 & 2 \\
  -i & 2 & 3
\end{pmatrix}$, then finally we obtain
\[
x_{11} = {\frac{{{\sum\limits_{\beta \in J_{2,\,3} {\left\{ {i}
\right\}}} { \left| {\left( {{\rm {\bf A}}^{ *} {\rm {\bf A}}}
\right)_{\,.\,1} \left( {{{\rm {\bf d}}}\,_{.\,1}^{{\rm {\bf B}}}}
\right)\, _{\beta} ^{\beta}} \right| } } }}{{{\sum\limits_{\beta
\in J_{2,3}} {{\left| {\left( {{\rm {\bf A}}^{ *} {\rm {\bf A}}}
\right)_{\beta} ^{\beta} } \right|}} \sum\limits_{\alpha \in
I_{1,2}}{{\left| {\left( {{\rm {\bf B}}{\rm {\bf B}}^{ *} }
\right) _{\alpha} ^{\alpha} } \right|}}} }}}=\frac{\det\begin{pmatrix}
  1 & 2i \\
  -i & 3
\end{pmatrix}+\det\begin{pmatrix}
  1 & 3i \\
  -i & 3
\end{pmatrix}}{60}=-\frac{1}{60}.\]
Similarly,
\[
  x_{12} =\frac{\det\begin{pmatrix}
  -i & 2i \\
  -1 & 3
\end{pmatrix}+\det\begin{pmatrix}
  -i & 3i \\
  -1 & 3
\end{pmatrix}}{60}=-\frac{i}{60},\]    \[x_{21} =\frac{\det\begin{pmatrix}
  3 & 1 \\
  -2i & -i
\end{pmatrix}+\det\begin{pmatrix}
  -i & 2 \\
  -i & 3
\end{pmatrix}}{60}=-\frac{2i}{60},\]
\[
   x_{22} =\frac{\det\begin{pmatrix}
  3 & -i \\
  -2i & -1
\end{pmatrix}+\det\begin{pmatrix}
  -1 & 2 \\
  -1 & 3
\end{pmatrix}}{60}=-\frac{2}{60},\] \[   x_{31} =\frac{\det\begin{pmatrix}
  3 & 1 \\
  -3i & -i
\end{pmatrix}+\det\begin{pmatrix}
  3 & -i \\
  2 & -i
\end{pmatrix}}{60}=-\frac{i}{60},\]
\[
    x_{32} =\frac{\det\begin{pmatrix}
  3 & -i \\
  -3i & -1
\end{pmatrix}+\det\begin{pmatrix}
  3 & -1 \\
  2 & -1
\end{pmatrix}}{60}=-\frac{1}{60}.
\]

 2. Let
us consider the matrix equation (\ref{eq_ex:AXB=D}),
where
\[{\bf A}=\begin{pmatrix}
  2 & 0 & 0 \\
  -i & i & i \\
  -i & -i & -i
\end{pmatrix},\,\, {\bf B}=\begin{pmatrix}
  1 & -1 & 1 \\
  i & -i & i \\
  -1 & 1 & 2
\end{pmatrix},\,\, {\bf D}=\begin{pmatrix}
  1 & i & 1 \\
  i & 0 & 1 \\
  1 & i & 0
\end{pmatrix}.\]
We shall find the Drazin inverse solution of (\ref{eq_ex:AXB=D}) by
(\ref{eq:d^B}). We obtain
\[{\bf A}^{2}=\begin{pmatrix}
  4 & 0 & 0 \\
  2-2i & 0 & 0 \\
  -2-2i & 0 & 0
\end{pmatrix},\,{\bf A}^{3}=\begin{pmatrix}
  8 & 0 & 0 \\
  4-4i & 0 & 0 \\
  -4-4i & 0 & 0
\end{pmatrix},\] \[ {\bf B}^{2}=\begin{pmatrix}
  -i & i & 3-i \\
  1 & -1 & 1+3i \\
  -3+i & 3-i & 3+i
\end{pmatrix}.\]
Since ${\rm rank}\,{\rm {\bf A}} = 2$ and ${\rm rank}\,{\rm {\bf A}}^{2}= {\rm rank}\,{\rm {\bf A}}^{2} = 1$, then $k_{1}={\rm Ind}\,{\rm {\bf A}}=2$ and $r_{1}=1$. Since ${\rm rank}\,{\rm {\bf B}} ={\rm rank}\,{\rm {\bf B}}^{2} = 2$, then $k_{2}={\rm Ind}\,{\rm {\bf B}}=1$ and $r_{2}=2$.
 Then we have
\[ {\rm {\bf \widetilde{D}}}= {\rm {\bf
A}}^{2}{\rm {\bf D}}{\rm {\bf B}}=\begin{pmatrix}
  -4 & 4 & 8\\
  -2+2i & 2-2i & 4-4i\\
  2+2i & -2-2i & -4-4i
  \end{pmatrix},\]
  and ${{{\sum\limits_{\beta \in J_{1,\,3}}
{{\left| {\left( {{\rm {\bf A}}^{ 3} } \right) {\kern
1pt} _{\beta} ^{\beta} } \right|}}} }}=
8+0+0=8,$
\[\begin{array}{l}
    {{{\sum\limits_{\alpha \in I_{2,\,3}}  {{\left| {\left( {{\rm {\bf B}}^{2} } \right) {\kern 1pt} _{\alpha} ^{\alpha} }
\right|}}} }}=\\ \det\begin{pmatrix}
  -i & i \\
  1 & -1
\end{pmatrix}+\det\begin{pmatrix}
  -1 & 1+3i \\
  3-i & 3+i
\end{pmatrix}+\det\begin{pmatrix}
  -i & 3-i \\
  -3+i & 3+i
\end{pmatrix}= \\
    0+(-9-9i)+(9-9i)=-18i.
  \end{array}
\]
By (\ref{eq:def_d^B_m}), we can get \[{\bf d}_{.1}^{{\bf
B}}=\begin{pmatrix}
  12-12i \\
  -12i \\
  -12
\end{pmatrix},\,\,{\bf d}_{.2}^{{\bf  B}}=\begin{pmatrix}
  -12+12i \\
  12i \\
  12
\end{pmatrix},\,\,{\bf d}_{.3}^{{\bf  B}}=\begin{pmatrix}
  8 \\
  -12-12i\\
  -12+12i
\end{pmatrix}.\]
Since ${\rm {\bf A}}^{ 3}_{\,.\,1}
\left( {{{\rm {\bf d}}}\,_{.\,1}^{{\rm {\bf B}}}} \right)=
\begin{pmatrix}
  12-12i & 0 & 0 \\
   -12i & 0 & 0 \\
  -12 & 0 & 0
\end{pmatrix}$, then finally we obtain
\[
x_{11} = {\frac{{{\sum\limits_{\beta \in J_{1,\,3} {\left\{ {1}
\right\}}} { \left| {{\rm {\bf A}}^{ 3}_{\,.\,1} \left( {{{\rm {\bf d}}}\,_{.\,1}^{{\rm {\bf B}}}}
\right)\, _{\beta} ^{\beta}} \right| } } }}{{{\sum\limits_{\beta
\in J_{1,3}} {{\left| {\left( {{\rm {\bf A}}^{ 3} }
\right)_{\beta} ^{\beta} } \right|}} \sum\limits_{\alpha \in
I_{2,3}}{{\left| {\left( {{\rm {\bf B}}^{2} }
\right) _{\alpha} ^{\alpha} } \right|}}} }}}=\frac{12-12i}{8\cdot (-18i)}=\frac{1+i}{12}.\]
Similarly,
\[
  x_{12} =\frac{-12+12i}{8\cdot (-18i)}=\frac{-1-i}{12},\,\,\,
  x_{13} =\frac{8}{8\cdot (-18i)}=\frac{i}{18},\]

   \[x_{21} =\frac{-12i}{8\cdot (-18i)}=\frac{1}{12},\,\,
   x_{22} =\frac{12i}{8\cdot (-18i)}=-\frac{1}{12},\,\,
  x_{23} =\frac{-12-12i}{8\cdot (-18i)}=\frac{1-i}{12},\]

   \[   x_{31} =\frac{12}{8\cdot (-18i)}=-\frac{i}{12},\,\,
    x_{32} =\frac{-12}{8\cdot (-18i)}=\frac{i}{12}.
\,\,
  x_{33} =\frac{-12+12i}{8\cdot (-18i)}=\frac{-1-i}{12}.\]
  Then
\[{\rm {\bf X}}=\left(
                                    \begin{array}{ccc}
                                      \frac{1+i}{12} & \frac{-1-i}{12} & \frac{i}{18} \\
                                      \frac{1}{12} & -\frac{1}{12} & \frac{1-i}{12} \\
                                      -\frac{i}{12} & \frac{i}{12} & \frac{-1-i}{12}\\
                                    \end{array}
                                  \right)
                                   \]
is the Drazin inverse solution of (\ref{eq_ex:AXB=D}).

\section{An application of the determinantal representations of the Drazin inverse to some differential matrix equations}
In this section we  demonstrate an  application of the determinantal representations of the Drazin inverse to following differential matrix equations, ${\bf X}'+ {\bf A}{\bf X}={\bf B}$ and ${\bf X}'+{\bf X}{\bf A}={\bf B}$, where the matrix ${\bf A}$ is singular.

Consider the
 matrix differential equation
 \begin{equation}\label{eq:left_dif}
 {\bf X}'+ {\bf A}{\bf X}={\bf B}
\end{equation}
where $ {\bf
A}\in{\rm {\mathbb{C}}}^{n \times n}$, $ {\bf B}\in{\rm
{\mathbb{C}}}^{n \times n}$ are given, $ {\rm {\bf X}}\in{\rm
{\mathbb{C}}}^{n \times n}$ is unknown. It's well-known that the general solution of
 (\ref{eq:left_dif})
is found to be
 \[
 {\bf X}(t)=\exp^{ -{\bf A}t}\left(\int \exp^{ {\bf A}t}dt\right){\bf B}
\]
If ${\bf A}$ is invertible, then
\[\int \exp^{ {\bf A}t}dt={\bf A}^{-1}\exp^{ {\bf A}t}+{\bf G},
\]
where ${\bf G}$ is an arbitrary $n\times n$ matrix. If ${\bf A}$ is singular, then the following theorem gives an answer.
 \begin{theorem}(\cite{ca1}, Theorem 1) \label{theor:int_to_draz}If ${\bf A}$ has index $k$, then
 \[
 \int \exp^{ {\bf A}t}dt={\bf A}^{D}\exp^{ {\bf A}t}+({\bf I}-{\bf A}{\bf A}^{D})t\left[{\bf I}+\frac{{\bf A}}{2}t+\frac{{\bf A}^{2}}{3!}t^{2}+...+\frac{{\bf A}^{k-1}}{k!}t^{k-1}\right]+{\bf G}. \]

\end{theorem}
Using Theorem \ref{theor:int_to_draz} and the power series expansion
of $\exp^{ -{\bf A}t}$, we get an explicit form for a general solution of (\ref{eq:left_dif})
 \[\begin{array}{l}
     {\bf X}(t)= \\
      \left\{{\bf A}^{D}+({\bf I}-{\bf A}{\bf A}^{D})t\left({\bf I}-\frac{{\bf A}}{2}t+\frac{{\bf A}^{2}}{3!}t^{2}-...(-1)^{k-1}\frac{{\bf A}^{k-1}}{k!}t^{k-1}\right)+{\bf G}\right\}{\bf B}.
   \end{array}
 \]
If we put ${ \bf G}={\bf 0}$, then we obtain the following  partial solution of (\ref{eq:left_dif}),

 \begin{equation}\label{eq:left_dif_part_sol}\begin{array}{c}
                                               {\bf X}(t)={\bf A}^{D}{\bf B}+({\bf B}-{\bf A}^{D}{\bf A}{\bf B})t-\frac{1}{2}({\bf A}{\bf B}-{\bf A}^{D}{\bf A}^{2}{\bf B})t^{2}+... \\
                                                \frac{(-1)^{k-1}}{k!}({\bf A}^{k-1}{\bf B}-{\bf A}^{D}{\bf A}^{k}{\bf B})t^{k}.
                                             \end{array}
\end{equation}

Denote ${\rm {\bf A}}^{l}{\rm {\bf B}}=:\widehat{{{\rm
{\bf B}}}}^{(l)}
= (\widehat{b}^{(l)}_{ij})\in {\mathbb{C}}^{n\times n}$ for all $l=\overline {1,2k}$.

 \begin{theorem} The partial solution (\ref{eq:left_dif_part_sol}), ${\bf X}(t)=(x_{ij})$, possess the following determinantal representation,
  \begin{equation}\label{eq:left_dif_repr}
\begin{array}{l}
   x_{ij}=
  {\frac{{{\sum\limits_{\beta \in
J_{r,\,n} {\left\{ {i} \right\}}} {{\left| \left( {{\rm
{\bf A}}^{k+1}_{\,.\,i} \left( {\widehat{{\rm
{\bf b}}}^{(k)}_{.j}} \right)} \right){\kern 1pt} {\kern 1pt} _{\beta}
^{\beta} \right|}} } }}{{{\sum\limits_{\beta \in J_{r,\,\,n}}
{{\left| {\left( {{\rm {\bf A}}^{ k+1} } \right){\kern
1pt} {\kern 1pt} _{\beta} ^{\beta} }  \right|}}} }}}+
 \left ({b_{ij}- {\frac{{{\sum\limits_{\beta \in
J_{r,\,n} {\left\{ {i} \right\}}} {{\left| \left( {{\rm
{\bf A}}^{k+1}_{\,.\,i} \left( {\widehat{{\rm
{\bf b}}}^{(k+1)}_{.j}} \right)} \right){\kern 1pt} {\kern 1pt} _{\beta}
^{\beta} \right|}} } }}{{{\sum\limits_{\beta \in J_{r,\,\,n}}
{{\left| {\left( {{\rm {\bf A}}^{k+1} } \right){\kern
1pt} {\kern 1pt} _{\beta} ^{\beta} }  \right|}}} }}}} \right)t\\
-\frac{1}{2} \left ({\widehat{b}^{(1)}_{ij}- {\frac{{{\sum\limits_{\beta \in
J_{r,\,n} {\left\{ {i} \right\}}} {{\left| \left( {{\rm
{\bf A}}^{k+1}_{\,.\,i} \left( {\widehat{{\rm
{\bf b}}}^{(k+2)}_{.j}} \right)} \right){\kern 1pt} {\kern 1pt} _{\beta}
^{\beta} \right|}} } }}{{{\sum\limits_{\beta \in J_{r,\,\,n}}
{{\left| {\left( {{\rm {\bf A}}^{k+1} } \right){\kern
1pt} {\kern 1pt} _{\beta} ^{\beta} }  \right|}}} }}}} \right)t^{2}+...\\\frac{(-1)^{k}}{k!}
\left ({\widehat{b}^{(k-1)}_{ij}- {\frac{{{\sum\limits_{\beta \in
J_{r,\,n} {\left\{ {i} \right\}}} {{\left| \left( {{\rm
{\bf A}}^{k+1}_{\,.\,i} \left( {\widehat{{\rm
{\bf b}}}^{(2k)}_{.j}} \right)} \right){\kern 1pt} {\kern 1pt} _{\beta}
^{\beta} \right|}} } }}{{{\sum\limits_{\beta \in J_{r,\,\,n}}
{{\left| {\left( {{\rm {\bf A}}^{k+1} } \right){\kern
1pt} {\kern 1pt} _{\beta} ^{\beta} }  \right|}}} }}}} \right)t^{k}
\end{array}
\end{equation}
for all $i,j=\overline {1,n}$.

\end{theorem}
{\textit{Proof}}.
Using the determinantal representation of the  identity ${\bf A}^{D}{\bf A}$ (\ref{eq:AdA}), we obtain the following determinantal representation of the matrix ${\bf A}^{D}{\bf A}^{m}{\bf B}:=(y_{ij})$,

\[y_{ij}=\sum\limits_{s = 1}^{n}p_{is}\sum\limits_{t = 1}^{n}a^{(m-1)}_{st}b_{tj}=
{\sum\limits_{\beta \in J_{r,n} {\left\{ {i} \right\}}}}\frac{{ \sum\limits_{s = 1}^{n}{{\left|
{\left( {{\rm {\bf A}}_{.\,i}^{k+1} \left( {{\rm {\bf
a}}_{.s}}^{(k+1)} \right)} \right)_{\beta} ^{\beta} } \right|}}}\cdot \sum\limits_{t = 1}^{n}a^{(m-1)}_{st}b_{tj}}{{\sum\limits_{\beta \in J_{r,n} } {{\left|
{\left( {{\rm {\bf A}}^{k+1} } \right)_{\beta} ^{\beta} } \right|}}}}=
\]
\[{\sum\limits_{\beta \in J_{r,n} {\left\{ {i} \right\}}}}\frac{{ \sum\limits_{t= 1}^{n}{{\left|
{\left( {{\rm {\bf A}}_{.\,i}^{k+1} \left( {{\rm {\bf
a}}_{.t}}^{(k+m)} \right)} \right)_{\beta} ^{\beta} } \right|}}}\cdot b_{tj}}{{\sum\limits_{\beta \in J_{r,n} } {{\left|
{\left( {{\rm {\bf A}}^{k+1} } \right)_{\beta} ^{\beta} } \right|}}}}=
  {\frac{{{\sum\limits_{\beta \in
J_{r,\,n} {\left\{ {i} \right\}}} {{\left| \left( {{\rm
{\bf A}}^{k+1}_{\,.\,i} \left( {\widehat{{\rm
{\bf b}}}^{(k+m)}_{.j}} \right)} \right){\kern 1pt} {\kern 1pt} _{\beta}
^{\beta} \right|}} } }}{{{\sum\limits_{\beta \in J_{r,\,\,n}}
{{\left| {\left( {{\rm {\bf A}}^{ k+1} } \right){\kern
1pt} {\kern 1pt} _{\beta} ^{\beta} }  \right|}}} }}}
\]
for all $i,j=\overline {1,n}$ and $m=\overline {1,k}$.
 From this and the determinantal representation of the Drazin inverse solution (\ref{eq:dr_AX}) and the identity  (\ref{eq:AdA}) it follows (\ref{eq:left_dif_repr}).
$\blacksquare$
\begin{corollary}  If $Ind {\bf A}=1$, then the partial solution of (\ref{eq:left_dif}),
\[ {\bf X}(t)=(x_{ij})={\bf A}^{g}{\bf B}+({\bf B}-{\bf A}^{g}{\bf A}{\bf B})t,\]
possess the following determinantal representation
 \begin{equation}\label{eq:left_dif_repr_gr}
x_{ij}= {\frac{{{\sum\limits_{\beta \in
J_{r,\,n} {\left\{ {i} \right\}}} {{\left| \left( {{\rm
{\bf A}}^{2}_{\,.\,i} \left( {\widehat{{\rm
{\bf b}}}_{.j}}^{(1)} \right)} \right){\kern 1pt} {\kern 1pt} _{\beta}
^{\beta} \right|}} } }}{{{\sum\limits_{\beta \in J_{r,\,\,n}}
{{\left| {\left( {{\rm {\bf A}}^{ 2} } \right){\kern
1pt} {\kern 1pt} _{\beta} ^{\beta} }  \right|}}} }}}+
 \left ({b_{ij}- {\frac{{{\sum\limits_{\beta \in
J_{r,\,n} {\left\{ {i} \right\}}} {{\left| \left( {{\rm
{\bf A}}^{2}_{\,.\,i} \left( {{\widehat{{\rm
{\bf b}}}}_{.j}}^{(2)} \right)} \right){\kern 1pt} {\kern 1pt} _{\beta}
^{\beta} \right|}} } }}{{{\sum\limits_{\beta \in J_{r,\,\,n}}
{{\left| {\left( {{\rm {\bf A}}^{ 2} } \right){\kern
1pt} {\kern 1pt} _{\beta} ^{\beta} }  \right|}}} }}}} \right)t.
\end{equation}
for all $i,j=\overline {1,n}$.
\end{corollary}

Consider the
 matrix differential equation
 \begin{equation}\label{eq:right_dif}
 {\bf X}'+ {\bf X}{\bf A}={\bf B}
\end{equation}
where $ {\bf
A}\in{\rm {\mathbb{C}}}^{n \times n}$, $ {\bf B}\in{\rm
{\mathbb{C}}}^{n \times n}$ are given, $ {\rm {\bf X}}\in{\rm
{\mathbb{C}}}^{n \times n}$ is unknown. The general solution of
 (\ref{eq:right_dif})
is found to be
 \[
 {\bf X}(t)={\bf B}\exp^{ -{\bf A}t}\left(\int \exp^{ {\bf A}t}dt\right)
\]
If ${\bf A}$ is singular, then an explicit form for a general solution of (\ref{eq:right_dif}) is
 \[\begin{array}{l}
     {\bf X}(t)=\\
     {\bf B}\left\{{\bf A}^{D}+({\bf I}-{\bf A}{\bf A}^{D})t\left({\bf I}-\frac{{\bf A}}{2}t+\frac{{\bf A}^{2}}{3!}t^{2}+...(-1)^{k-1}\frac{{\bf A}^{k-1}}{k!}t^{k-1}\right)+{\bf G}\right\}.
   \end{array}
 \]
If we put ${ \bf G}={\bf 0}$, then we obtain the following  partial solution of (\ref{eq:right_dif}),
  \begin{equation}\label{eq:right_dif_part_sol}\begin{array}{c}
                                                 {\bf X}(t)={\bf B}{\bf A}^{D}+({\bf B}-{\bf B}{\bf A}{\bf A}^{D})t-\frac{1}{2}({\bf B}{\bf A}-{\bf B}{\bf A}^{2}{\bf A}^{D})t^{2}+... \\
                                                  \frac{(-1)^{k-1}}{k!}({\bf B}{\bf A}^{k-1}-{\bf B}{\bf A}^{k}{\bf A}^{D})t^{k}.
                                               \end{array}
\end{equation}
Denote ${\rm {\bf B}}{\rm {\bf A}}^{l}=:\check{{{\rm
{\bf B}}}}^{(l)}
= (\check{b}^{(l)}_{ij})\in {\mathbb{C}}^{n\times n}$ for all $l=\overline {1,2k}$.
Using the determinantal representation of the Drazin inverse solution (\ref{eq:dr_XA}), the group inverse (\ref{eq:AAg}) and the identity  (\ref{eq:AAd})
we evidently obtain the following theorem.
 \begin{theorem} The partial solution (\ref{eq:right_dif_part_sol}), ${\bf X}(t)=(x_{ij})$, possess the following determinantal representation,
  \[
\begin{array}{l}
   x_{ij}=
  {\frac{{{\sum\limits_{\alpha \in I_{r,n}
{\left\{ {j} \right\}}} {{\left| \left( {{\rm
{\bf A}}^{k+1}_{j\,.} \left( {\check{{\rm
{\bf b}}}^{(k)}_{.\,i}} \right)} \right){\kern 1pt} {\kern 1pt} _{\alpha} ^{\alpha}  \right|}} } }}{{{\sum\limits_{\alpha \in I_{r,n} }
{{\left| {\left( {{\rm {\bf A}}^{ k+1} } \right){\kern
1pt} {\kern 1pt} _{\alpha} ^{\alpha}  }  \right|}}} }}}+
 \left ({b_{ij}- {\frac{{{\sum\limits_{\alpha \in I_{r,n}
{\left\{ {j} \right\}}} {{\left| \left( {{\rm
{\bf A}}^{k+1}_{j\,.} \left( {\check{{\rm
{\bf b}}}^{(k+1)}_{i\,.}} \right)} \right){\kern 1pt} {\kern 1pt} _{\alpha} ^{\alpha}  \right|}} } }}{{{\sum\limits_{\alpha \in I_{r,n} } {{\left| {\left( {{\rm {\bf A}}^{k+1} } \right){\kern
1pt} {\kern 1pt} _{\alpha} ^{\alpha}  }  \right|}}} }}}} \right)t \\
-\frac{1}{2} \left ({\check{b}^{(1)}_{ij}- {\frac{{{\sum\limits_{\alpha \in I_{r,n}
{\left\{ {j} \right\}}} {{\left| \left( {{\rm
{\bf A}}^{k+1}_{j\,.} \left( {\check{{\rm
{\bf b}}}^{(k+2)}_{i\,.}} \right)} \right){\kern 1pt} {\kern 1pt} _{\alpha} ^{\alpha}  \right|}} } }}{{{\sum\limits_{\alpha \in I_{r,n} } {{\left| {\left( {{\rm {\bf A}}^{k+1} } \right){\kern
1pt} {\kern 1pt} _{\alpha} ^{\alpha}  }  \right|}}} }}}} \right)t^{2}+...\\\frac{(-1)^{k}}{k!}
\left ({\check{b}^{(k-1)}_{ij}- {\frac{{{\sum\limits_{\alpha \in I_{r,n}
{\left\{ {j} \right\}}} {{\left| \left( {{\rm
{\bf A}}^{k+1}_{j\,.} \left( {\check{{\rm
{\bf b}}}^{(2k)}_{i\,.}} \right)} \right){\kern 1pt} {\kern 1pt} _{\alpha} ^{\alpha} \right|}} } }}{{{\sum\limits_{\alpha \in I_{r,n} }
{{\left| {\left( {{\rm {\bf A}}^{k+1} } \right){\kern
1pt} {\kern 1pt} _{\alpha} ^{\alpha} }  \right|}}} }}}} \right)t^{k}
\end{array}
\]
for all $i,j=\overline {1,n}$.
\end{theorem}
\begin{corollary}  If $Ind {\bf A}=1$, then the partial solution of (\ref{eq:right_dif}),
\[ {\bf X}(t)=(x_{ij})={\bf B}{\bf A}^{g}+({\bf B}-{\bf B}{\bf A}{\bf A}^{g})t,\]
possess the following determinantal representation
 \[
x_{ij}= {\frac{{{\sum\limits_{\alpha \in I_{r,n} {\left\{ {j} \right\}}} {{\left| \left( {{\rm
{\bf A}}^{2}_{j\,.} \left( {\widehat{{\rm
{\bf b}}}_{i\,.}}^{(1)}  \right)} \right){\kern 1pt} {\kern 1pt} _{\alpha}
^{\alpha} \right|}} } }}{{{\sum\limits_{\alpha \in I_{r,n}}
{{\left| {\left( {{\rm {\bf A}}^{ 2} } \right){\kern
1pt} {\kern 1pt} _{\alpha} ^{\alpha} }  \right|}}} }}}+
 \left ({b_{ij}- {\frac{{{\sum\limits_{\alpha \in I_{r,n} {\left\{ {j} \right\}}} {{\left| \left( {{\rm
{\bf A}}^{2}_{j\,.} \left( {{\widehat{{\rm
{\bf b}}}}_{i\,.}}^{(2)} \right)} \right){\kern 1pt} {\kern 1pt} _{\alpha}
^{\alpha} \right|}} } }}{{{\sum\limits_{\alpha  \in I_{r,n}}
{{\left| {\left( {{\rm {\bf A}}^{ 2} } \right){\kern
1pt} {\kern 1pt} _{\alpha}
^{\alpha} }  \right|}}} }}}} \right)t.
\]
for all $i,j=\overline {1,n}$.
\end{corollary}

\subsection{Example}
1.
 Let
us consider the differential matrix equation
\begin{equation}\label{eq_ex:X'=AX+B}
{\bf X}'+ {\rm {\bf A}} {\bf X} = {\rm {\bf
B}},
\end{equation}
where
\[ {\bf A}=\begin{pmatrix}
  1 & -1 & 1 \\
  i & -i & i \\
  -1 & 1 & 2
\end{pmatrix},\,\, {\bf B}=\begin{pmatrix}
  1 & i & 1 \\
  i & 0 & 1 \\
  1 & i & 0
\end{pmatrix}.\]
Since ${\rm rank}\,{\rm {\bf A}} ={\rm rank}\,{\rm {\bf A}}^{2} = 2$, then $k={\rm Ind}\,{\rm {\bf A}}=1$ and $r=2$.
The matrix ${\bf A}$ is the group inverse. We shall find the partial solution of (\ref{eq_ex:X'=AX+B}) by (\ref{eq:left_dif_repr_gr}). We have
\[ {\bf A}^{2}=\begin{pmatrix}
  -i & i & 3-i \\
  1 & -1 & 1+3i \\
  -3+i & 3-i & 3+i
\end{pmatrix},\,\, \widehat{{{\rm
{\bf B}}}}^{(1)}={\bf A}{\bf B}= \begin{pmatrix}
 2 -i & 2i & 0 \\
  1+2i & -2 & 0 \\
  1+i & i & 0
\end{pmatrix},\] \[\widehat{{{\rm
{\bf B}}}}^{(2)}={\bf A}^{2}{\bf B}= \begin{pmatrix}
 2 -2i & 2+3i & 0 \\
  2+2i & -3+2i & 0 \\
  1+5i & -2 & 0
\end{pmatrix}.\]
 and
\[\begin{array}{l}
    {{{\sum\limits_{\alpha \in J_{2,\,3}}  {{\left| {\left( {{\rm {\bf A}}^{2} } \right) {\kern 1pt} _{\beta} ^{\beta} }
\right|}}} }}=\\ \det\begin{pmatrix}
  -i & i \\
  1 & -1
\end{pmatrix}+\det\begin{pmatrix}
  -1 & 1+3i \\
  3-i & 3+i
\end{pmatrix}+\det\begin{pmatrix}
  -i & 3-i \\
  -3+i & 3+i
\end{pmatrix}= \\
    0+(-9-9i)+(9-9i)=-18i.
  \end{array}
\]
Since $\left( {{\rm {\bf A}}^{ 2} } \right)_{\,.\,1}
 \left( {{\widehat{{\rm
{\bf b}}}}_{.1}}^{(1)} \right)=
\begin{pmatrix}
  2-i & i & 3-i \\
   1+2i &-1 & 1+3i \\
 1+i&  3-i & 3+i
\end{pmatrix}$ and \[\left( {{\rm {\bf A}}^{ 2} } \right)_{\,.\,1}
 \left( {{\widehat{{\rm
{\bf b}}}}_{.1}}^{(2)} \right)=
\begin{pmatrix}
  2-2i & i & 3-i \\
   2+2i &-1 & 1+3i \\
 1+5i&  3-i & 3+i
\end{pmatrix},\] then finally we obtain
\[\begin{array}{c}
    x_{11}= {\frac{{{\sum\limits_{\beta \in
J_{2,3} {\left\{ {1} \right\}}} {{\left| \left( {{\rm
{\bf A}}^{2}_{\,.\,1} \left( {\widehat{{\rm
{\bf b}}}_{.1}}^{(1)} \right)} \right){\kern 1pt} {\kern 1pt} _{\beta}
^{\beta} \right|}} } }}{{{\sum\limits_{\beta \in J_{2,3}}
{{\left| {\left( {{\rm {\bf A}}^{ 2} } \right){\kern
1pt} {\kern 1pt} _{\beta} ^{\beta} }  \right|}}} }}}+
 \left ({b_{11}- {\frac{{{\sum\limits_{\beta \in
J_{2,3} {\left\{ {1} \right\}}} {{\left| \left( {{\rm
{\bf A}}^{2}_{\,.\,1} \left( {{\widehat{{\rm
{\bf b}}}}_{.1}}^{(2)} \right)} \right){\kern 1pt} {\kern 1pt} _{\beta}
^{\beta} \right|}} } }}{{{\sum\limits_{\beta \in J_{2,3}}
{{\left| {\left( {{\rm {\bf A}}^{ 2} } \right){\kern
1pt} {\kern 1pt} _{\beta} ^{\beta} }  \right|}}} }}}} \right)t=\\
    \frac{3-3i}{-18i}+ \left(1-\frac{-18i}{-18i}\right)t=\frac{1+i}{6}.
  \end{array}
\]
Similarly,
\[x_{12}=\frac{-3+3i}{-18i}+ \left(i-\frac{9+9i}{-18i}\right)t=\frac{-1-i}{6}+\frac{1+i}{2}t,\,\,x_{13}=0+ \left(1-0\right)t=t,\]
\[x_{21}=\frac{3+3i}{-18i}+ \left(i-\frac{-18}{-18i}\right)t=\frac{-1+i}{6},\]
\[x_{22}=\frac{-3-3i}{-18i}+ \left(0-\frac{-9+9i}{-18i}\right)t=\frac{1-i}{6}+\frac{1+i}{2}t,\,\,x_{23}=0+ \left(1-0\right)t=t,\]
\[\]
\[x_{31}=\frac{-12i}{-18i}+ \left(1-\frac{-18i}{-18i}\right)t=\frac{2}{3},\]
\[x_{32}=\frac{9+3i}{-18i}+ \left(i-\frac{-18}{-18i}\right)t=\frac{-1+3i}{6},\,\,x_{33}=0+ \left(0-0\right)t=0.\]
Then
\[{\rm {\bf X}}= \frac{1}{6}\left(
                                    \begin{array}{ccc}
                                      1+i& -1-i+(3+3i)t& t\\
                                    -1+i& 1-i+(3+3i)t& t\\
                                      4 & -1+3i &0\\
                                    \end{array}
                                  \right)
                                   \]
is  the partial solution of (\ref{eq_ex:X'=AX+B}) .

\end{document}